\documentclass[twoside,11pt]{article}

\usepackage{jmlr2e}
\usepackage{epsf}
\usepackage{epsfig}
\usepackage{latexsym}
\usepackage{amssymb}
\usepackage{amsmath}
\usepackage{graphicx}
\usepackage{bbm}
\usepackage{hyperref}
\usepackage{subcaption}
\usepackage{enumerate}
\usepackage{theorem}
\usepackage{multirow}
\usepackage{multicol}
\usepackage{array,booktabs}
\usepackage{algorithm,algpseudocode}

\newcommand{\Gc}{\mathcal{G}}
\newcommand{\Lc}{\mathcal{L}}
\newcommand{\Sc}{\mathcal{S}}

\DeclareMathOperator*{\argmin}{argmin}
\DeclareMathOperator*{\minimize}{minimize}
\DeclareMathOperator{\diag}{diag}
\DeclareMathOperator{\trace}{trace}
\DeclareMathOperator{\sign}{sign}
\algnewcommand{\Inputs}[1]{%
  \State \textbf{Inputs}
  \Statex \hspace*{\algorithmicindent}\parbox[t]{.8\linewidth}{\raggedright #1}
}
\algnewcommand{\Initialize}[1]{%
  \State \textbf{Initialize}\hspace*{\algorithmicindent}\parbox[t]{.8\linewidth}{\raggedright #1}
}
\algnewcommand{\Iterate}[1]{%
  \State \textbf{Iterate}\hspace*{\algorithmicindent}\parbox[t]{.8\linewidth}{\raggedright #1}
}
\newtheorem{thm}{Theorem}

\newtheorem{defn}[thm]{Definition}

\newtheorem{rem}[thm]{Remark}
\newtheorem{lem}[thm]{Lemma}

\begin{document}

\title{
Estimation of Graphical Models through Structured Norm Minimization}

\author{\name Davoud Ataee Tarzanagh \email tarzanagh@ufl.edu \\
       \addr
       Department of Mathematics\\
       UF Informatics Institute \\
       University of Florida\\
       Gainesville, FL 32611-8105, USA
       \AND
       \name George Michailidis \email gmichail@ufl.edu \\
       \addr
       Department of Statistics\\
       UF Informatics Institute \\
       University of Florida\\
       Gainesville, FL 32611-8545, USA}

\editor{Bert Huang}

\maketitle

\begin{abstract}
Estimation of Markov Random Field and covariance models from high-dimensional data represents a canonical problem that has received
a lot of attention in the literature. A key assumption, widely employed, is that of {\em sparsity} of the underlying model. In this paper,
we study the problem of estimating such models exhibiting a more intricate structure comprising simultaneously of {\em sparse,
structured sparse} and {\em dense} components. Such structures naturally arise in several scientific fields, including molecular biology,
finance and political science. We introduce a general framework based on a novel structured norm that enables us to estimate
such complex structures from high-dimensional data. The resulting optimization problem is convex and we introduce a linearized multi-block alternating direction method of multipliers (ADMM) algorithm to solve it efficiently.
We illustrate the superior performance of the proposed framework on a number of synthetic data sets generated from both random and structured networks. Further, we apply the method to a number of real data sets and discuss the results.
\end{abstract}


\section{Introduction}

There is a substantial body of literature on methods for estimating network structures from high-dimensional data, motivated by important biomedical and social science applications; see  \citet{Barab99,Liljeros01,Robins07,Guo10,Danaher13,Friedman07,Tan14,Guo15}. Two powerful formalisms have been employed for this task, the Markov Random Field (MRF) model and the Gaussian covariance graph model (GCGM). The former captures statistical conditional dependence relationships amongst random variables that correspond to the network nodes, while the latter to marginal associations. Since in most applications the number of model parameters to be estimated far exceeds the available sample size, the assumption of sparsity is made and imposed through regularization. An $\ell_1$ penalty on the parameters encoding the network edges is the most common choice; see \citet{Friedman07,ElKaroui08,Cai11,Xue12}, which can also be interpreted from the Bayesian perspective as using an independent double-exponential prior distribution on each edge parameter. Consequently, this approach encourages sparse uniform network structures that may not be the most suitable choice for many real world applications, which in turn have {\em hub} nodes or {\em dense subgraphs}.
As argued in \citet{Barab99,Liljeros01,Newman00,Li05,For10,Newman12} many networks exhibit different
structures at different scales. An example includes a densely connected subgraph, also known as a {\em community} in the social networks literature. Such structures in social interaction networks may correspond to groups of people sharing common interests or being co-located (\citealp{Traud11,Newman04}), while in biological systems to groups of proteins responsible for regulating or synthesizing chemical products (\citealp{Gui05,Lewis10}; see, Figure~\ref{fig2} for an example). Hence, in many applications, simple sparsity or
alternatively, a dense structure fails to capture salient features of the true underlying mechanism that gave rise to the
available data.

In this paper, we introduce a framework based on a novel structured sparse norm that allows us to recover such complex structures.
Specifically, we consider Markov Random Field and covariance models where the parameter of interest, $\Theta$ can be expressed as the superposition of sparse, structured sparse and dense components as follows:
\begin{eqnarray}\label{parts}
  \Theta &=& \underbrace{Z_1 + Z_1^\top}_\text{sparse part} \quad + \quad \underbrace{ Z_2 + Z_2^\top+ \dots+ Z_n + Z_n^\top}_\text{structured sparse part} \quad + \quad \underbrace{E}_\text{dense part},
\end{eqnarray}
where $Z_1$ is a sparse matrix, $Z_2, \dots, Z_n$ are the set of $n-1$ structured sparse matrices (see, Figure~\ref{fig2} for an example of such structured matrices), and $E$ is a dense matrix having possibly very many small, non-zero entries. As shown in Figure~\ref{fig2}, the elements of $Z_1$ represent edges between non-structured nodes, and the non-zero parts of structured matrices $Z_2, \dots, Z_n$ correspond to densely connected subgraphs (communities).

We elaborate more on the decomposition proposed above. We start by discussing on the sparse and structured sparse component and then
elaborate on the dense component. Traditional sparse (lasso \citealp{tib96,Friedman07}) and group sparse (group lasso \citealp{Yuan07,Obo09,Obo11})
are tailor-made to estimate and recover sparse and structured sparse model structures, respectively. However, these methods can not accommodate different structures, unless users specify a {\em priori} the structure of interest (e.g. hub nodes and sparse components), thus severely limiting their application scope. On the other hand, the general framework introduced, is capable of estimating from high-dimensional data, {\em groups with overlaps}, {\em hubs} and {\em dense subgraphs}, with the size and location of such structures {\em not known a priori}.

Next, we discuss the role of the dense component $E$. In many applications, the data generation mechanism may correspond to a true sparse
or structured sparse structure, ''corrupted" by a dense component comprising of possible many small entries. A simple example of such a generating
mechanism in linear models would have the regression coefficient being sparse with a few large entries and a more dense component having possibly
many small, nonzero entries. In such instances, a pure sparse model formulation may not perform particularly well due to the presence of the dense
component and may require very careful tuning to recover the sparse component of interest.  This line of reasoning is also adopted in
\citet{cher15}. Note however, that the model may also be used in settings where there is a significant dense component; however, as
discussed in \citet{cher15} recovery of the individual component is not guaranteed. Hence, in this work we adopt the viewpoint
that $E$ represents a small ''perturbation" of the sparse+structured sparse structure.
To achieve these goals, it leverages a new structured norm that is used as the regularization term of the corresponding objective function.

The resulting optimization problem is solved through a multi-block ADMM algorithm. A key technical innovation is the development of a linearized ADMM
algorithm that avoids introducing auxiliary variables which is a common strategy in the literature. We establish the global convergence of the proposed
algorithm and illustrate its efficiency through numerical experimentation.  The algorithm takes advantage of the special structure of the problem formulation and thus is suitable for large instances of the problem. To the best of our knowledge, this is the first work that gives global convergence guarantees for linearized multi-block ADMM with Gauss-Seidel updates, which is of interest in its own accord.

The remainder of the paper is organized as follows:
In Section~2, we present the new structured norm used as the regularization term in the objective function of the Markov Random Field, covariance graph,
regression and vector auto-regression models. In Section~3, we introduce an efficient multi-block ADMM algorithm to solve the problem, and provide the convergence analysis of the algorithm. In Section~4, we illustrate the proposed framework on a number of synthetic and real data sets,
while some concluding remarks are drawn in Section 5.

\section{A General Framework for Learning under Structured Sparsity}\label{sec2}

We start by introducing key definitions and associated notation.

\subsection{Symmetric Structured Overlap Norm}

Let $X$ be an $m \times p$ data matrix, $\Theta$ be a $p\times p$ symmetric matrix containing the parameters of interest of the statistical
loss function $\Gc(X, \Theta)$. The most popular assumption used in the literature is that $\Theta$ is {\em sparse} and can be successfully
recovered from high-dimensional data by solving the following optimization problem
\begin{equation}\label{loss1}
\minimize_{\Theta \in \Sc} \quad \Gc(X, \Theta) + \lambda \big\| \Theta \big \|_1,
\end{equation}
where $\Sc$ is some set depending on the loss function; $\lambda$ is s a non-negative regularization constant; and $\|.\|_1$ denotes the $\ell_1$ norm or the sum of the absolute values of the matrix elements.

To explicitly model different structures in the parameter $\Theta$, we introduce the following {\em symmetric structured overlap norm (SSON)}:

\begin{defn}\label{def1} (\textbf{Symmetric Structured Overlap Norm}). Let $\Theta$ be a $p\times p$ symmetric matrix containing the model parameters of interest. The symmetric structured overlap norm for a set of partitioned matrices $Z_1, \dots, Z_n$ is given by,
\begin{eqnarray}\label{eq:3}
 \nonumber 
\minimize_{Z_1,\dots, Z_n ,~E} \quad \Omega(\Theta, Z_1, \dots,Z_n,E)&:=&\lambda_1 \|Z_1-\diag(Z_1)\|_1 \\
\nonumber
&+& \sum_{i=2}^{n} \hat{\lambda}_i \|Z_i-\diag(Z_i)\|_1
+ \lambda_i \sum_{j=1}^{l_i}\|(Z_i-\diag(Z_i))_j\|_{F} \\
\nonumber
&+&
\frac{\lambda_e}{2}\|E\|^2_{F},\\
\Theta &=& \sum_{i=1}^{n} \bigl(Z_i +Z_i^\top\bigr)+ E ,
\end{eqnarray}
where $\{\lambda_i\}_{i=1}^{n} $ and $\{\hat{\lambda}_i\}_{i=2}^{n}$ are nonnegative regularization constants; $l_i$ is the number of blocks of the partitioned matrix $Z_i$; $(Z_i-\diag(Z_i))_j$ is the $j$th block of the partitioned matrix $Z_i$; $E$ is an unstructured noise matrix; $\|.\|_1$ denotes the $\ell_1$ norm or the sum of the absolute values of the matrix elements;  and $\|.\|_F$ the Frobenius norm.
\end{defn}

We note that the overlap norm defined by \citet{mohan12,Tan14} encourages the recovery of matrices that can be expressed as a union of few rows and the corresponding columns (i.e. hub nodes). However, SSON represents a new symmetric and significantly more general variant of the overlap norm that promotes matrices that can be expressed as the sum of symmetric structured matrices. Moreover, unlike the previous group sparsity and the latent group lasso discussed in \citet{Yuan07,Obo09,Obo11} that require users to specify structures of interest \textit{a priori}, the SSON achieves a similar objective in an {\em agnostic manner}, relying only on how well such structures fit the observed data.

\begin{figure}
\centering
    \includegraphics[width=.45\textwidth]{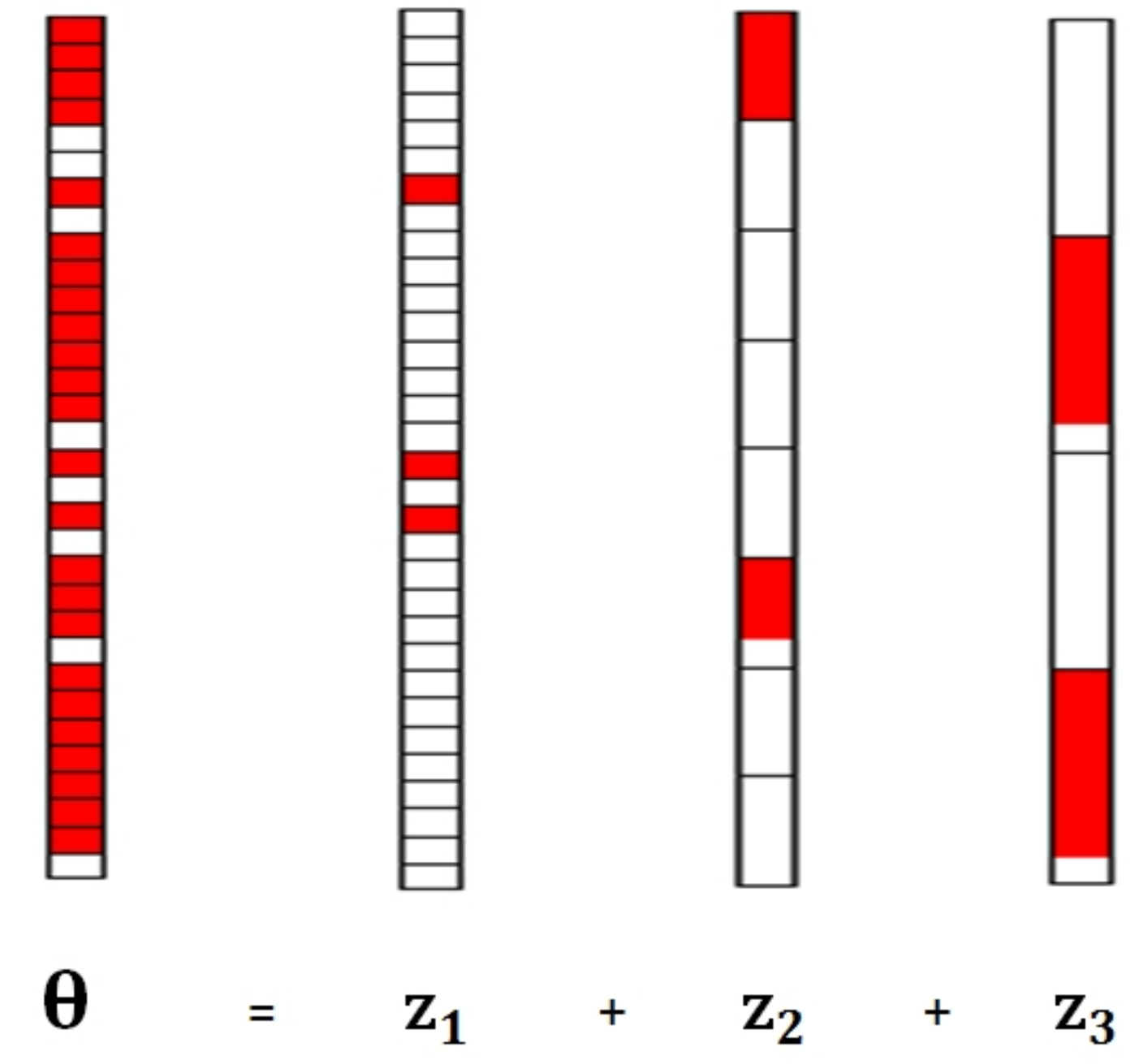}
 \caption{Decomposition of a vector $\theta$ into partitioned vectors $z_1$, $z_2$ and $z_3$,  where $z_1$ is sparse, $z_2$ and $z_3$  are structured sparse vectors. White and red elements are zero and non-zero in the model parameter vector $\theta$, respectively.
}
\label{partvec}
\end{figure}

In many applications, such as regression models, we are interested in modeling different structures in a parameter vector $\theta$. In these cases,  we have the following definition as a special case of SSON:

\begin{defn}\label{def2} Let $\theta$ be a $p\times 1$ vector containing the model parameters of interest. The structured overlap norm for a set of partitioned vectors $z_1, \dots, z_n$ is given by,
\begin{eqnarray}\label{eq:3-3}
 \nonumber 
\minimize_{z_1,\dots, z_n ,~e} \quad \omega(\theta, z_1, \dots,z_n,e)&:=&\lambda_1 \|z_1\|_1 +\sum_{i=2}^{n} \hat{\lambda}_i \|z_i\|_1+ \lambda_i \sum_{j=1}^{l_i}\|{z_i}_j\|_{2} +\frac{\lambda_e}{2}\|e\|^2_{2},\\
\theta &=&  z_1+ z_2 + \dots + z_n + e,
\end{eqnarray}
where $\{\lambda_i\}_{i=1}^{n} $ and $\{\hat{\lambda}_i\}_{i=2}^{n}$ are nonnegative regularization constants; $l_i$ is the number of blocks of the partitioned vector $z_i$; $z_{i_j}$ is the $j$th block of the partitioned vector $z_i$ (see, Figure~\ref{partvec}); $e$ is an unstructured noise vector; $\|.\|_1$ denotes the $\ell_1$ norm or the sum of the absolute values of the vector elements;  and $\|.\|_2$ the two norm.
\end{defn}

\begin{rem}
In the formulation of the problem, $\lambda_1$, $\{\hat{\lambda}_2, \dots, \hat{\lambda}_n, \lambda_2, \dots, \lambda_n\}$, and $\lambda_e$ are tuning parameters corresponding to the sparse component $Z_1$,  the structured components $\{Z_2, \dots Z_n\}$ and the dense (noisy) component $E$, respectively. While the nonzero components may be clustered into groups, the nonzero groups may also be sparse. The latter
can be achieved by \eqref{eq:3} when $\{\hat{\lambda}_2, \dots, \hat{\lambda}_n\}$ are positive constants.
\end{rem}

\begin{rem} The SSON admits the lasso \citep{tib96}, the group lasso with overlaps \citep{Obo09,Obo11} and the ridge shrinkage \citep{hoerl70} methods as three extreme cases, by respectively setting $\{\hat{\lambda}_2, \dots, \hat{\lambda}_n, \lambda_2 , \dots, \lambda_n, \lambda_e\}\rightarrow \infty$,  $\{\lambda_1, \hat{\lambda}_2, \dots, \hat{\lambda}_n, \lambda_e\} \rightarrow \infty$, and $\{\lambda_1, \dots, \lambda_n,  \hat{\lambda}_2, \dots, \hat{\lambda}_n\}\rightarrow \infty$\footnote{For example, with $\lambda_e \rightarrow \infty$, we set $ \frac{\lambda_e}{2}\|E\|^2_{F}=0$  when $E=0$,  so the problem is well-defined.}.
\end{rem}

Note that SSON is rather different from the sparse group lasso, which also uses a combination of $\ell_1$ and $\ell_G$ penalization, where $\|.\|_G$ is the group lasso norm. The sparse group lasso penalty is $\bar{\omega}(\theta)= \lambda_1 \|\theta\|_1+ \lambda_2\|\theta\|_G$, and thus the includes lasso and group lasso as extreme cases corresponding to $\lambda_2 = 0$ and $\lambda_1 = 0$, respectively. However, $\bar{\omega}(\theta)$ does not split $\theta$ into a sparse and a group sparse part and will produce a sparse solution as long as $\lambda_1 > 0$. Hence, the sparse group lasso method can be thought of as a sparsity-based method with additional shrinkage by $\|\theta\|_G$. The group sparsity processes data very differently from SSON and consequently has very different prediction risk behavior. The same argument illustrates the advantages of the proposed SSON penalty over the well-known elastic net penalty. The elastic net is a combination of lasso and ridge penalties \citep{zou05}. However, the elastic net does not split $\theta$ into a sparse and a dense component. Our results show that SSON tends to perform no worse than, and often performs significantly better than ridge, lasso, group lasso or elastic net with penalty levels chosen by cross-validation.

\begin{rem}
  In order to encourage different structures in the parameter matrix $\Theta$, we consider the Frobenius norm of blocks of partitioned matrices, which leads to recovery of dense subgraphs. Other values for the norm of such blocks are also possible; e.g. the $\ell_\infty$ norm.
\end{rem}

\begin{figure}[!ht]
\captionsetup[subfigure]{labelformat=empty}
 \begin{subfigure} {0.45\textwidth}
    {\includegraphics[width=\textwidth]{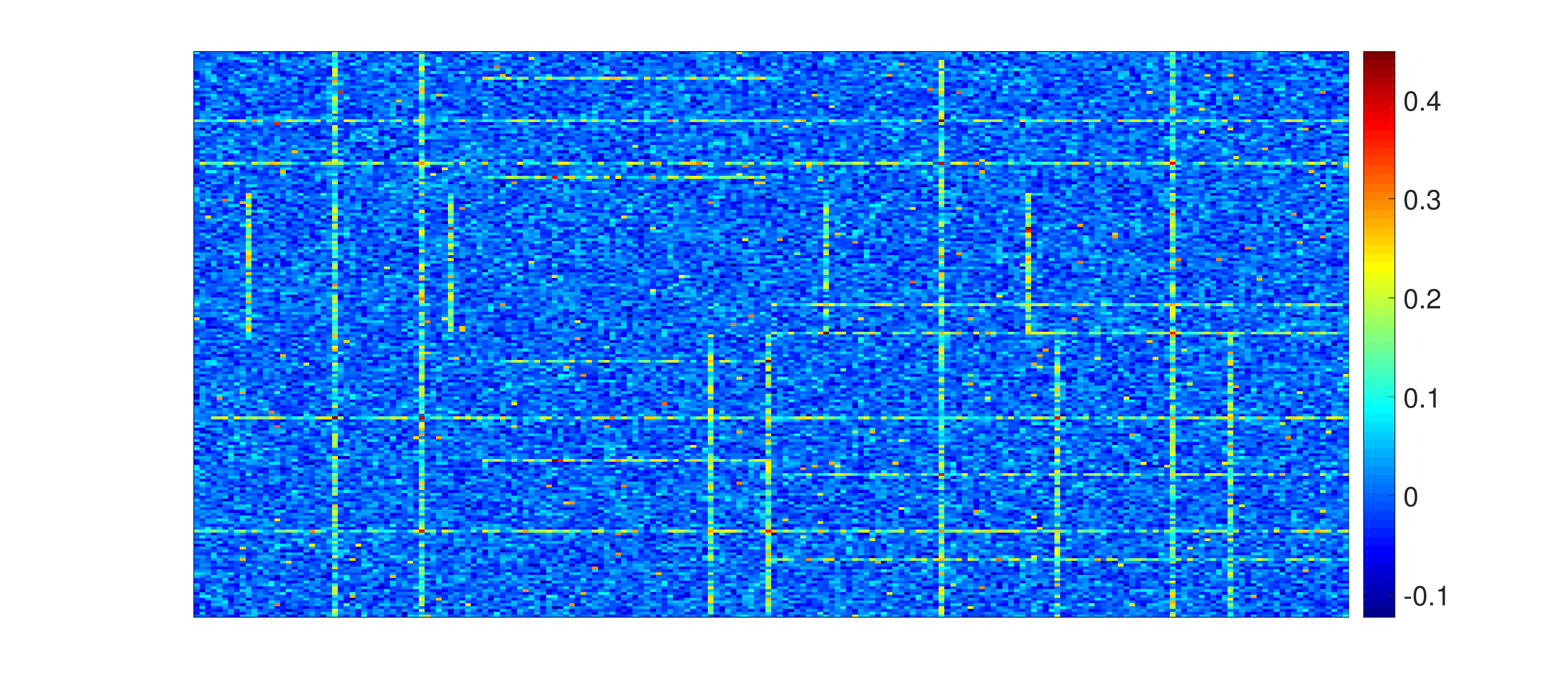}}
    \caption{Ground truth+noise.}
  \end{subfigure}
  \hfill
  \begin{subfigure} {0.45\textwidth}
  {\includegraphics[width=\textwidth]{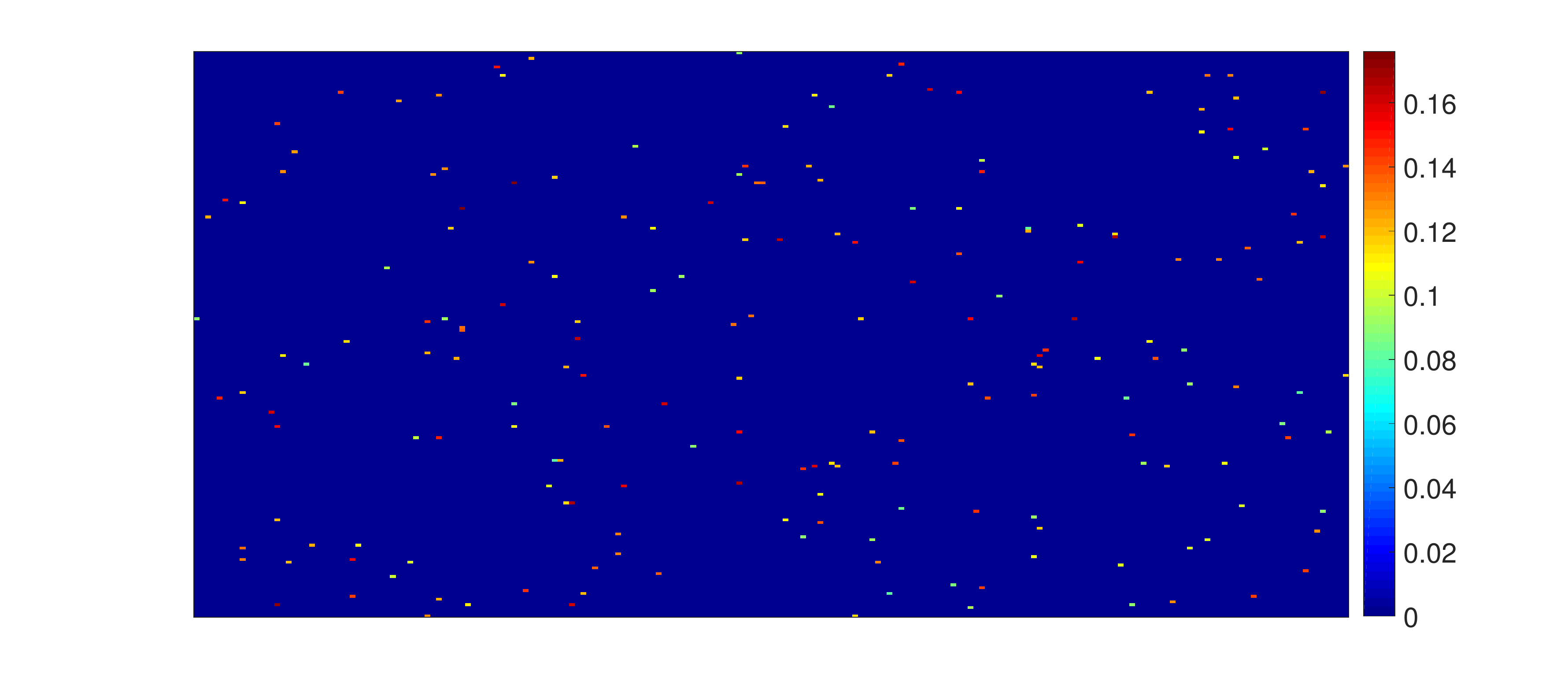}}
    \caption{Sparse part.}
    \end{subfigure}
    \hfill
  \begin{subfigure} {0.45\textwidth}
  {\includegraphics[width=\textwidth]{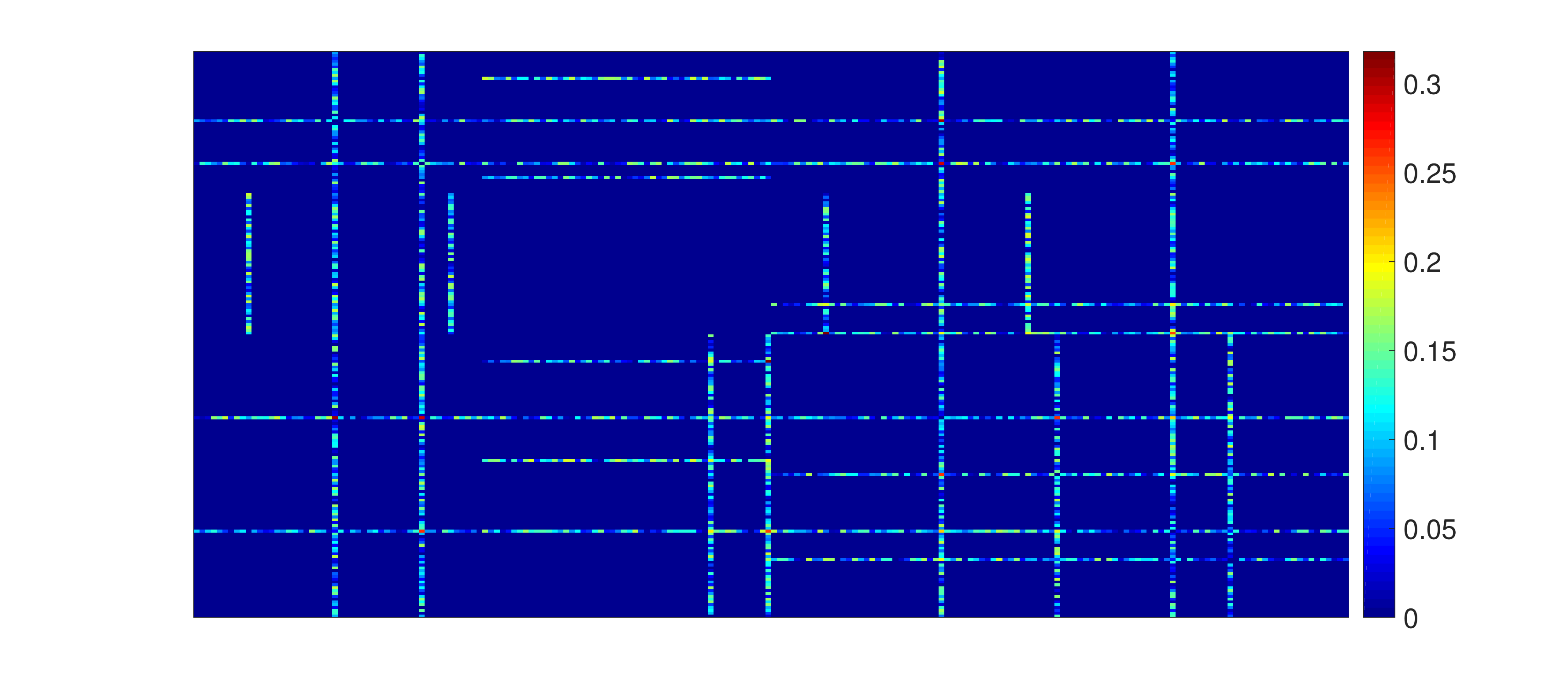}}
    \caption{Structured sparse part.}
  \end{subfigure}
 \hfill
 \begin{subfigure} {0.45\textwidth}
    {\includegraphics[width=\textwidth]{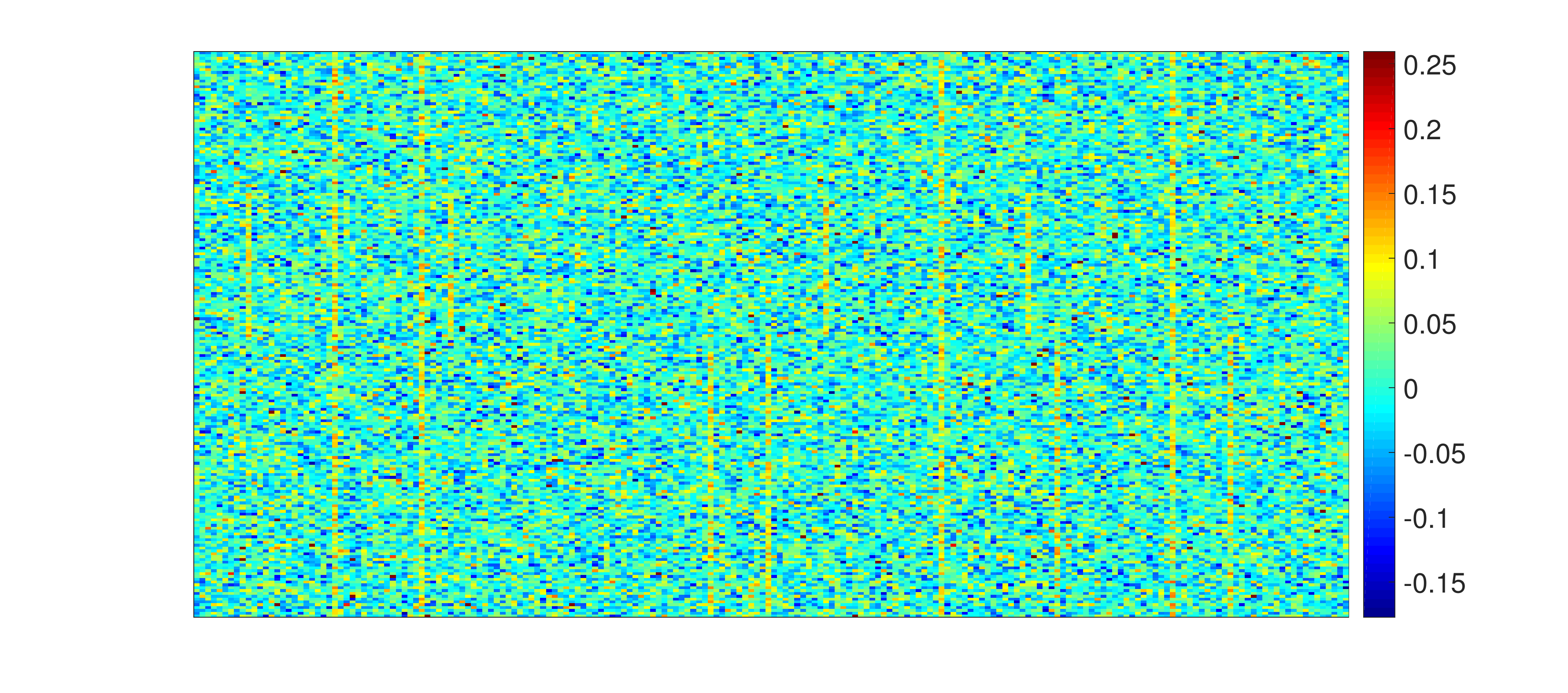}}
    \caption{Noisy part.}
  \end{subfigure}
\caption{Heat map of the covariance matrix $\Theta_3$ decomposed into sparse and structured sparse parts in the presence of noise, estimated by SSON using problem \eqref{eq:7}.}
\label{fignoisdecom}
\end{figure}

\begin{rem}
The matrix $E$ is an important component of the SSON framework.
\begin{description}
\item It enables to develop a {\em convergent multi-block ADMM} to solve the problem of estimating a structured Markov Random Field or covariance model. Note that in general, a direct extension of ADMM to multi-block convex minimization problems is {\em not} necessarily convergent even without linearization of the corresponding subproblems as shown in~\citet{chen16}.
\item From a performance standpoint, our results show that adding a ridge penalty term $\frac{\lambda_e}{2} \|E\|_F^2$  to the structured norm is provably effective in correctly identifying the underlying structures in the presence of noisy data \citep{zou05,cher15} (see, Figure~\ref{fignoisdecom} for an example of decomposition \eqref{parts} in the presense of noise for covariance matrix estimation.)
\end{description}
\end{rem}

Next, we discuss the use of the SSON as a regularizer for maximum likelihood estimation of the following popular statistical models:
(i) members of the Markov Random Field family including the Gaussian graphical model, the Gaussian graphical model with latent variables and the binary Ising model, (ii) the Gaussian covariance graph model and (iii) the classical regression and the vector auto-regression models.
For the sake of completeness, we provide a complete, but succinct description of the corresponding models and the proposed regularization.
\begin{figure}
\begin{minipage}[c][6.5cm][t]{.35\textwidth}
  \vspace*{\fill}
  \centering
  \includegraphics[width=\textwidth]{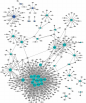}
  \subcaption{\tiny{Postpartum NAC Gene Network (\citealp{zhao14}).}}
  \label{fig:test1}
\end{minipage}%
  \hfill
\begin{minipage}[c][6.5cm][t]{.45\textwidth}
  \centering
\begin{subfigure}[b]{0.15\textwidth}
    \includegraphics[width=\textwidth]{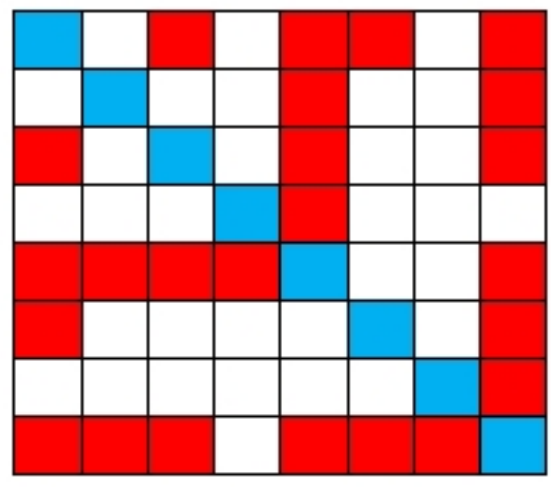}
    \subcaption*{$\Theta$}
  \end{subfigure}
  \vfill
\begin{subfigure}[b]{0.15\textwidth}
    \includegraphics[width=\textwidth]{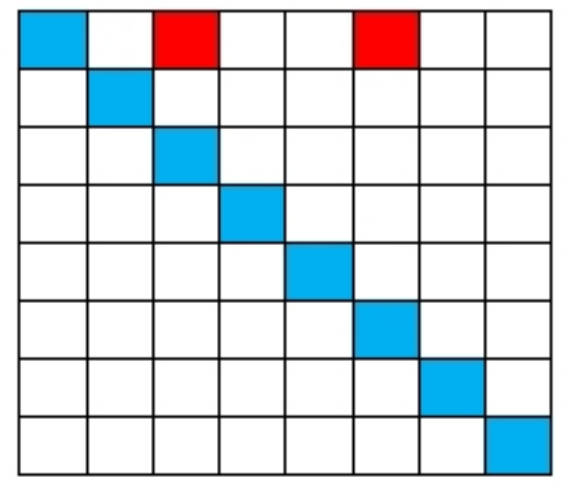}
    \subcaption*{$Z_1$}
  \end{subfigure}
\hspace{2ex}
\begin{subfigure}[b]{0.15\textwidth}
    \includegraphics[width=\textwidth]{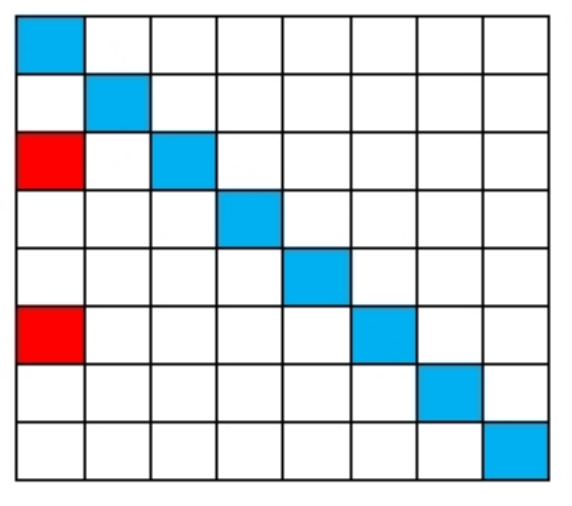}
    \subcaption*{$Z^T_1$}
  \end{subfigure}
  \vfill
\begin{subfigure}[b]{0.15\textwidth}
    \includegraphics[width=\textwidth]{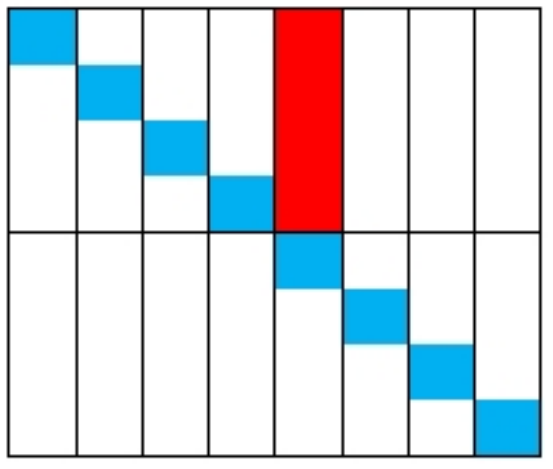}
    \subcaption*{$Z_2$}
  \end{subfigure}
\hspace{2ex}
\begin{subfigure}[b]{0.15\textwidth}
    \includegraphics[width=\textwidth]{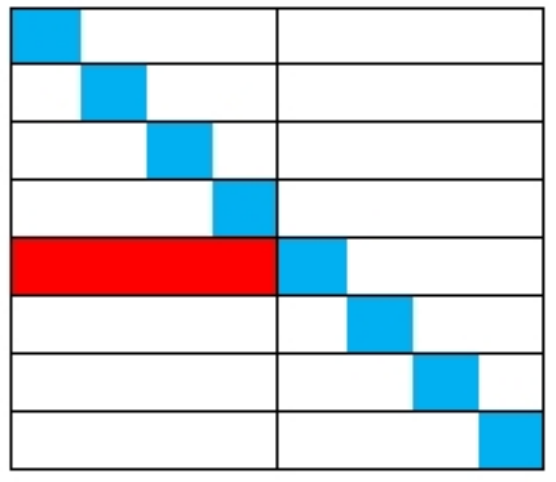}
    \subcaption*{$Z^T_2$}
  \end{subfigure}
  \vfill
\begin{subfigure}[b]{0.15\textwidth}
    \includegraphics[width=\textwidth]{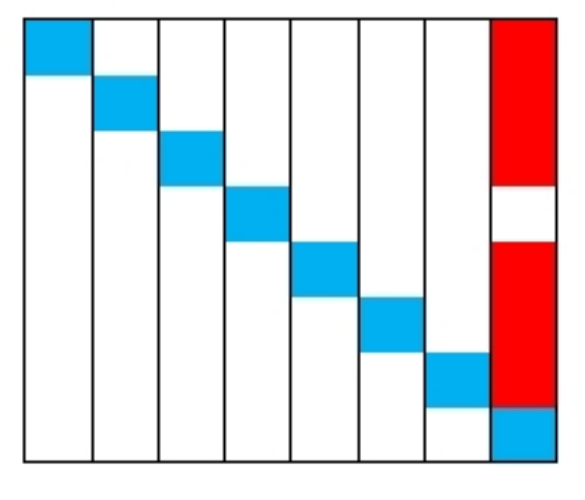}
    \subcaption*{$Z_3$}
  \end{subfigure}
\hspace{2ex}
  \begin{subfigure}[b]{0.15\textwidth}
    \includegraphics[width=\textwidth]{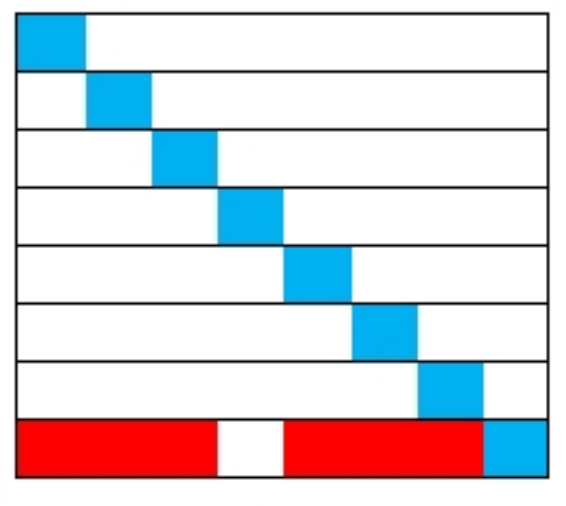}
    \subcaption*{$Z^T_3$}
  \end{subfigure}
 \subcaption{\tiny{Examples of partitioned matrices for the underlying network in (a).}}
\end{minipage}\label{fffg}
\vfill
\vspace{6ex}
\caption{
The figure illustrates that block partitions through structured matrices could be set based on a desire for interpretability of the resulting estimated network structure. Panel (a) shows example of structured gene network, while panel (b) provides decomposition into structured matrices for the network in (a).  Blue elements are diagonal ones, white elements are zero and red elements are non-zero in the model parameter matrix $\Theta$.  The structured penalty function \eqref{eq:3} is then applied to each block for matrices $\{Z_i\}_{i=1}^{n}$.}
\label{fig2}
\end{figure}

\subsection{Structured Gaussian Graphical Models}

Let $X$ be a data matrix consisting of $p$-dimensional samples from a zero mean Gaussian distribution,
$$
x_1, \ldots, x_m \stackrel{\text{i.i.d.}}{\sim} \mathcal{N}(0, \Sigma).
$$
In order to obtain a sparse and interpretable estimate of the precision matrix $\Sigma^{-1}$ that captures conditional dependence relationships, many authors have considered the well-known \textit{graphical lasso} problem \citep{Friedman07,Rothman08} in the form of \eqref{loss1} with loss function
\begin{equation}\label{lgaus}
  \Gc_1(X, \Theta_1):= \trace (\hat{\Sigma}\Theta_1)- \log \det \Theta_1, \qquad \Theta_1 \in \Sc,
\end{equation}
where $\hat{\Sigma}$ is the empirical covariance matrix of $X$; $\Theta_1$ is the estimate of the precision matrix $\Sigma^{-1}$; and $\Sc$ is the set of $p \times p$ symmetric positive definite matrices.

As is well known, the norm penalty in \eqref{loss1} encourages zeros (sparsity) in the solution. However, as previously argued, many biological and social network applications exhibit more complex structures than mere sparsity. Using the proposed SSON, we define the following objective function for the problem at hand:
\begin{eqnarray} \label{eq:5}
\nonumber
\minimize_{\Theta_1, Z_1,\dots, Z_n \in \Sc,~E} && \Gc_1(X, \Theta_1)+ \Omega(\Theta_1, Z_1, \dots,Z_n,E),\\
\Theta_1  &=& \sum_{i=1}^{n} \bigl(Z_i + Z_i^\top\bigr) + E,
\end{eqnarray}
where $\Theta_1$ is the model parameter matrix and $\Omega(\Theta_1, Z_1, \dots, Z_n, E)$ the corresponding SSON defined in \eqref{eq:3}.

\begin{figure}[!htb]
 \begin{subfigure}[b]{0.29\textwidth}
    {\includegraphics[width=\textwidth]{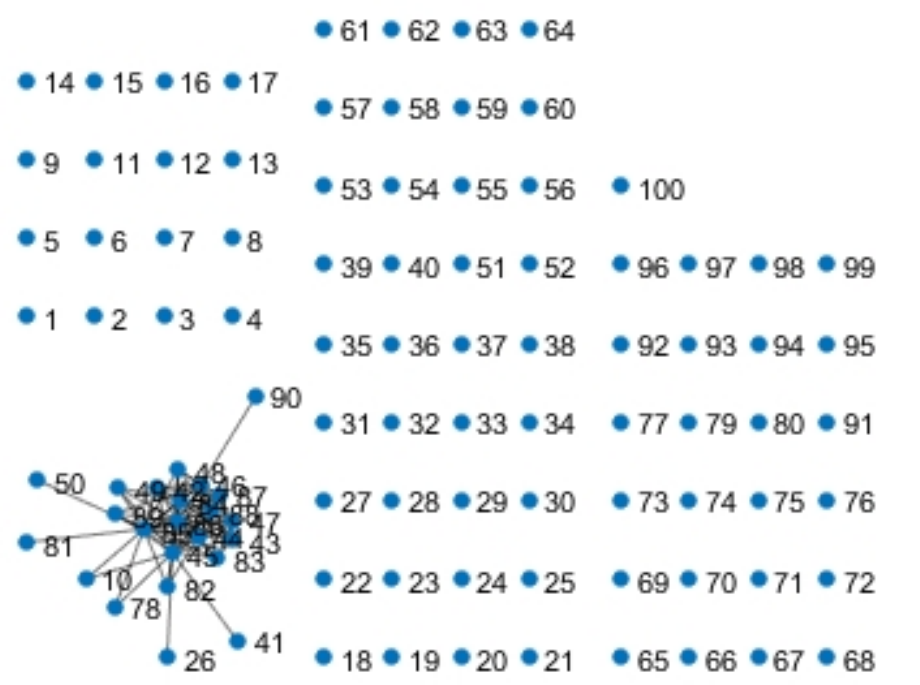}}
    \caption{Graphical lasso regularization. }
    \label{fig3d}
  \end{subfigure}
  \hfill
  \begin{subfigure}[b]{0.29\textwidth}
  {\includegraphics[width=\textwidth]{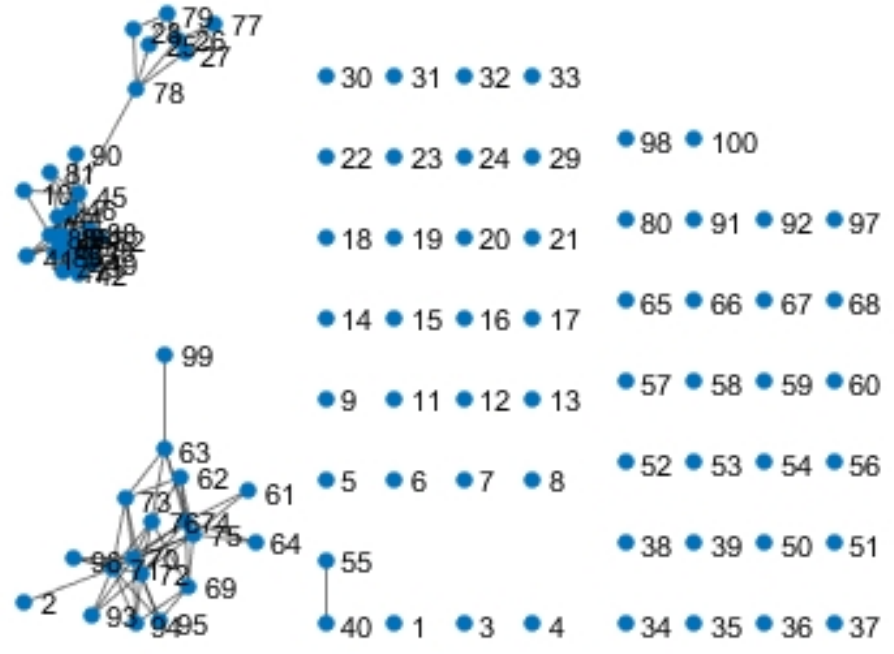}}
    \caption{SSON based regularization.}
    \label{fig3e}
  \end{subfigure}
  \hfill
  \begin{subfigure}[b]{0.29\textwidth}
  {\includegraphics[width=\textwidth]{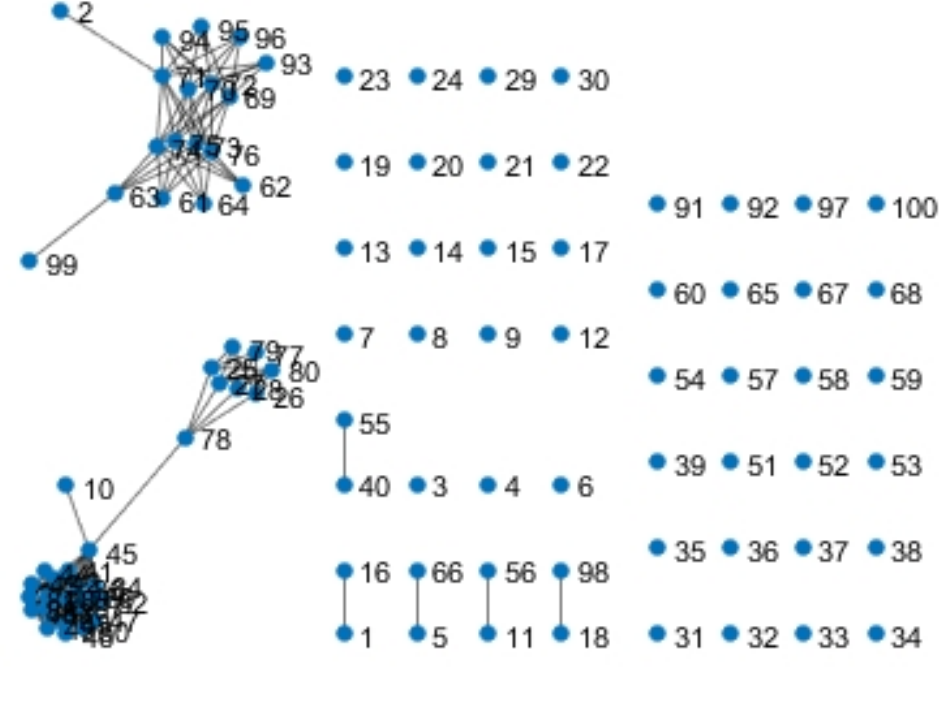}}
    \caption{Ground truth.}
    \label{fig3f}
  \end{subfigure}

 \begin{subfigure}[b]{0.31\textwidth}
    {\includegraphics[width=\textwidth]{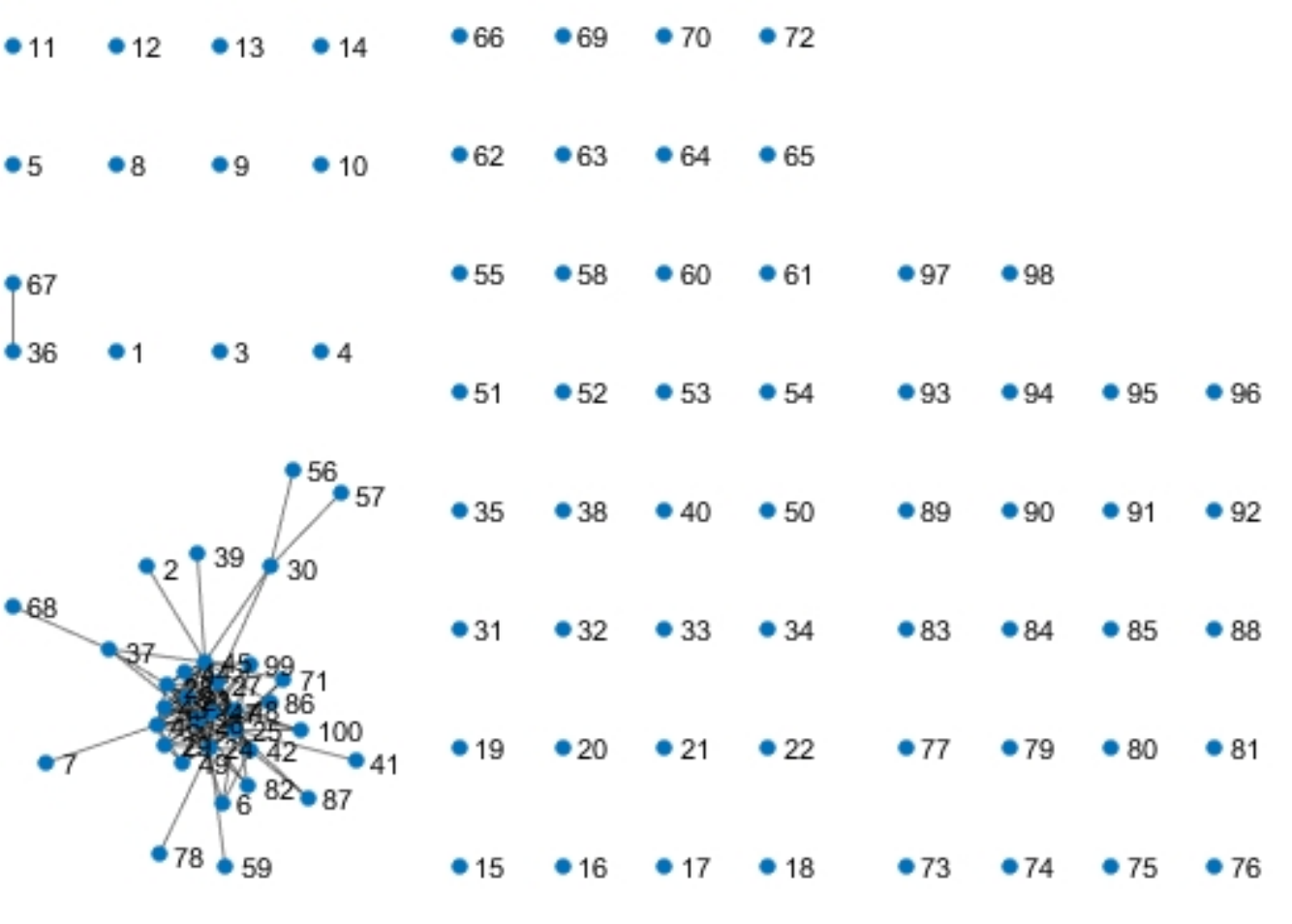}}
    \caption{Graphical lasso regularization.}
    \label{fig3g}
  \end{subfigure}
  \hfill
  \begin{subfigure}[b]{0.29\textwidth}
  {\includegraphics[width=\textwidth]{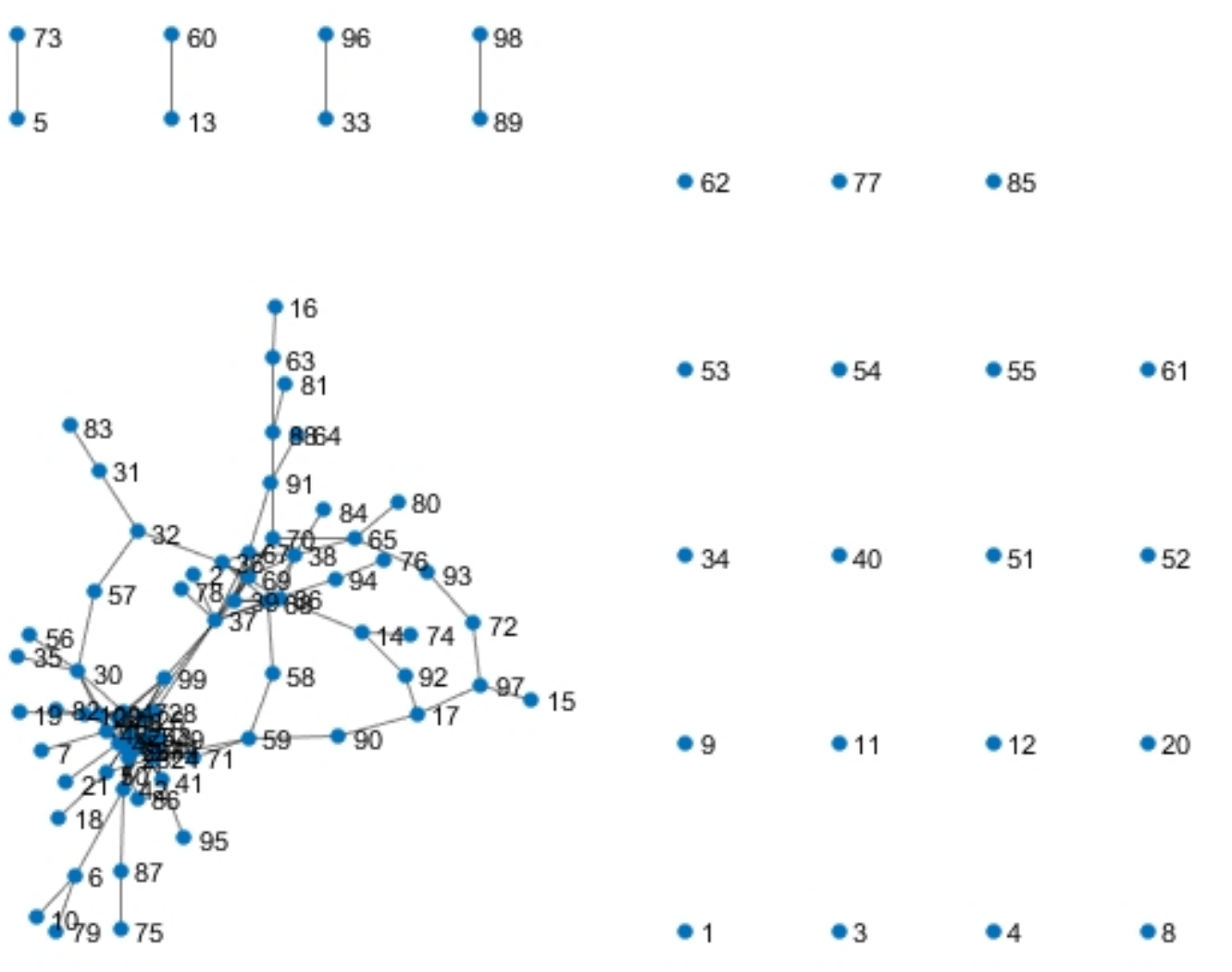}}
    \caption{SSON based regularization.}
    \label{fig3h}
  \end{subfigure}
  \hfill
  \begin{subfigure}[b]{0.29\textwidth}
  {\includegraphics[width=\textwidth]{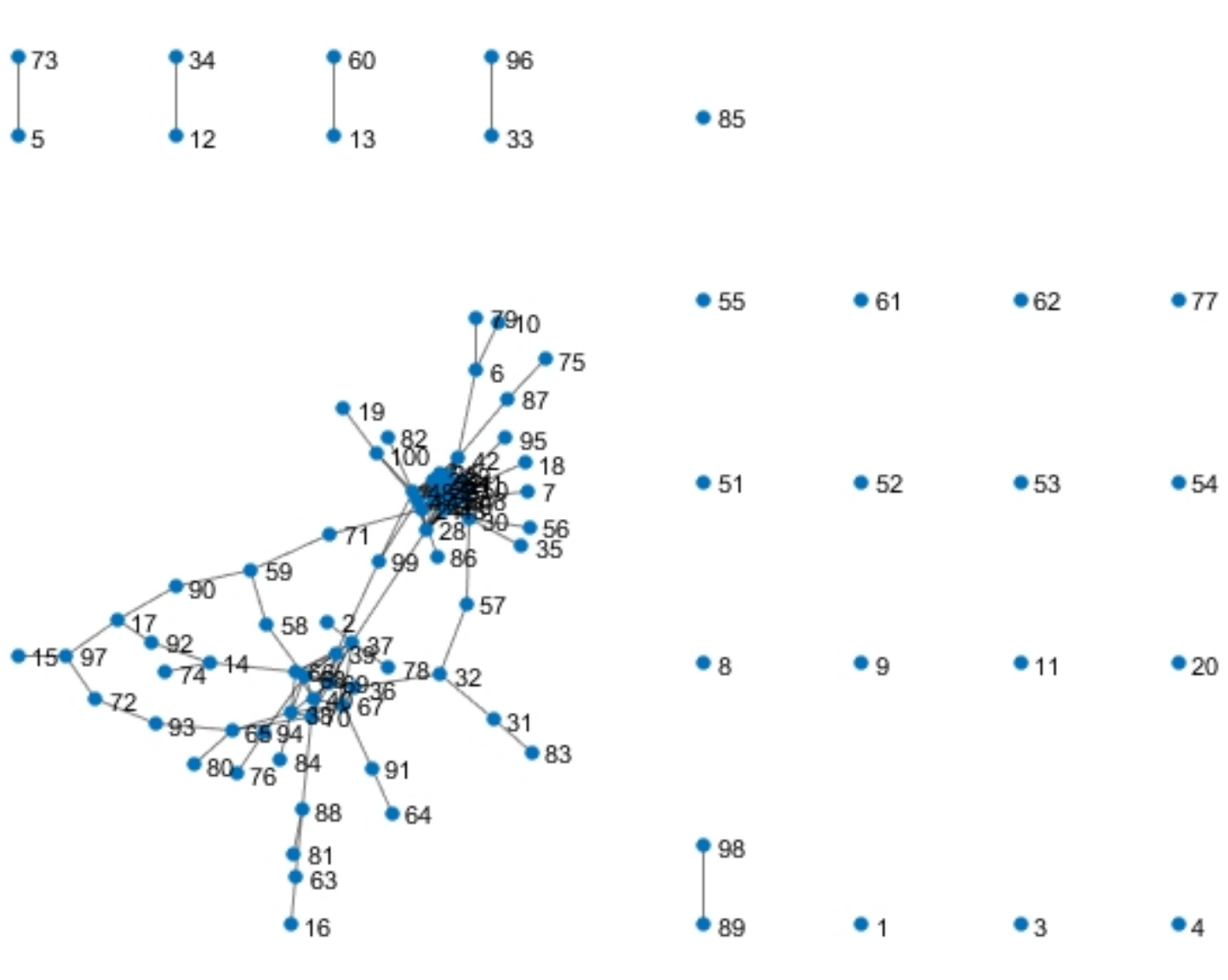}}
    \caption{Ground truth.}
    \label{fig3i}
  \end{subfigure}
\caption{Estimates from the SSON based regularization on two examples of Gaussian graphical models comprising of $p=100$ nodes,
using in \eqref{fig3e} three structured matrices and in \eqref{fig3h} four structured matrices.}
\label{fig3}
\end{figure}

Formulation \eqref{eq:5} allows us to obtain more accurate and compact network estimates than conventional methods whenever the network exhibits different structures. Moreover, our formulation does not require \textit{a priori} knowledge of the underlying network structure
(i.e. which nodes in the network form densely connected subgraphs (see, Figure~\ref{fig3})).

In Figure \ref{fig3}, the performance of our proposed approach is illustrated on two simulated data sets exhibiting different structures
 (sub-figures~\eqref{fig3f} and \eqref{fig3i}); it can be seen that the proposed SSON based graphical lasso (sub-figures~\eqref{fig3e} and \eqref{fig3h}) can recover the network structure much better than the popular graphical lasso based estimator \citep{Friedman07} (subfigures~\eqref{fig3d} and \eqref{fig3g}).

\subsection{Structured Ising Model}

Another popular graphical model, suitable for binary or categorical data, is the Ising one \citep{Ising25}. It is assumed that
observations $x_1, \dots, x_m$ are independent and identically distributed from
\begin{equation}\label{eq:11}
  f(x,\Theta_2)= \frac{1}{\mathbb{W}(\Theta_2)} \exp \Big( \sum_{j=1}^{p} \theta_{jj} x_j + \sum_{1\leq j< j'\leq p} \theta _{jj'} x_j x_{j'}\Big),
\end{equation}
where $\mathbb{W}(\Theta_2)$ is the partition function, which ensures that the density sums to one. Here, $\Theta_2$ is a $ p\times p$ symmetric matrix that specifies the network structure: $\theta_{jj'}=0$ implies that the $j$th and $j'$th variables are conditionally independent given the remaining ones.

Several papers proposing estimation procedures for this model have been published. \citet{lee06} considered maximizing an $\ell_1$-penalized log-likelihood for this model. Due to the difficulty in computing the log-likelihood with the expensive partition function, several authors have considered alternative approaches. For instance, \citet{Ravi11} proposed a neighborhood selection approach. The latter proposal involves solving $p$ logistic regressions separately (one for each node in the network), which leads to an estimated parameter matrix that is in general not symmetric. In contrast, several authors considered maximizing an $\ell_1$-penalized pseudo-likelihood with a symmetric constraint on $\Theta_2$ \citep{Guo10,Guo11}. Under the model \eqref{eq:11}, the $\log$-pseudo-likelihood for $m$ observations takes the form
\begin{equation}\label{eq:12}
   \Gc_2(X, \Theta_2):= \sum_{j=1}^{p}\sum_{j'=1}^{p} \theta_{jj'} (X^T X)_{jj'}-\sum_{i=1}^{m}\sum_{j=1}^{p} \log \Big(1 + \exp \big( \theta_{jj}+ \sum_{j'\neq j}\theta_{jj'}x_{ij'}\big)\Big),
\end{equation}

We propose instead to impose the SSON on $\Theta_2$ in \eqref{eq:12} in order to estimate a binary network with different structures.
This leads to the following optimization problem
\begin{eqnarray} \label{eq:13}
\nonumber
\minimize_{\Theta_2, Z_1,\dots, Z_n \in \Sc,~E}&&  \Gc_2(X, \Theta_2)
+ \Omega(\Theta_2, Z_1, \dots, Z_n, E),\\
\Theta_2 &=& \sum_{i=1}^{n} \bigl(Z_i + Z_i^\top\bigr) + E,
\end{eqnarray}
where $\Theta_2$ is the model parameter matrix and $\Omega(\Theta_2, Z_1, \dots, Z_n, E)$ the corresponding SSON defined in \eqref{eq:3}.

An interesting connection can be drawn between our technique and the Ising block model discussed in \citet{bert16}, which is
a perturbation of the mean field approximation of the Ising model known as the Curie-Weiss model: the sites are partitioned into two blocks of equal size and the interaction between those within the same block is stronger than across blocks, to account for more order within each
block. However, one can easily seen that the Ising block model is a special case of \eqref{eq:13}.

\subsection{Structured Gaussian Covariance Graphical Models}

Next, we consider estimation of a covariance matrix under the assumption that $$x_1, \ldots, x_m \stackrel{\text{i.i.d.}}{\sim} \mathcal{N}(0, \Sigma).$$
This is of interest because the sparsity pattern of $\Sigma$ specifies the structure of the marginal independence graph \citep{Drton03,Drton08}.

Let $\Theta_3$ be a $p \times p$ symmetric matrix containing the parameters of interest. Setting the loss function $\Gc_3(X, \Theta_3) := \dfrac{1}{2}\| \Theta_3 - \hat{\Sigma}\|^2_{F}$, \citet{Xue12} proposed to estimate the positive definite covariance matrix, $\Sigma$ by solving
\begin{equation} \label{eq:6}
\minimize_{\Theta_3 \in \Sc}  \quad \Gc_3(X, \Theta_3) + \lambda \| \Theta_3 \|_1,
\end{equation}
where $\hat{\Sigma}$ is the empirical covariance matrix, $\Sc = \{ \Theta_3 : \Theta_3\succeq \varepsilon I \text{ and } \Theta_3 = \Theta^T_3\}$, and $\varepsilon$ is a small positive constant. We extend \eqref{eq:6} to accommodate structures of the covariance graph by imposing the SSON on $\Theta_3$. This results in the following optimization problem
\begin{eqnarray} \label{eq:7}
\nonumber
\minimize_{\Theta_3, Z_1,\dots, Z_n \in \Sc,~E} && \Gc_3(X, \Theta_3)+\Omega(\Theta_3, Z_1, \dots, Z_n, E),\\
\Theta_3 &=& \sum_{i=1}^{n} \bigl(Z_i + Z_i^\top\bigr) + E.
\end{eqnarray}
where $\Theta_3$ is the model parameter matrix and $\Omega(\Theta_3, Z_1, \dots, Z_n, E)$ the corresponding SSON defined in \eqref{eq:3}.

\subsection{Structured Gaussian Graphical Models with latent variables}

In many applications throughout science and engineering, it is often the case that some relevant variables are not observed. For the Gaussian Graphical model, \citet{Chandrasekaran} proposed a convex optimization problem to estimate it in the presence of latent variables. Let $\Theta_4$ be a $p \times p$ symmetric matrix containing the parameters of interest. Setting $\Gc_4(X, \Theta_4):= \langle \Theta_4, \Sigma_O \rangle - \log \det \Theta_4$, their objective function is given by
\begin{eqnarray}\label{eq:8}
\minimize_{\Theta_4, Z_1, Z_2 \in \Sc} &&
\Gc_4(X, \Theta_4) + \alpha \|Z_1\|_1 + \beta \trace(Z_{n+1}) + \mathbbm{1}_{Z_{n+1} \succeq 0}, \nonumber \\
\Theta_4 &=& Z_1 - Z_{n+1},
\end{eqnarray}
where $\Sigma _O$ is the sample covariance matrix of the observed variables; $\alpha$ and $\beta$ are positive constants; and the indicator function
$ \mathbbm{1}_{Z_{n+1} \succeq 0}$
is defined as
\begin{eqnarray*} \label{eq:9}
\mathbbm{1}_{Z_{n+1} \succeq 0}: = \begin{cases}
0 , & \text{if} \; Z_{n+1} \succeq 0, \\
+ \infty, & \text{otherwise.}
\end{cases}
\end{eqnarray*}
This convex optimization problem aims to estimate an inverse covariance matrix that can be decomposed into a sparse matrix $Z_1$ minus a low-rank matrix $Z_{n+1}$ based on high-dimensional data.

Next, we extend the SSON to solve the latent variable graphical model selection. Problem \eqref{eq:8} can be rewritten in the following equivalent form by introducing new variables $\{Z_i\}_{i=1}^{n}$:
\begin{eqnarray} \label{eq:10}
\nonumber
\minimize_{\Theta_4, Z_1,\dots, Z_n \in \Sc,~E} && \Gc_4(X, \Theta_4) + \Omega(\Theta_4, Z_1, \dots, Z_n , E)+\lambda_{n+1} \trace(Z_{n+1}) + \mathbbm{1}_{Z_{n+1} \succeq 0},\\
\Theta_4 &=& \sum_{i=1}^{n} \bigl(Z_i + Z_i^\top\bigr) -Z_{n+1}+ E,
\end{eqnarray}
where $\Theta_4$ is the model parameter matrix and $\Omega(\Theta_4, Z_1, \dots, Z_n, E)$ the corresponding SSON defined in \eqref{eq:3}.

\subsection{Structured Linear Regression and Vector Auto-Regression}

The proposed SSON is also applicable to structured regression problems. Although this is not the main focus on this paper, nevertheless, we include a brief discussion, especially for lag selection in vector autoregressive models that are of prime interest in the analysis of high-dimensional time series data.
The canonical formulation of the regularized regression problem is given by:
\begin{equation}\label{reg1}
  \min_{\theta \in \mathbb{R}^{p}} \|y -X \theta\|_2 + \lambda \Psi(\theta).
\end{equation}
where $\{(y_i,x_i)\}_ {i=1}^{m}$, $y_i \in \mathbb{R}$, $x_i \in \mathbb{R}^ p$, with  $y =[y_1, \dots, y_m]^\top$ being the
response variable and $X = [x_1^\top, \dots x_m^\top]$  a set of $p$-predictors that are assumed to be independently and identically distributed (i.i.d.);
 $\lambda>0$ is a regularization parameter and $\Psi(\theta)$ is a suitable norm. Specific choices of $\Psi(.)$ lead to popular regularizers including the
lasso -$\Psi(\theta)= \|\theta\|_1$- and the group lasso.

We propose instead to impose the SSON on $\theta$ in \eqref{reg1} in order to solve structured regression problems. Problem \eqref{reg1} can be rewritten in the following form by introducing new variables $\{z_i\}_{i=1}^{n}$ and $e$:
\begin{eqnarray} \label{reg2}
\nonumber
\minimize_{\theta, z_1,\dots, z_n, e} && \Gc(X, \theta) + \omega(\theta, z_1, \dots, z_n , e),\\
\theta_4 &=& z_1 + z_2 + \dots + z_{n}+ e,
\end{eqnarray}
where  $\Gc(X, \theta)=\|y -X \theta\|_2$; $\theta$ is the model parameter vector and $\omega(\theta, z_1, \dots, z_n, e)$ the corresponding structured norm defined in \eqref{eq:3-3}.

Problem \eqref{reg2} can equivalently be thought of as a generalization of subspace clustering \citep{elham09}. Indeed, in order to segment the data into their respective subspaces, we need to compute an affinity vector $\theta$ that encodes the pairwise affinities between data vectors.

An interesting application of the SSON for multivariate regression problems is on structured estimation of  \emph{vector autoregression} (VAR) models (\citealp{lut05}), a popular model for economic and financial time series data \citep{tsay05}, dynamical systems \citep{ljung99} and more recently
 brain function connectivity \citep{vald57}. The model captures both temporal and cross-dependencies between stationary time series.
Formally, let $\{ x_1, \dots, x_m \}$ be a $p$-dimensional time series set of observations that evolve over time according to a lag-$d$ model:
\begin{eqnarray*}
x_{t+1}= \sum_{k=1}^d \Theta_k^\top x_{t-k}+ \epsilon_t,
\qquad \epsilon_1, \ldots, \epsilon_{m-1} \stackrel{\text{i.i.d.}}{\sim} \mathcal{N}(0, \Sigma), \qquad t=1,\dots, m-1,
\end{eqnarray*}
where  $\{\Theta\}_{k=1}^d \in \mathbb{R}^{p \times p}$ are \emph{transition
matrices} for different lags, and $\{\epsilon_1, \ldots, \epsilon_{m-1}\}$  independent multivariate Gaussian \emph{white noise} processes. The VAR process is assumed to be stable and stationary (bounded spectral density), while the noise covariance matrix $\Sigma$ is assumed to be positive definite with bounded largest eigenvalue \citep{basu15}.

Given $m$ observations $\{ x_1, x_2, \dots, x_m \}$ from a stationary VAR process, the lag-$m$ VAR can be written
is given by
\begin{equation}\label{auto1}
\underbrace{
\begin{bmatrix}
    x_{m} \\
    x_{m-1} \\
    \vdots \\
    x_{2}
\end{bmatrix}
}_{\textbf{Y}}
=
\underbrace{\begin{bmatrix}
    {x_{m-1}}^\top \\
    {x_{m-2}}^\top \\
    \vdots \\
    {x_{1}}^\top
\end{bmatrix}}_{\textbf{X}} \Theta
+
\underbrace{
\begin{bmatrix}
    {\epsilon_{m-1}}^\top\\
    {\epsilon_{m-2}}^\top \\
    \vdots \\
    {\epsilon_{1}}^\top
\end{bmatrix}
}_{\textbf{$\varepsilon$}}.
\end{equation}
It can be seen that to estimate $\Theta$ one can solve the following least squares problem
\begin{equation}\label{autolasso}
  \min_{\Theta \in \mathbb{R}^{p\times p}} \|Y -X \Theta\|_F.
\end{equation}
However, as the number of component time series increases, the number of parameters to be estimated grows as $dp^2$ ; hence,
structural assumptions are imposed to estimate them from limited sample size. A popular choice is the lasso
 \citep{basu15}, that leads to sparse estimates. However, it does not incorporate the notion of lag selection, which could lead to
certain spurious coefficients coming from further lags in the past. To address this problem, \citet{basu2015jmlr} proposed a thresholded lasso
estimate. However, our SSON can be used for lag selection, that guarantees that more recent lags are favored over further in the past ones.

Let $\Theta_5$ be a $mp \times mp$ symmetric matrix containing the parameters of interest for all $m$ lages of the problem. Setting the loss function $\Gc_5(X, \Theta_5):= \|Y -X \Theta_5\|$, we propose to estimate the transition matrix, $\Theta$ by solving the following optimization problem:
\begin{eqnarray} \label{eqaut:7}
\nonumber
  \min_{\Theta_5, Z_1,\dots, Z_n, E \in \mathbb{R}^{p \times p}} && \Gc_5(X, \Theta_5)+\Omega(\Theta_5, Z_1, \dots, Z_n, E),\\
\Theta_5 &=& \sum_{i=1}^{n} \bigl(Z_i + Z_i^\top\bigr) + E,
\end{eqnarray}
where $\Theta_5$ is the estimate of the covariance matrix and $\Omega(\Theta_5, Z_1, \dots, Z_n, E)$ the corresponding SSON defined in \eqref{eq:3}.

\section{Multi-Block ADMM for Estimating Structured Network Models}\label{sec3}

Objective functions \eqref{eq:5}, \eqref{eq:13}, \eqref{eq:7}, \eqref{eq:10}, \eqref{reg2}, and \eqref{eqaut:7} involve separable convex functions, while the constraint is simply linear, and therefore they are suitable for ADMM based algorithms. We next introduce a linearized multi-block ADMM algorithm  to solve these problems and  establish its global convergence properties.

The alternating direction method of multipliers (ADMM) is widely used in solving structured convex optimization problems due to its superior performance in practice; see  \citet{Scheinberg10,Boyd11,Hong12,Lin15,Lint16,Sun15,Davis15,Haj15,Haj16}. On the theoretical side, \citet{chen16} provided a counterexample showing that the ADMM may {\em fail to converge} when the number of blocks exceeds {\em two}. Hence, many authors reformulate the problem of estimating a Markov Random Field model to a two block ADMM algorithm by grouping the variables and introducing auxiliary variables \citep{ma13,mohan12,Tan14}. However, in the context of large-scale optimization problems, the grouping ADMM method becomes expensive due to its high memory requirements. Moreover, despite lack of convergence guarantees under standard convexity assumptions, it has been observed by many researchers that the unmodified multi-block ADMMs with Gauss-Seidel updates often outperform all its modified versions in practice \citep{Wang2013,Sun15,Davis15}.

Next, we present a {\em convergent multi-block} ADMM with Gauss-Seidel updates to solve convex problems
\eqref{eq:5}, \eqref{eq:13}, \eqref{eq:7}, \eqref{eq:10}, and \eqref{eqaut:7}. The ADMM is constructed for an augmented Lagrangian function defined by
\begin{eqnarray}\label{eq:3.1}
   \Lc_{\gamma} (\Theta, Z_1, \dots, Z_n,E; \Lambda) &=& \Gc(X,\Theta)+ f_1(Z_1)+\dots +f_n(Z_n)+ \frac{\lambda_e}{2}\|E\|^2_{F}\\
   \nonumber 
 &-& \langle \Lambda, \Theta-\sum_{i=1}^{n}Z_i+Z_i^\top -E \rangle + \frac{\gamma}{2} \| \Theta-\sum_{i=1}^{n}Z_i+Z_i^\top -E \|^2_{F},
\end{eqnarray}
where $\Lambda$ is the Lagrange multiplier, $\gamma$ a penalty parameter, $\Gc(X,\Theta)$ the loss function of interest and
\begin{eqnarray} \label{eq:3.2}
 \nonumber 
f_1(Z_1) &:=& \lambda_1 \|Z_1-\diag(Z_1)\|_1,\\
f_i(Z_i) &:=& \hat{\lambda}_i \|Z_i-\diag(Z_i)\|_1+\lambda_i \sum_{j=1}^{l_i}\|(Z_i-\diag(Z_i))_j\|_{F}, \quad i=2, \dots, n.
\end{eqnarray}

In a typical iteration of the ADMM for solving \eqref{eq:3.1}, the following updates are implemented:
\begin{eqnarray} \label{eq:3.3}
\Theta^{k+1} &=& \argmin_{\Theta} \quad  \Gc(X,\Theta) + \frac{\gamma}{2} \| \Theta- B_{0}\|^2_{F},  \label{eq:3.3a0}\\
Z^{k+1}_i &=& \argmin_{Z_i} \quad f_i(Z_i) + \frac{\gamma}{2} \| Z_i+Z_i^\top - B_{i}\|^2_{F},  \quad  \quad i=1, \ldots n, \label{eq:3.3a1}\\
E^{k+1} &=& \argmin_{E} \quad  f_e(E) +  \frac{\gamma}{2} \| E - B_{n+1}\|^2_{F},\label{eq:3.3a2}\\
\Lambda^{k+1} &=& {\Lambda}^k - \gamma({\Theta}^{k+1}- \sum_{i=1}^{n} Z_i^{k+1}+{Z_i^{k+1}}^\top- E^{k+1}).
\end{eqnarray}
where
\begin{eqnarray}\label{eq:3.4-2}
\nonumber
  B_{0} &=& \sum_{i=1}^{n}Z_i^{k}+{Z_i^{k}}^\top+ E^k +\frac{1}{\gamma} \Lambda^k,\\
\nonumber
  B_{1} &=& \Theta^{k+1}- (\sum_{i=2}^{n}Z_i^{k}+{Z_i^{k}}^\top+ E^k +\frac{1}{\gamma} \Lambda^k),\\
  \nonumber
  B_{i} &=& \Theta^{k+1}-(\sum_{j=1}^{i-1}Z_j^{k+1}+{Z_j^{k+1}}^\top \\
  \nonumber
  &+&\sum_{j=i+1}^{n} Z_j^{k}+{Z_j^{k}}^\top + E^k +\frac{1}{\gamma}\Lambda^k), \quad  \quad i=2, \ldots n-1,\\
  \nonumber
  B_{n} &=& \Theta^{k+1}-(\sum_{i=1}^{n-1} Z_i^{k+1}+{Z_i^{k+1}}^\top + E^k + \frac{1}{\gamma}\Lambda^k),\\
  B_{n+1} &=& \Theta^{k+1}-(\sum_{i=1}^{n} Z^{k+1}_i + {Z^{k+1}_i}^\top + \frac{1}{\gamma}\Lambda^k) .
\end{eqnarray}

To avoid introducing auxiliary variables and still solve subproblems \eqref{eq:3.3a1} efficiently,
 we propose to approximate the subproblems \eqref{eq:3.3a1} by linearizing the quadratic term of its objective function (see also
\citealp{Botle14,Linls11,yang13}). With this linearization, the resulting approximation to \eqref{eq:3.3a1} is then simple enough to have a closed-form solution. More specifically, letting $H_i(Z_i) = \frac{\gamma}{2} \|Z_i+ Z_i^\top - B_i\|_F^2$, we define  the following majorant function of $H_i(Z_i)$ at point $Z_i^k$,
\begin{eqnarray}\label{app1}
 H_i(Z_i) \leq \gamma \Big(\frac{1}{2} \|Z^{k}_i + {Z^{k}_i}^\top -B_i\|_F^2 + \langle \nabla H_i(Z^{k}_i), Z_i-Z^{k}_i \rangle + \frac{ \varrho}{2 }\| Z_i-Z^{k}_i\|_F^2 \Big),
\end{eqnarray}
where $\varrho$ is a proximal parameter, and
\begin{equation}\label{app2}
 \nabla H_i(Z^{k}_i) :=  2(Z^{k}_i + {Z^{k}_i}^\top)- (B_i+ B_i^\top),
\end{equation}

Plugging \eqref{app1} into \eqref{eq:3.3a1}, with simple algebraic manipulations, we obtain:
\begin{eqnarray} \label{eq:3.3lin}
Z^{k+1}_i &=& \argmin_{Z_i} \quad f_i(Z_i) + \frac{ \varrho \gamma}{2} \| Z_i - C_{i}\|_{F}^2, \quad \quad i=1, \ldots n,
\end{eqnarray}
where $ C_i = Z_i^k - \frac{1}{\varrho} \nabla H_i(Z^{k}_i)$.

The next result establishes the sufficient decrease property of the objective function given in \eqref{eq:3.3a1}, after a proximal map step computed in \eqref{eq:3.3lin}.

\begin{lem}\label{lemdec}(Sufficient decrease property).
Let $\varrho > \frac{L_{H_i}}{\gamma}$, where  $L_{H_i}$ is a Lipschitz constant of the gradient $\nabla H_i(Z_i)$ and $\gamma$ is a penalty parameter defined in \eqref{eq:3.1}. Then, we have
\begin{eqnarray*}
  f_i(Z_i^{k+1})+ H_i(Z_i^{k+1})&\leq& f_i(Z_i^{k})+ H_i(Z_i^{k}) - \frac{(\varrho \gamma - L_{H_i})}{2}\|Z_i^{k+1}-Z_i^{k}\|_F^2,   \qquad i=1, \dots, n,
\end{eqnarray*}
 where $Z_i^{k+1} \in \mathbb{R}^{n\times n}$ defined by \eqref{eq:3.3lin}.
\end{lem}
\textbf{Proof}. The proof of this Lemma follows along similar lines to the proof of Lemma~3.2 in \citet{Botle14}.

It is well known that \eqref{eq:3.3lin} has a closed-form solution that is given by the shrinkage operation \citep{Boyd11}:
\begin{eqnarray}\label{eq:3.5}
\nonumber
 Z^{k+1}_1 &=& \text{Shrink}\Big(C_1, \frac{\lambda_1 }{\varrho\gamma}\Big),\\
 Z^{k+1}_{i_j} &=& \max \big(1 - \dfrac{\lambda_i}{\varrho\gamma \| \text{Shrink}(C_{i_j}, \frac{\hat{\lambda}_i}{\varrho\gamma})\|_F}, 0 \big)\cdot \text{Shrink}(C_{i_j}, \frac{\hat{\lambda}_i}{\varrho\gamma}), \quad   \substack{ i=2, \ldots n, \\j = 1, \ldots, l_i,}
\end{eqnarray}
where $\text{Shrink}(\cdot, \cdot)$ in \eqref{eq:3.5} denotes the soft-thresholds operator, applied element-wise to a matrix A \citep{Boyd11}: $$\text {Shrink}(A_{ij},b):=\sign(A_{ij})\max \big(|A_{ij}|-b,0\big) \quad \quad \substack{i=1, \ldots p, \\j = 1, \ldots, p,}.$$

\begin{rem}
Note that in the case of solving problem \eqref{eq:10}, one needs to add another block function $f_{n+1}(Z_{n+1}) := \lambda_{n+1} \trace(Z_{n+1}) + \mathbbm{1}_{Z_{n+1} \succeq 0}$ to the augmented Lagrangian function \eqref{eq:3.1} and update $\{C_i\}_{i=1}^{n}$.
In this case, the proximal mapping of $f_{n+1}$ is
\begin{equation}\label{eq:34}
\text{prox}(f_{n+1}, \gamma, Z_{n+1}) := \argmin  \limits_{Z_{n+1}} f_{n+1}(Z_{n+1}) + \frac{\gamma}{2} \|Z_{n+1}- C_{n+1}\|^2_{F},
\end{equation}
where $C_{n+1} =\Theta^{k+1}-(\sum_{i=1}^{n} Z^{k+1}_i + {Z^{k+1}_i}^\top+ E^k+\frac{1}{\gamma}\Lambda^k).$
It is easy to verify that \eqref{eq:34} has a closed-form solution given by
$$
Z_{n+1} = U~\max(D-\frac{\lambda_{n+1}}{\gamma},0)~U^T,
$$
where $UDU^T $ is the eigenvalue decomposition of $C_{n+1}$ (see, \citealp{Chandrasekaran,ma13} for more details).
\end{rem}

The discussions above suggest that the following unmodified ADMM for solving \eqref{eq:3.1} gives rise to an efficient algorithm.

\begin{algorithm}[H]

\caption{Multi-Block ADMM Algorithm for Solving \eqref{eq:3.1}.}\label{alg:1}

\begin{algorithmic}[1]
    \Initialize{The parameters:}
    \begin{enumerate}[(a)]
    \item
        Primal variables $\Theta$, $Z_1$, $\dots$, $Z_n$,  $E$, to the $p \times p$ identity matrix.
    \item
        Dual variable $\Lambda$ to the $p \times p$ zero matrix.
    \item
        Constants $\varrho, \lambda_e, \tau > 0$,  and $\gamma\geq \sqrt{2} \lambda_e $.
    \item
        Nonnegative regularization constants $\lambda_1, \dots, \lambda_n,$ $\hat{\lambda}_2, \dots, \hat{\lambda}_n$.
    \end{enumerate}
    \Iterate{Until the stopping criterion $\| \Theta^k - \Theta^{k- 1}\|_F^2/\| \Theta^{k - 1} \|_{F} \leq \tau$ is met:}
    \begin{enumerate}[(a)]
    \item
    Update $\Theta$:
       \begin{enumerate}[]
    \item
    $\Theta^{k+1} = \argmin \limits_{\Theta \in \Sc}~ \Gc(X,\Theta) + \dfrac{\gamma}{2} \| \Theta - B_0\|^2_{F},$
        \end{enumerate}
     where $B_0$ is defined in \eqref{eq:3.4-2}.
     \item
     Update $Z_i$:
     \begin{enumerate}[i.]
     \item
          $  Z^{k+1}_1 = \text{Shrink}\Big(C_1, \frac{\lambda_1 }{\varrho\gamma}\Big),$
      \item $
 Z^{k+1}_{i_j} = \max \big(1 - \dfrac{\lambda_i}{\varrho\gamma \| \text{Shrink}(C_{i_j}, \frac{\hat{\lambda}_i}{\varrho\gamma})\|_F}, 0 \big)\cdot \text{Shrink}(C_{i_j}, \frac{\hat{\lambda}_i}{\varrho\gamma}), \quad   \substack{ i=2, \ldots n, \\j = 1, \ldots, l_i,}
$
     \end{enumerate}
     where $C_{i}$ is defined in \eqref{eq:3.3lin}.
     \item
     Update $E$:
       \begin{enumerate}[]
    \item
    $
 E^{k+1} = \argmin \limits_{E} ~ \frac{\lambda_e}{2}\|E\|^2_{F} + \frac{\gamma}{2} \|E - B_{n+1}\|^2_{F}$
        \end{enumerate}
        where $B_{n+1}$ is defined in \eqref{eq:3.4-2}.
     \item
     Update $\Lambda$:
     \begin{enumerate}[]
       \item   $ \Lambda^{k+1} = \Lambda^k -\gamma (\Theta^{k+1} -\sum_{i=1}^{n} Z^{k+1}_i + {Z^{k+1}_i}^\top-E^{k+1})
       $
     \end{enumerate}
  \end{enumerate}
  \end{algorithmic}
\end{algorithm}

\begin{rem}
The complexity of Algorithm \ref{alg:1} is of the same order as the graphical lasso \citep{Friedman07}, the method in \citet{Tan14} for hub node discovery and the algorithm used for estimation of sparse covariance matrices introduced by \citet{Xue12}. Indeed, one can easily see that with any set of structured matrices $\{Z_i\}_{i=1}^n$, the complexity of Algorithm \ref{alg:1} is equal to $O(p^3)$, which is the complexity of the eigen-decomposition for updating $\Theta$ in step~2(a).
\end{rem}

Since both the objective function and constraints of \eqref{eq:3.1} become separable after using the linearization technique introduced in \eqref{app1}, the problem can be decomposed into $n+2$ smaller subproblems; the latter can be solved in a parallel and distributed manner with a small modification in Algorithm~\ref{alg:1}. Indeed, we can apply a Jacobian ADMM to solve \eqref{eq:3.1} with the following updates,
\begin{eqnarray} \label{eqjacob:3.3}
\nonumber
\Theta^{k+1} &=& \argmin_{\Theta} \quad  \Gc(X,\Theta) + \frac{\gamma}{2} \| \Theta- B_{0}\|^2_{F}, \\
\nonumber
Z^{k+1}_i &=& \argmin_{Z_i} \quad f_i(Z_i) + \frac{ \varrho \gamma}{2} \| Z_i - C_{i}\|_{F}^2, \quad \quad i=1, \ldots n,
\\
\nonumber
E^{k+1} &=& \argmin_{E} \quad  f_e(E) +  \frac{\gamma}{2} \| E - B_{n+1}\|^2_{F},\\
\Lambda^{k+1} &=& {\Lambda}^k - \gamma({\Theta}^{k+1}- \sum_{i=1}^{n} Z_i^{k+1}+{Z_i^{k+1}}^\top- E^{k+1}).
\end{eqnarray}
where $ C_i $ is defined in \eqref{eq:3.3lin} with

\begin{eqnarray}\label{eqjacob:3.4-2}
\nonumber
  B_{0} &=& \sum_{i=1}^{n}Z_i^{k}+{Z_i^{k}}^\top+ E^k +\frac{1}{\gamma} \Lambda^k,\\
\nonumber
  B_{i} &=& \Theta^{k}-(\sum_{j=1}^{i-1}Z_j^{k}+{Z_j^{k}}^\top \\
  \nonumber
  &+&\sum_{j=i+1}^{n} Z_j^{k}+{Z_j^{k}}^\top + E^k +\frac{1}{\gamma}\Lambda^k), \quad  \quad i=2, \ldots n-1,\\
  \nonumber
  B_{n} &=& \Theta^{k}-(\sum_{i=1}^{n-1} Z_i^{k}+{Z_i^{k}}^\top + E^k + \frac{1}{\gamma}\Lambda^k),\\
  B_{n+1} &=& \Theta^{k}-(\sum_{i=1}^{n} Z^{k}_i + {Z^{k}_i}^\top + \frac{1}{\gamma}\Lambda^k).
\end{eqnarray}

Intuitively, the performance of the Jacobian ADMM should be worse than the Gauss-Seidel version, because the latter always uses the latest information of the primal variables in the updates. We refer to \citet{Liu15,Lin15} for a detailed discussion on the convergence analysis of the Jacobian ADMM and its variants. On the positive side, we obtain a parallelizable version of the multi-block ADMM algorithm.

\subsection{Convergence analysis}
The next result establishes the global convergence of the standard multi-block ADMM for solving SSON based statistical learning problems, by using the Kurdyka- Lojasiewicz (KL) property of the objective function in \eqref{eq:3.1}.
\begin{thm}\label{thm1}
The sequence $ U^k:=(\Theta^k, Z^k_1,\dots, Z^k_n,E^k,\Lambda^k)$ generated by Algorithm~\ref{alg:1} from any starting point converges to a stationary point of the problem given in~\eqref{eq:3.1}.
\end{thm}
\textit{Proof}. A detailed exposition is given in Appendix~\ref{apend2}.

\section{Experimental Results}

In this section, we present numerical results for Algorithm \ref{alg:1} (henceforth called SSONA), on both synthetic and real data sets.
The results are organized in the following three sub-sections: in Section \ref{Exp2}, we present numerical results on synthetic data comparing the performance of SSONA to that of grouping variables ADMM and also for assessing the accuracy in recovering a multi-layered structure in Markov Random Field and covariance graph models that constitute the prime focus in this paper.  In Section \ref{Exp1} we use the proposed SSONA for feature selection in classification problems involving two real data sets in order to calibrate SSON performance with respect to an independent validation set. Finally, in Section \ref{realdata}, we analyze using SSONA on some other interesting real data sets from the social and biological sciences.

\subsection{Experimental results for the SSON algorithm on graphical models based on synthetic data}\label{Exp2}

Next, we evaluate the performance of SSONA on ten synthetic graphical model problems, comprising of $ p = 100$, 500 and 1000 variables.  The underlying network structure corresponds to an Erd\H{o}s-R\'{e}nyi model graph, a nearest neighbor graph and a scale-free random graph, respectively. The CONTEST \footnote[1]{CONTEST is available at
{http://www.mathstat.strath.ac.uk/outreach/contest/}}
 package is used to generate the synthetic graphs, and the UGM \footnote[2]{UGM is available at  {http://www.di.ens.fr/~mschmidt/Software/UGM.html}} package to implement Gibbs sampling for estimating the Ising Model. Based on the generated graph topologies, we consider the following settings for generating synthetic data sets:
\begin{itemize}
  \item [I.]\textbf{Gaussian graphical models}: \\
For a given number of variables $p$, we first create a symmetric matrix $E \in \mathbb{R}^{p\times p}$ by using CONTEST in a MATLAB environment. Given matrix $E$, we set $\Sigma^{-1}$ equal to $E+(0.1- \bar{\Lambda}_{\min}(E))~I$, where $\bar{\Lambda}_{\min}(E)$ is the smallest eigenvalue of $E$ and $I$ denotes the identity matrix. We then draw $N=5p$ i.i.d. vectors $x_1,\dots,x_m$ from the gasserian distribution $\mathcal{N}(0,\Sigma)$ by using the \textit{mvnrnd} function in MATLAB, and then compute a sample covariance matrix of the variables.
\item [II.]\textbf{Gaussian graphical models with latent variables:}\\
For a given number of variables $p$, we first create a matrix $\Sigma^{-1} \in \mathbb{R}^{(p+r)\times (p+r)}$ by using CONTEST as described in I. We then choose the sub-matrix $\Theta_O =\Sigma^{-1}(1:p, 1:p)$ as the ground truth matrix of the matrix $\Theta_4$ and chose
\small{
\begin{eqnarray*}
  \Theta_U &=& \Sigma^{-1}(1:p,~ p+1:p+r)\big(\Sigma^{-1}(p+1:p+r,~p+r:p+r)\big)^{-1} \\
    && \Sigma^{-1}(p+1:p+r,~1:p)
\end{eqnarray*}
}
as the ground truth matrix of the low rank matrix $U$. We then draw $N=5p$ i.i.d. vectors $x_1,\dots,x_m$ from the Gaussian distribution $\mathcal{N}(0, (\Theta_O-\Theta_U)^{-1})$, and compute the sample covariance matrix of the variables $\Sigma_O$.
\item [III.] \textbf{The Binary Network}:\\
To generate the parameter matrix $\Sigma$, we create an adjacency matrix as in Setup I by using CONTEST. Then, each of $N= 5p$ observations is generated through Gibbs sampling. We take the first 100000 iterations as our burn-in period, and then collect observations, so that they are nearly independent.
\end{itemize}

 We compare SSONA to the following competing methods:
\begin{itemize}
\item \textbf{CovSel},  designed to estimate a sparse Gaussian graphical model \citep{Friedman07};
\item \textbf{HGL},     focusing on learning a Gaussian graphical model having hub nodes \citep{Tan14};
\item \textbf{PGADM},   designed to learn a Gaussian graphical model with some latent nodes \citep{ma13};
\item \textbf{Pseudo-Exact}, designed to learn a binary Ising graphical model \citep{Hof09};
\item \textbf{glasso-SF}, Learning Scale Free Networks by reweighted $\ell_1$ Regularization \citep{Liu11};
\item \textbf{GADMM}, A two block ADMM method with grouping variables.
\end{itemize}

All the algorithms have been implemented in the MATLAB R2015b environment on a PC with a 1.8 GHz processor and 6GB RAM memory.  Further, all the algorithms are being terminated either when
$$
\dfrac{\| \Theta^k - \Theta^{k-1}\|^2_{F}}{\| \Theta^{k-1} \|^2_{F}} \leq \tau, \qquad \tau = 1e-5,
$$
or the number of iterations and CPU times exceed 1,000 and 10 minutes, respectively.

We found that in practice the computation cost for SSONA increases with the size of structured matrices. Therefore, we use a limited memory version of SSONA in our experimental results to obtain good accuracy. Block sizes in Figure~\ref{fig2} could be set based on a desire for interpretability of the resulting estimates. In this section, we choose four structured matrices with blocks of size
\begin{eqnarray*}
(Z_2)_j &=& [1,\frac{p}{2}],  \qquad j=1 \dots, l_2,\\
(Z_3)_j &=& [1,\frac{p}{5}], \qquad j=1 \dots, l_3,\\
(Z_4)_j &=& [1,\frac{p}{10}], \qquad j=1 \dots, l_4,\\
(Z_5)_j &=& [1,\frac{p}{20}],  \qquad j=1 \dots, l_5,
\end{eqnarray*}
where $l_i$ is determined based on size of the adjacency matrix, $p$ (see, Figure~\ref{fig2}).

\begin{figure}[!ht]
 \begin{subfigure}[h]{.46\textwidth}
    {\includegraphics[width=\textwidth]{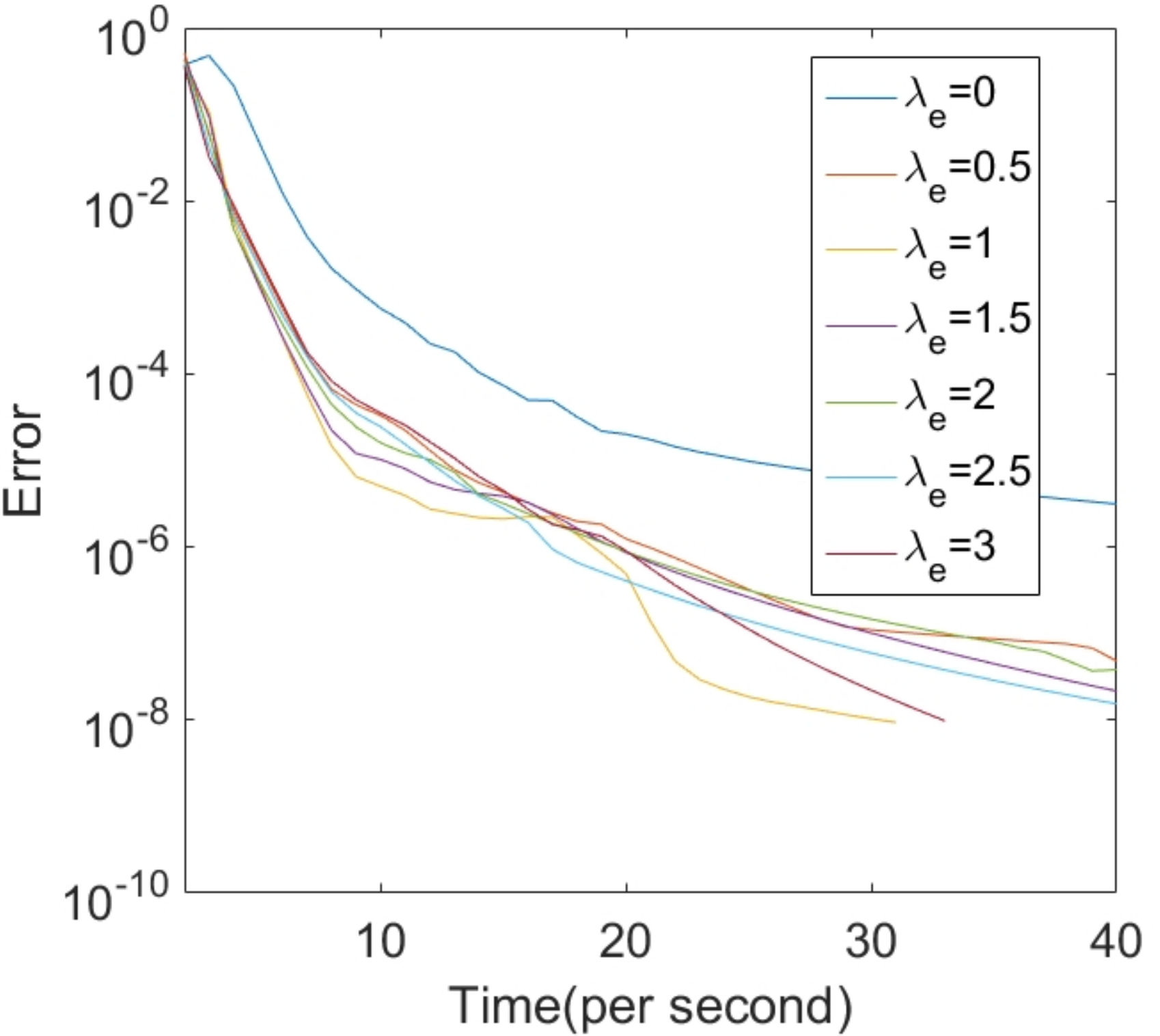}}
    \caption{}
  \end{subfigure}
  \hfill
   \begin{subfigure}[h]{.5\textwidth}
  {\includegraphics[width=\textwidth]{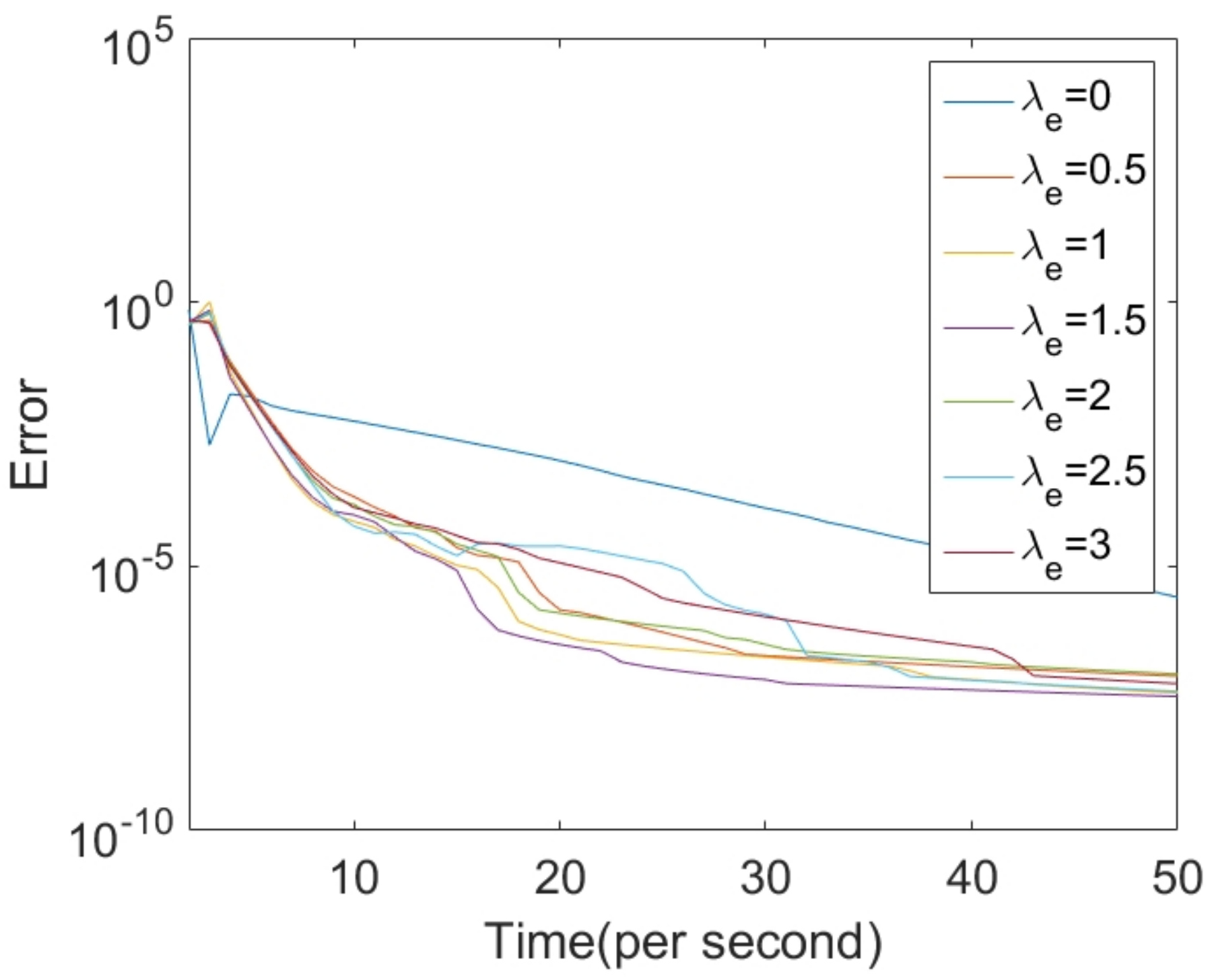}}
    \caption{}
 \end{subfigure}
\caption{Learning turning parameter $\lambda_e$ for two covariance estimation problems. Comparison of the absolute errors produced by the algorithms based on CPU time for different choices of $\lambda_e$.}
\label{figlame}
\end{figure}

\begin{figure}[!ht]
 \begin{subfigure}[h]{.47\textwidth}
    {\includegraphics[width=\textwidth]{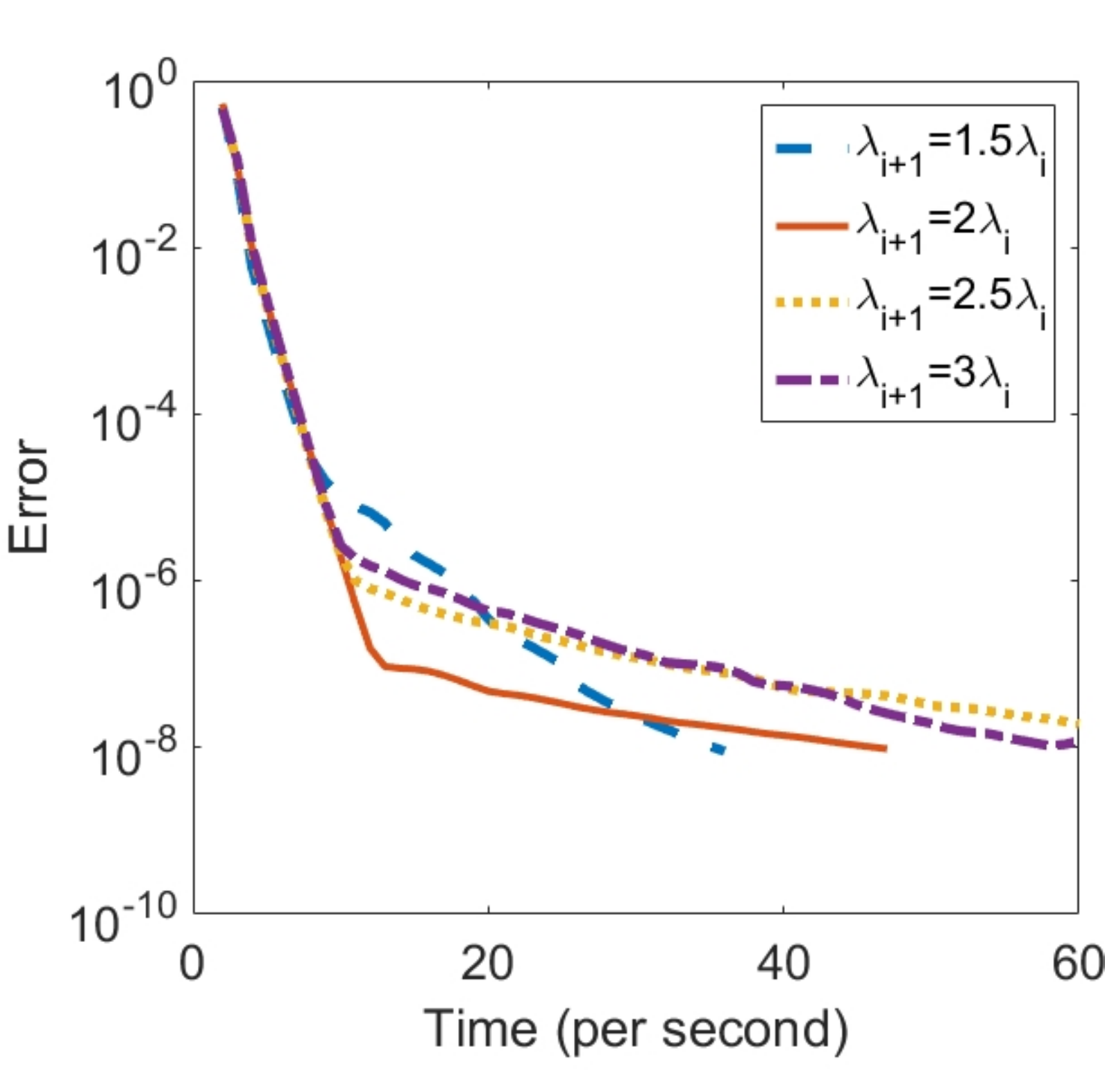}}
    \caption{}
  \end{subfigure}
  \hfill
   \begin{subfigure}[h]{.5\textwidth}
  {\includegraphics[width=\textwidth]{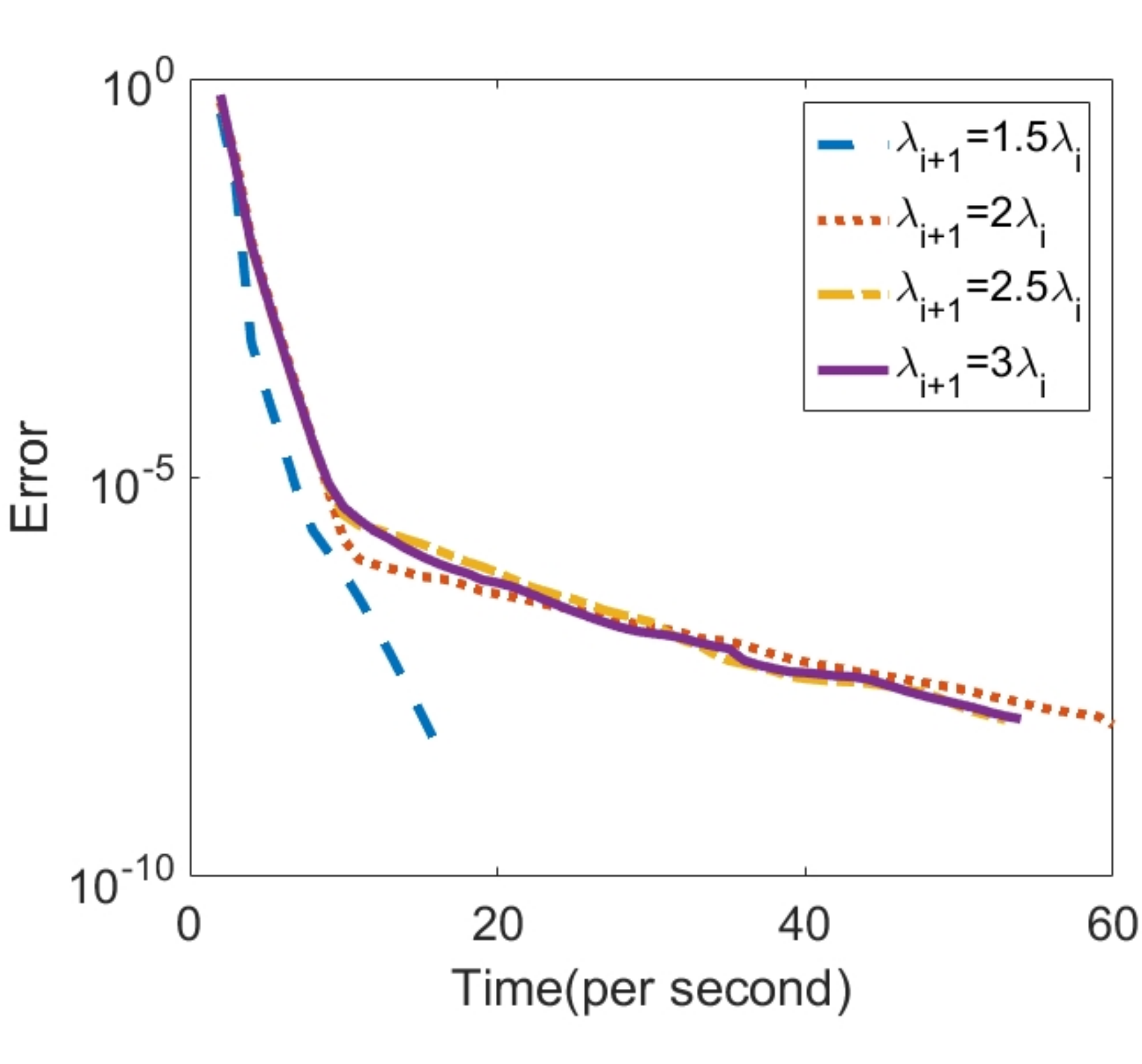}}
    \caption{}
 \end{subfigure}
\caption{Learning turning parameter $\lambda_i$ for $i=2, \dots, n$ for two covariance estimation problems for different choices of $\lambda_i$ for $i=2, \dots, n$.}
\label{figlamle}
\end{figure}

The penalty parameters $\lambda_e$ and $\{\lambda_i\}_{i=1}^{n}$ play an important rule for the convex decomposition to be successful. We learn them through numerical experimentation (see Figures \ref{figlame} and \ref{figlamle}) and set them respectively to
$$
\varrho =4, \quad \lambda_e=1,  \quad  \lambda_1,\lambda_2=0.5\lambda_e, \quad \hat{\lambda}_i=  0.25\lambda_e, \quad \text{and} \quad \lambda_{i+1}= 2\lambda_i \quad \text{for} \quad i=2, \dots, n.
$$

It can be seen from Figure~\ref{figlame} that with the addition of the ridge penalty term $\frac{\lambda_e}{2}\|E\|^2_{F}$ the algorithm clearly outperforms its unmodified counterpart in terms of CPU time for any fixed number of iterations. Indeed, when the model becomes more dense, SSONA is more effective to recover the network structure.

\begin{figure}[!ht]
 \begin{subfigure}[h]{.46\textwidth}
    {\includegraphics[width=\textwidth]{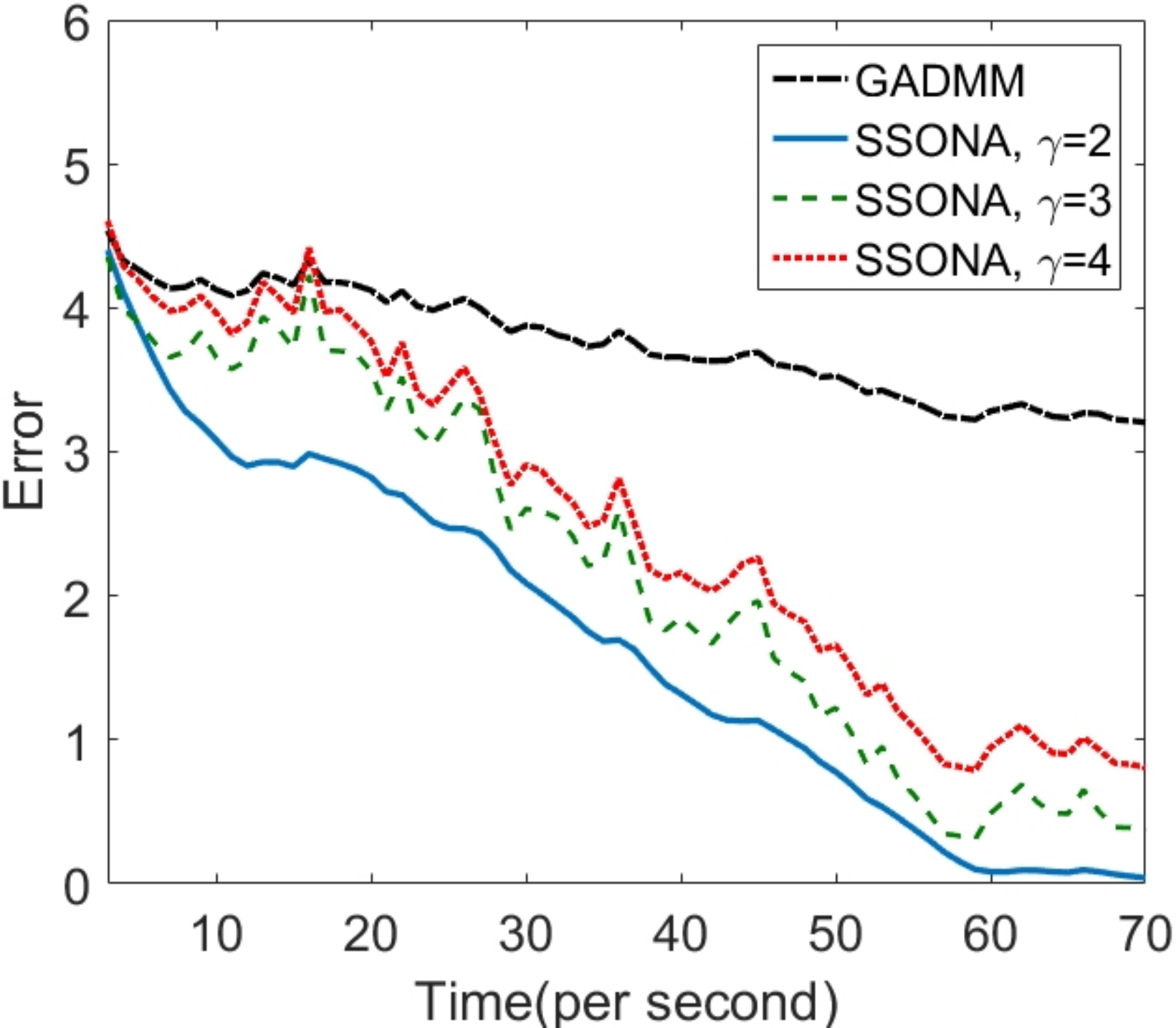}}
    \caption{}
  \end{subfigure}
  \hfill
   \begin{subfigure}[h]{.48\textwidth}
  {\includegraphics[width=\textwidth]{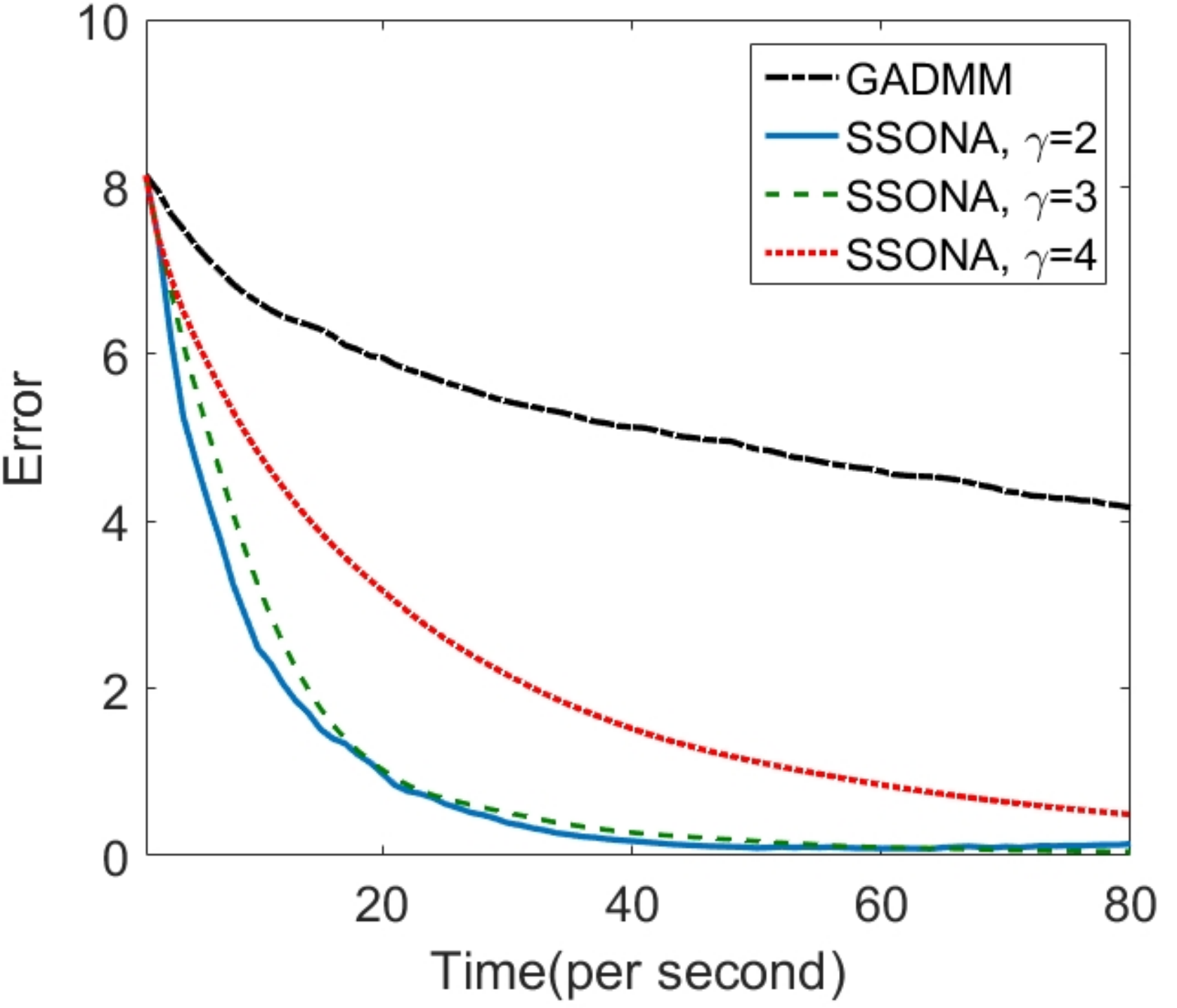}}
    \caption{}
 \end{subfigure}
\caption{Comparison of the absolute errors produced by the algorithms based on CPU time for different choices of $\gamma$.}
\label{figadmm}
\end{figure}

Next, we conduct experiments to assess the performance of the developed multi-block ADMM algorithm (SSONA) vis-a-vis the GADMM for solving two covariance graph estimation problems of dimension 1000 in the presence of noise.  Figure~\ref{figadmm} depicts the absolute error of the objective function for different choices of the regularization parameter $\gamma$ of the augmented Lagrangian and that of the dense noisy component $\lambda_e$; note that the latter is key for the convergence of the proposed algorithm.

We define the following two performance measures, as proposed in \citet{Tan14}:
\begin{itemize}
  \item  Number of correctly estimated edges, $n_e$  :
  $$
\sum_{j<j'}\Big( 1_{\{|\hat{\Theta}|> 1e-4 ~\text{and}~ |\Theta_{jj'}|\neq0  \}}\Big).$$
\item  Sum of squared errors, $s_e$:
$$
\sum_{j < j'}\Big( | \hat{\Theta}_{jj'}-\Theta_{jj'}| \Big)^2.$$
\end{itemize}

\begin{figure}[!ht]
 \begin{subfigure}[h]{0.47\textwidth}
    {\includegraphics[width=\textwidth]{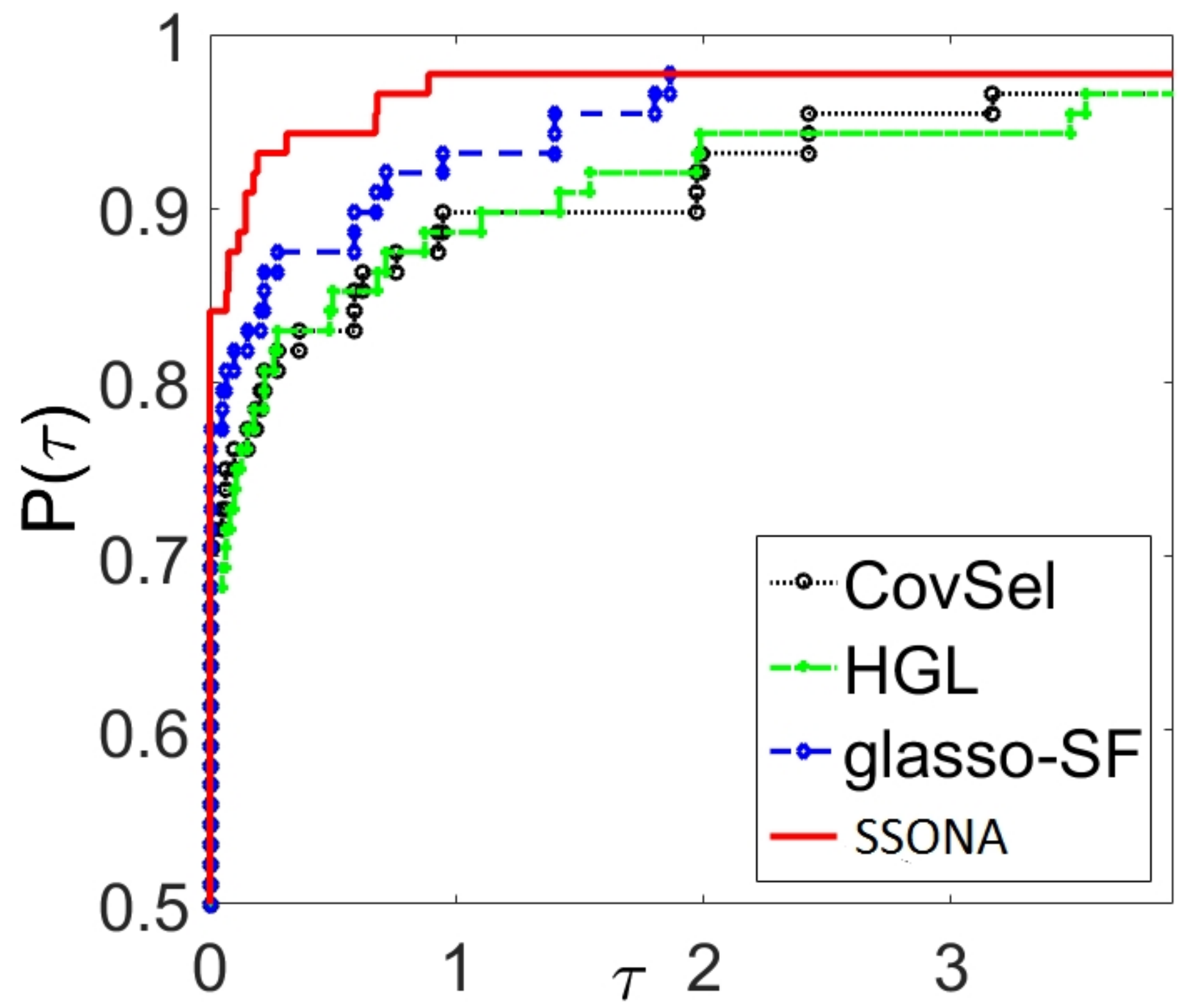}}
    \caption{Performance based on $n_e$.}
  \end{subfigure}
  \hfill
  \begin{subfigure}[h]{0.47\textwidth}
  {\includegraphics[width=\textwidth]{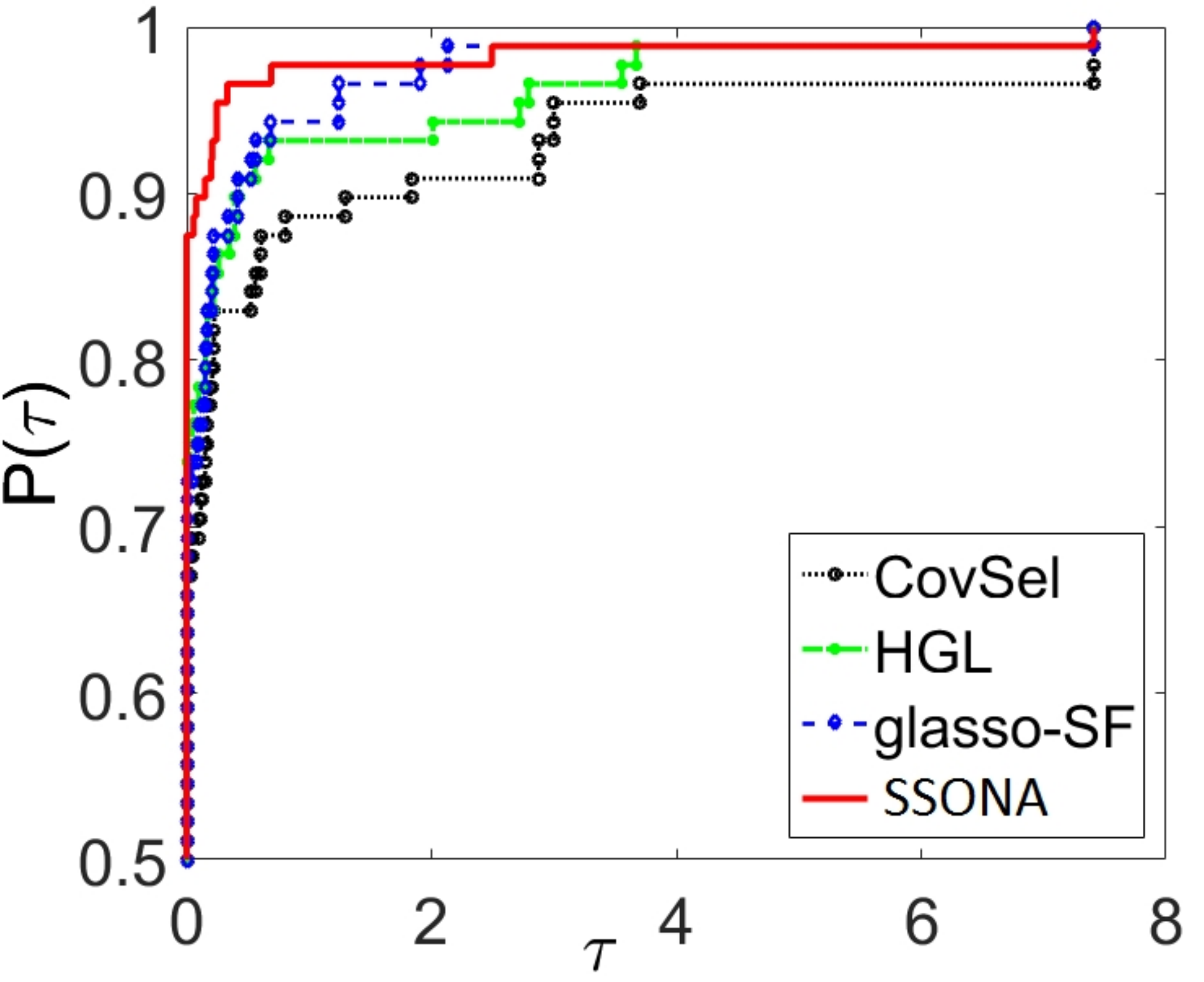}}
    \caption{Performance based on $s_e$.}  \end{subfigure}
\caption{Performance profiles of CovSel, HGL, glasso-SF and SSONA}
\label{figpr1}
\end{figure}

\begin{figure}[!ht]
 \begin{subfigure}[h]{0.46\textwidth}
    {\includegraphics[width=\textwidth]{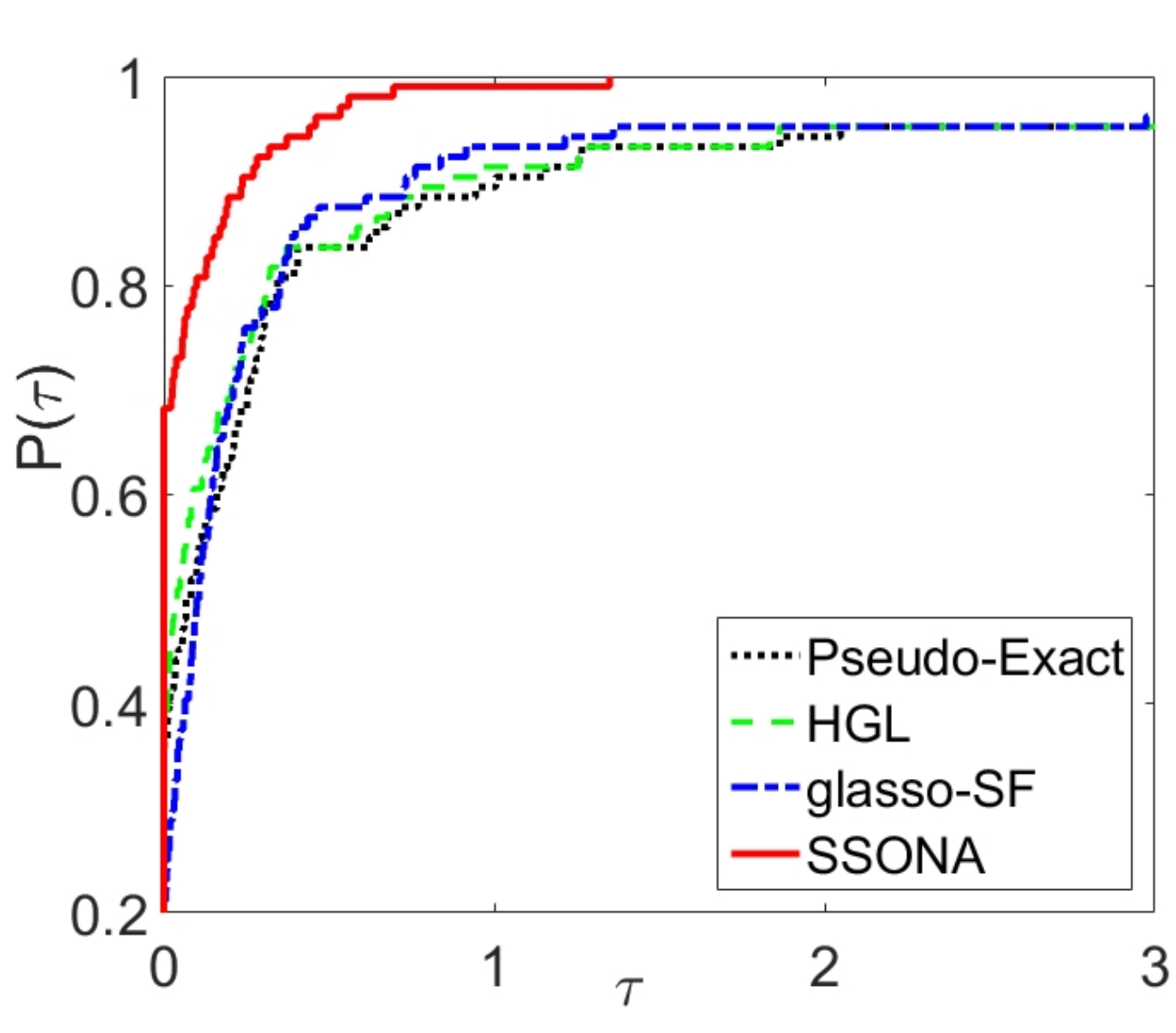}}
    \caption{Performance based on $n_e$.}
  \end{subfigure}
  \hfill
  \begin{subfigure}[h]{0.46\textwidth}
  {\includegraphics[width=\textwidth]{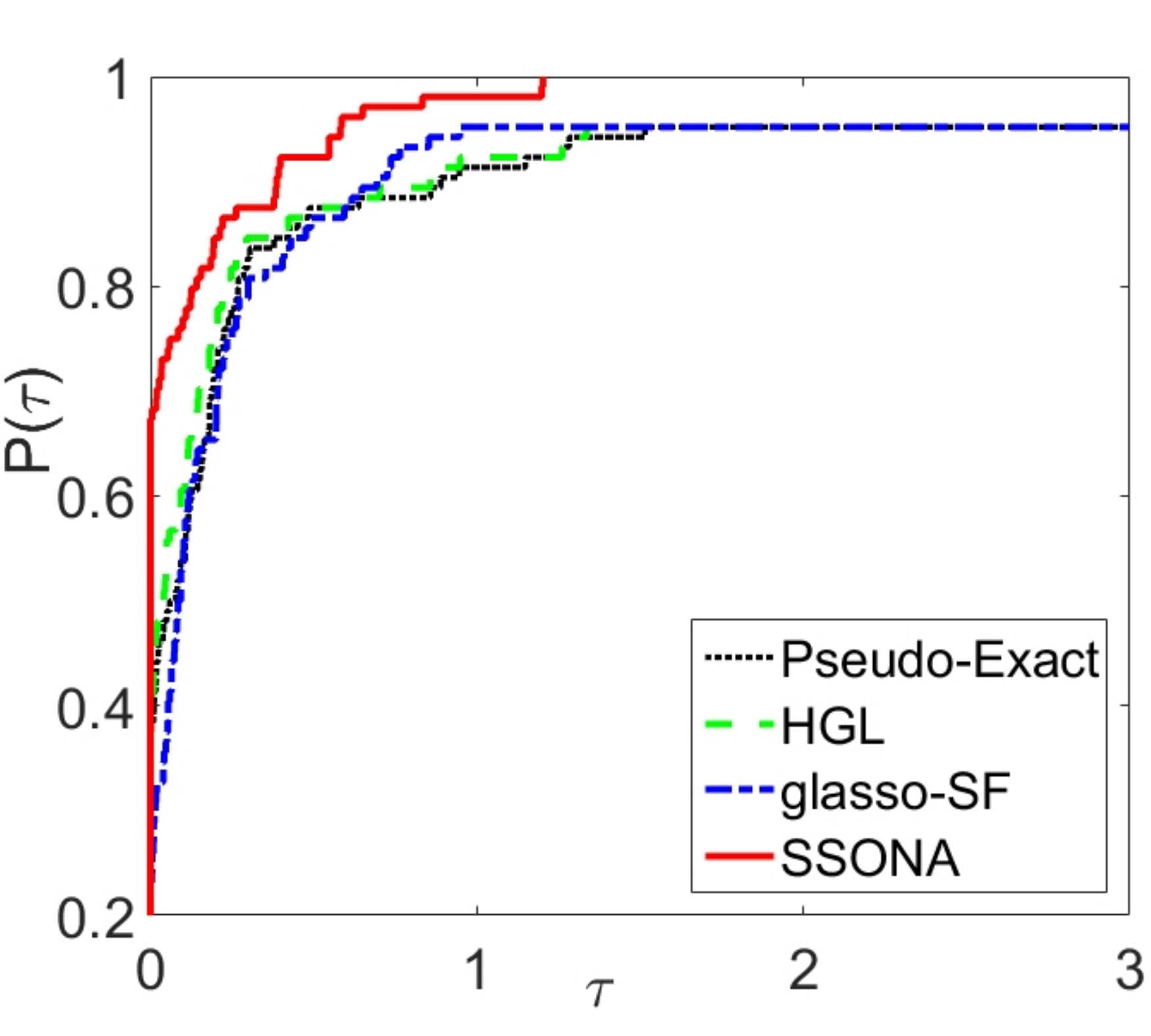}}
    \caption{Performance based on $s_e$.}  \end{subfigure}
\caption{Performance profiles of Pseudo-Exact, HGL, glasso-SF and SSONA.}
\label{figpr2}
\end{figure}

\begin{figure}[!ht]
 \begin{subfigure}[h]{0.47\textwidth}
    {\includegraphics[width=\textwidth]{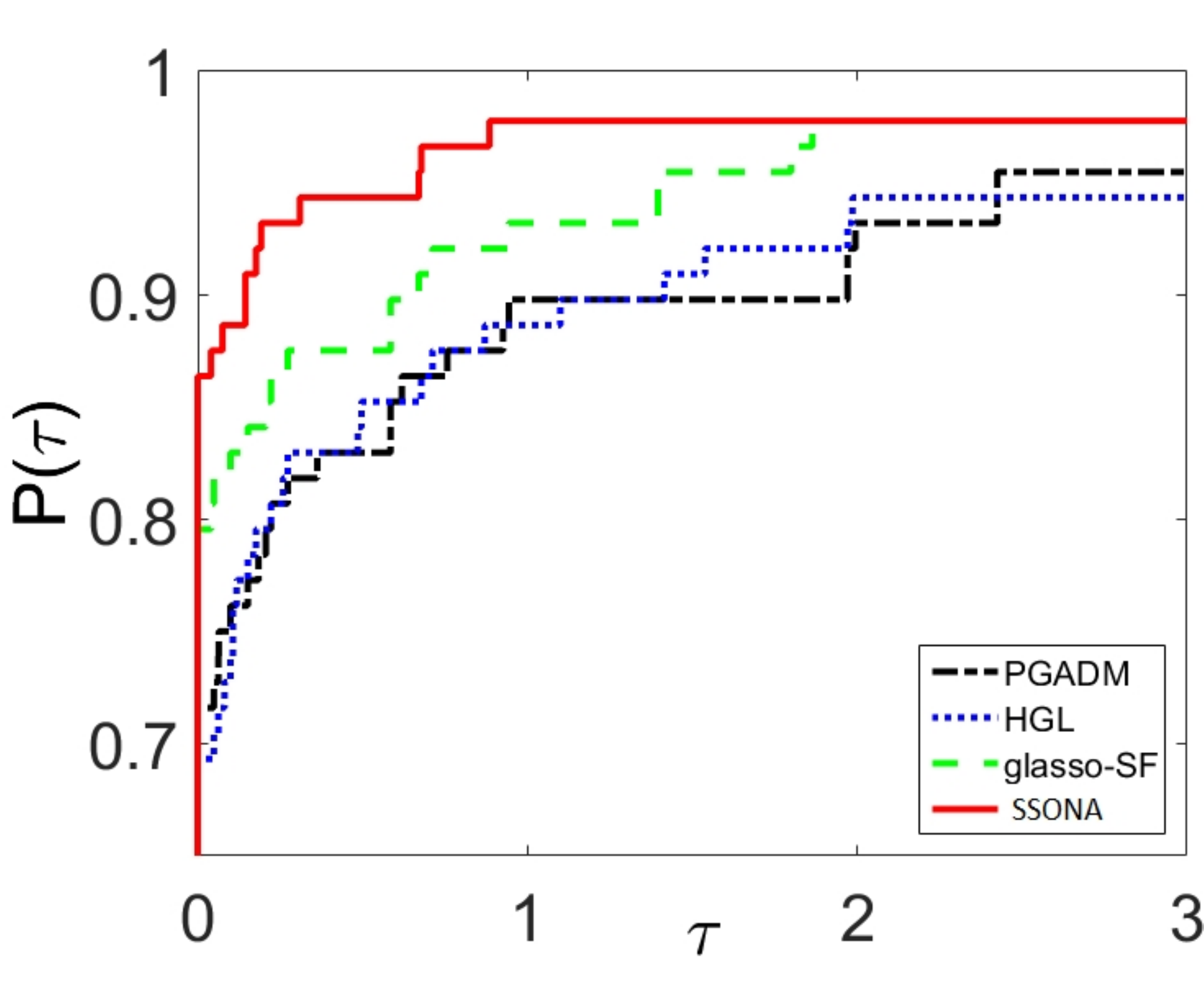}}
    \caption{Performance based on $n_e$.}
  \end{subfigure}
  \hfill
  \begin{subfigure}[h]{0.47\textwidth}
  {\includegraphics[width=\textwidth]{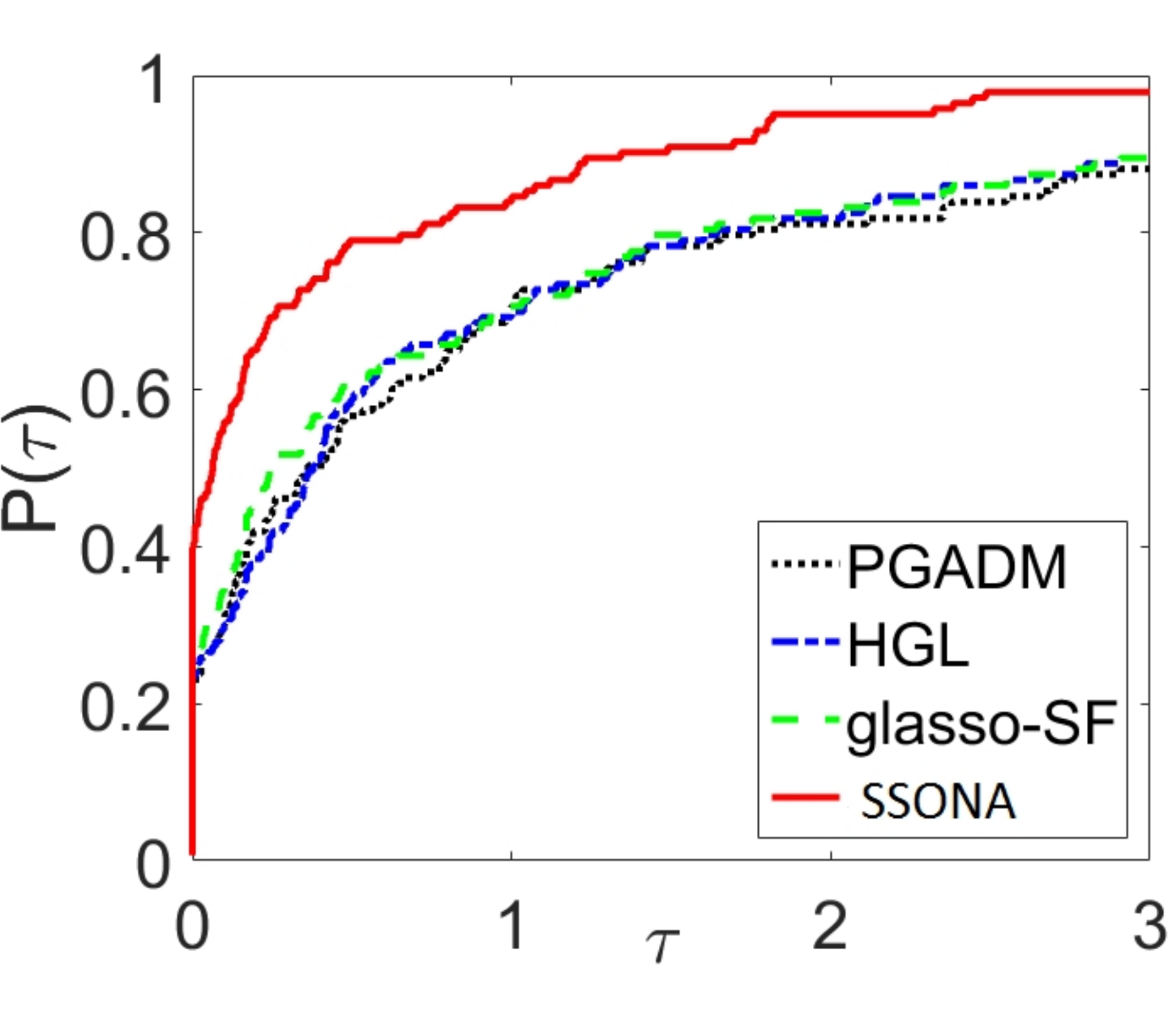}}
    \caption{Performance based on $s_e$.}
     \end{subfigure}
\caption{Performance profiles of PGADM, HGL, glasso-SF and SSONA.}
\label{figpr3}
\end{figure}

The experiment is repeated ten times and the average number of correctly estimated edges, $n_e$  and sum of squared errors, $s_e$ are considered for comparison. We have used the {\em performance profile}, as proposed in \citet{Dolan02}, to display the efficiency of the algorithms considered, in terms of $n_e$  and $s_e$. As stated in \citet{Dolan02}, this profile provides a wealth of information such as solver efficiency, robustness and probability of success in compact form and eliminates the influence of a small number of problems on the evaluating process and the sensitivity of results associated with the ranking of solvers. Indeed, the performance profile plots the fraction of problem instances for which any given method is within a factor of the best solver. The horizontal axis of the figure gives the percentage of the test problems for which a method is efficient, while the vertical axis gives the percentage of the test problems that were successfully solved by each method (robustness). The performance profiles of the considered algorithms in log2 scale are depicted in Figures~\ref{figpr1},\ref{figpr2} and \ref{figpr3}.

Figures~\ref{figpr1},\ref{figpr2} and \ref{figpr3} show the performance profiles of the considered algorithms for estimation of graphical models in terms of number of correctly estimated edges and sum of squared errors, respectively. The left and right panel are drawn in terms of $n_e$ and $s_e$, respectively. The results in these figures clearly demonstrate the superior performance of the proposed method,
since it solves all test problems without exhibiting any failure. Moreover, the SSONA algorithm is the best algorithm among the considered ones, as it solves more than 80 \% of the test problems achieving the maximum number of correctly estimated edge $n_e$ and minimum value of
estimation loss $s_e.$ Further, the performance index of SSONA grows up rapidly in comparison with the other considered algorithms. The latter implies that whenever SSONA is not the best algorithm, its performance index is close to the index of the best one.

\subsubsection{Experiments on structured graphical models}
\begin{figure}[!ht]
 \begin{subfigure}[h]{0.3\textwidth}
    {\includegraphics[width=\textwidth]{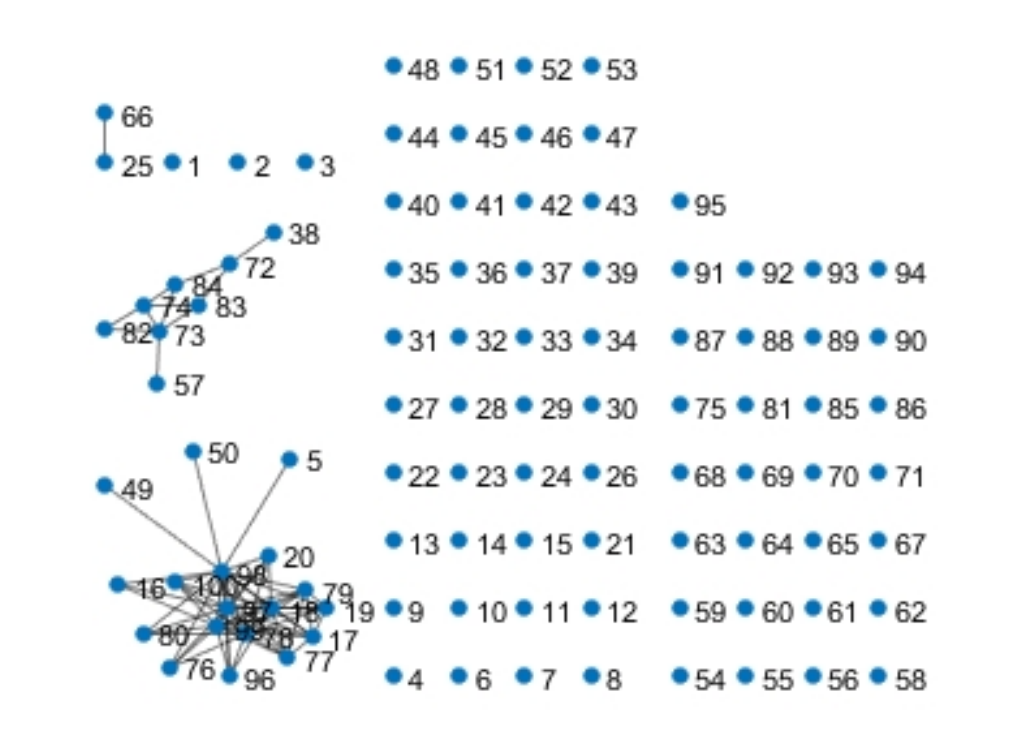}}
    \caption{Graphical lasso.}
  \end{subfigure}
  \hfill
  \begin{subfigure}[b]{0.28\textwidth}
  {\includegraphics[width=\textwidth]{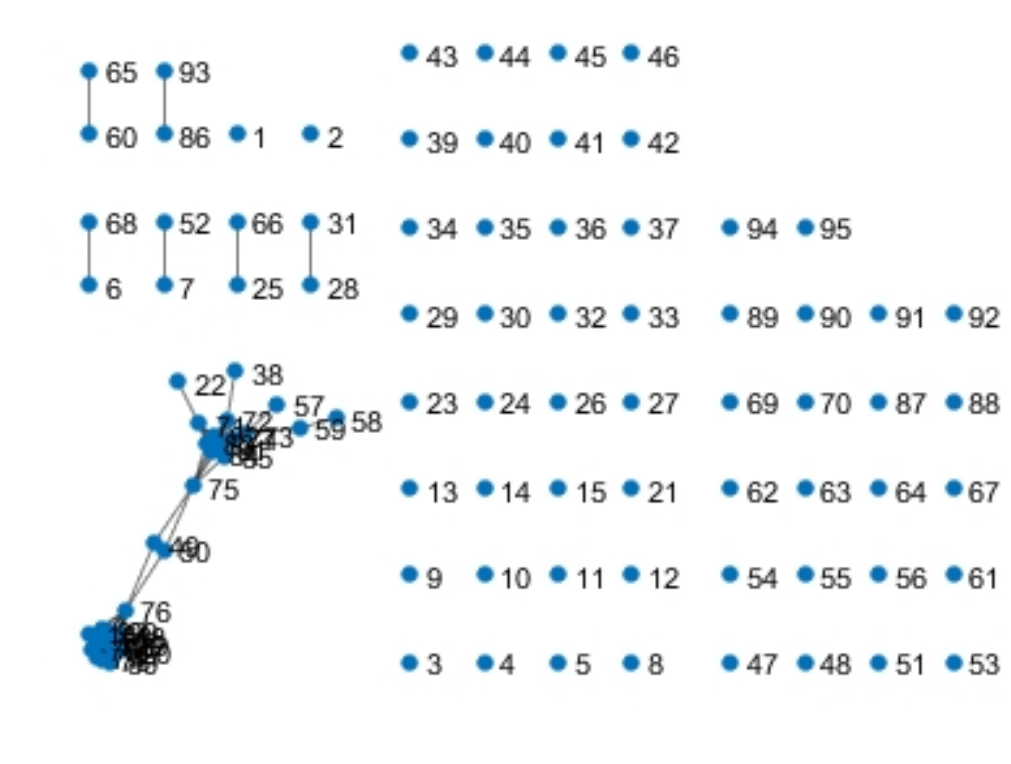}}
    \caption{SSONA.}
    \end{subfigure}
  \hfill
  \begin{subfigure} {0.28\textwidth}
  {\includegraphics[width=\textwidth]{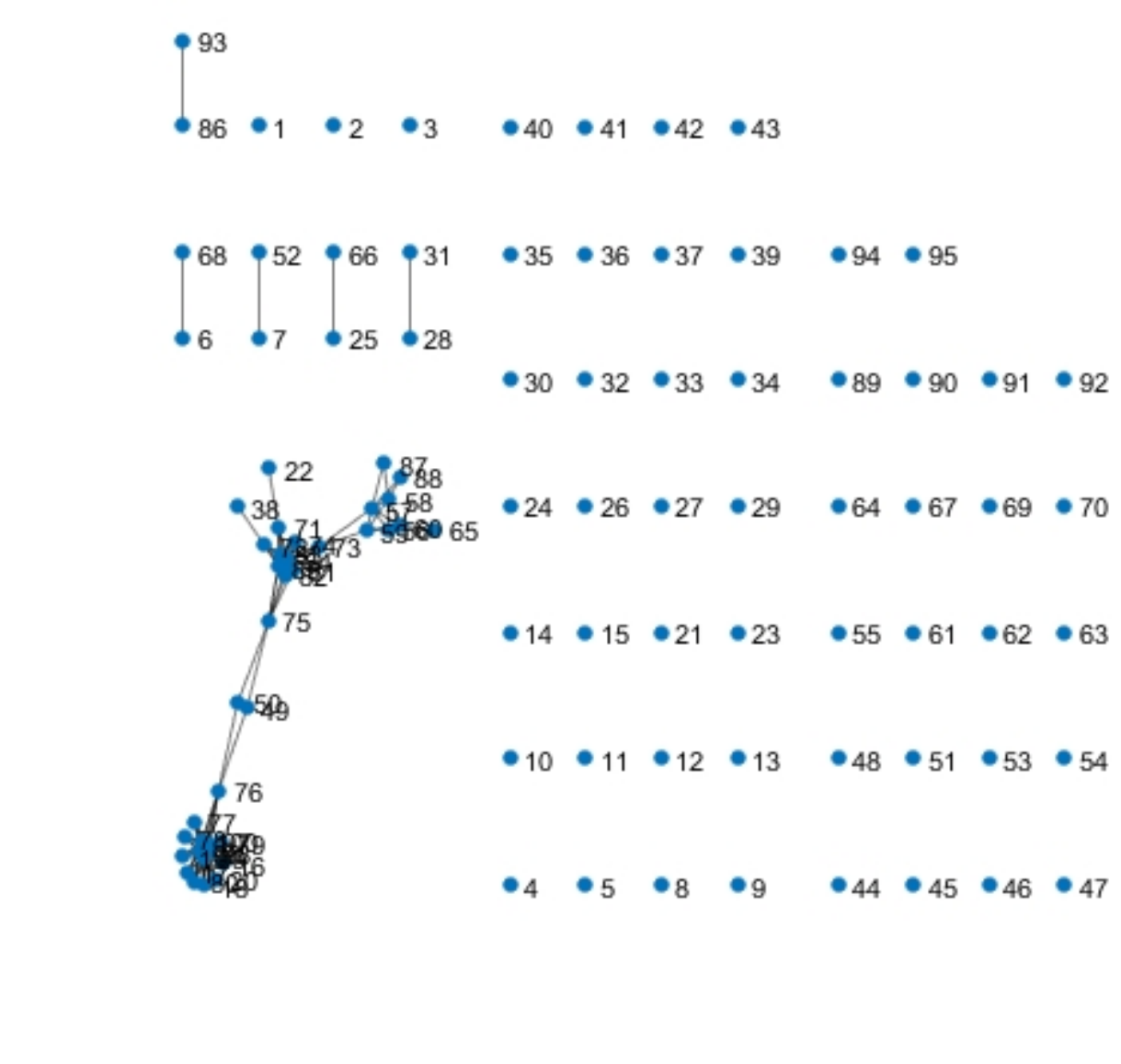}}
    \caption{Ground truth.}
  \end{subfigure}
 \begin{subfigure} {0.3\textwidth}
    {\includegraphics[width=\textwidth]{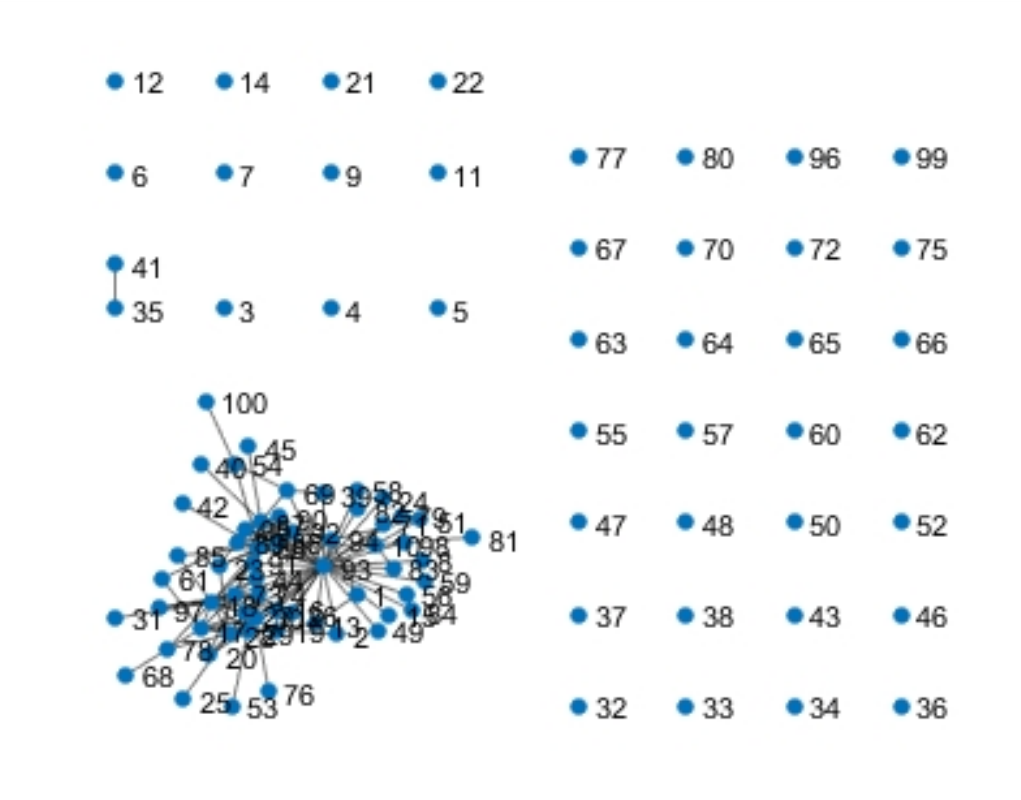}}
    \caption{Graphical lasso.}
  \end{subfigure}
  \hfill
  \begin{subfigure} {0.3\textwidth}
  {\includegraphics[width=\textwidth]{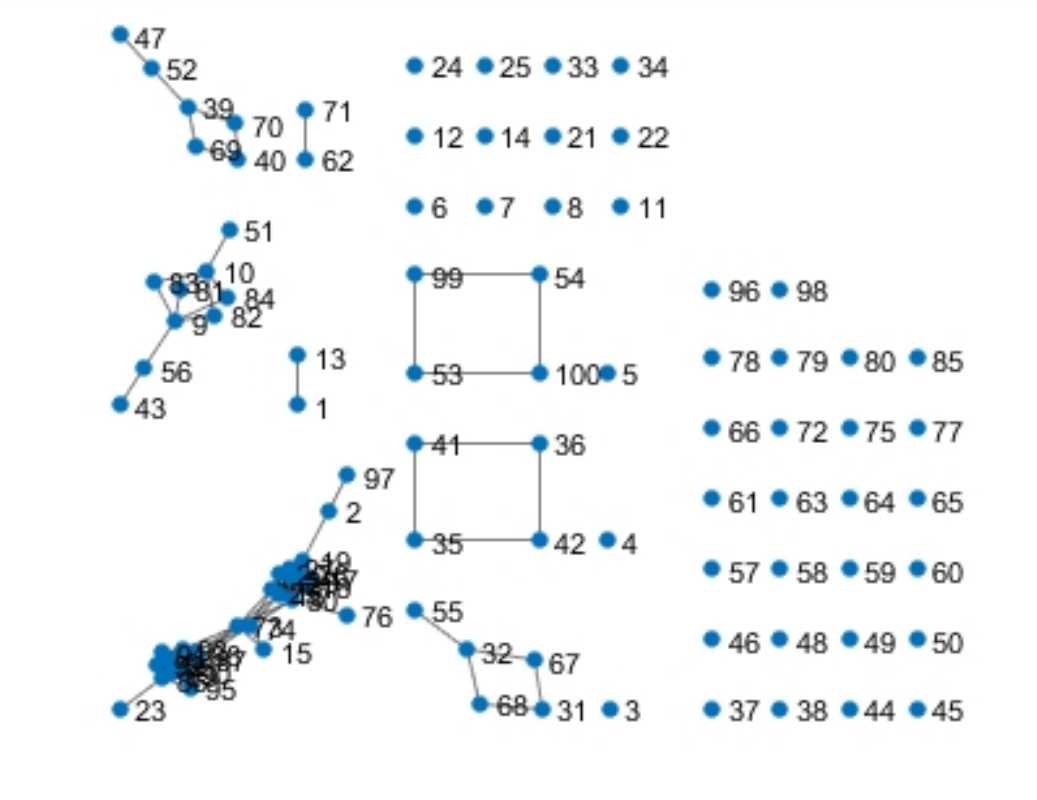}}
    \caption{SSONA.}
  \end{subfigure}
  \hfill
  \begin{subfigure} {0.28\textwidth}
  {\includegraphics[width=\textwidth]{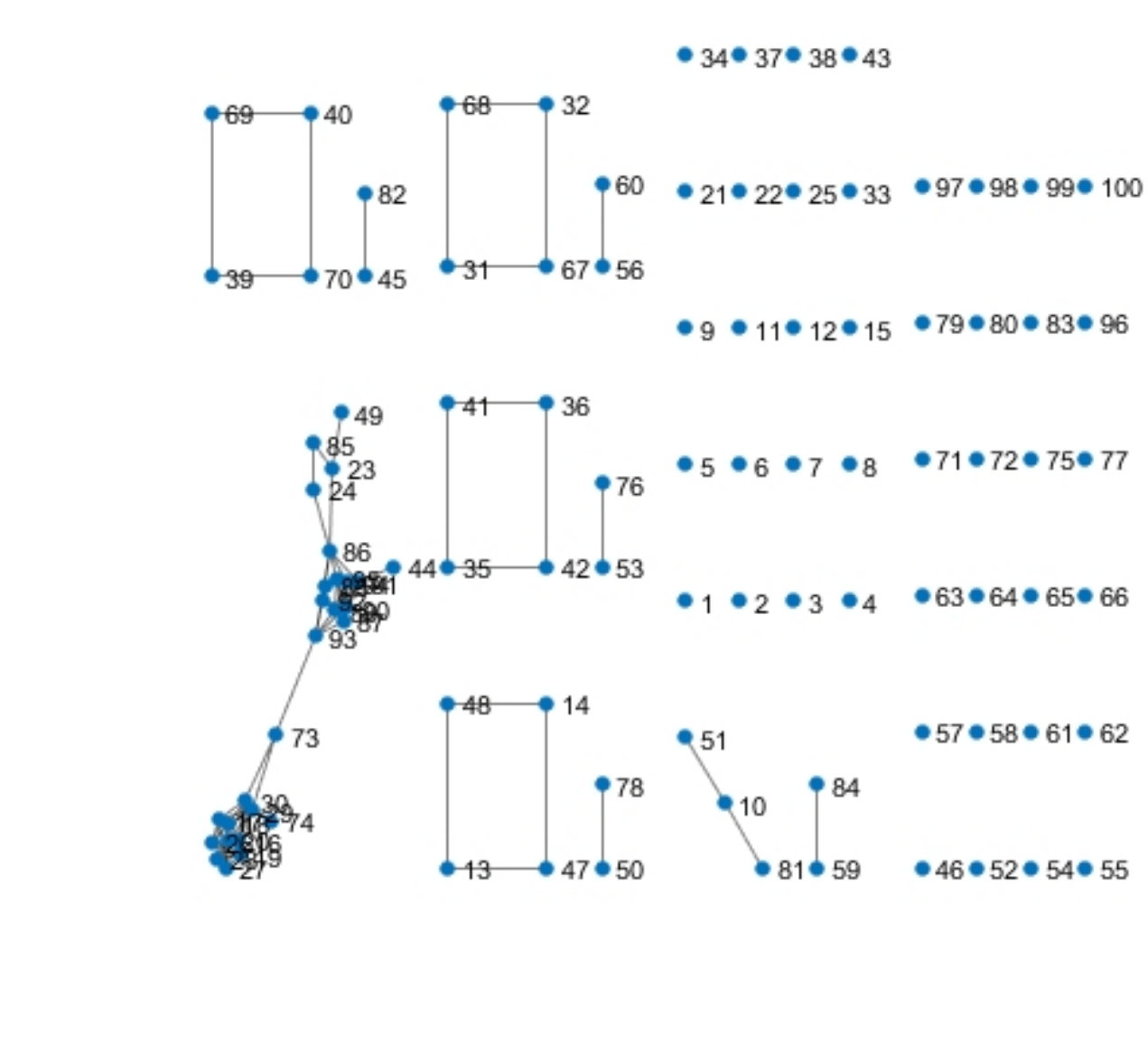}}
    \caption{Ground truth.} 
  \end{subfigure}
\caption{Simulation for the Gaussian graphical model. Row~I: Results for $p=100$ and $m=200$. Row~II:
Results for $p=100$ and $m=100$.}\label{figs1}
\end{figure}

\begin{figure}[!ht]
 \begin{subfigure} {0.3\textwidth}
    {\includegraphics[width=\textwidth]{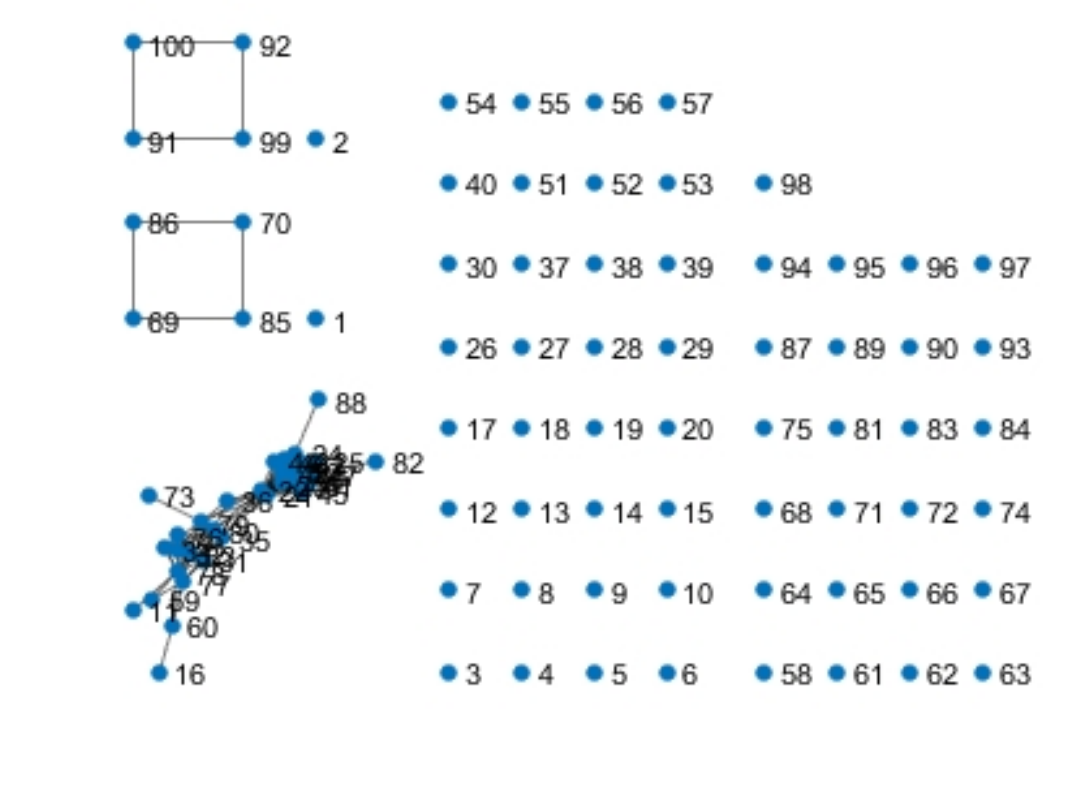}}
    \caption{Graphical lasso.}
  \end{subfigure}
  \hfill
  \begin{subfigure} {0.29\textwidth}
  {\includegraphics[width=\textwidth]{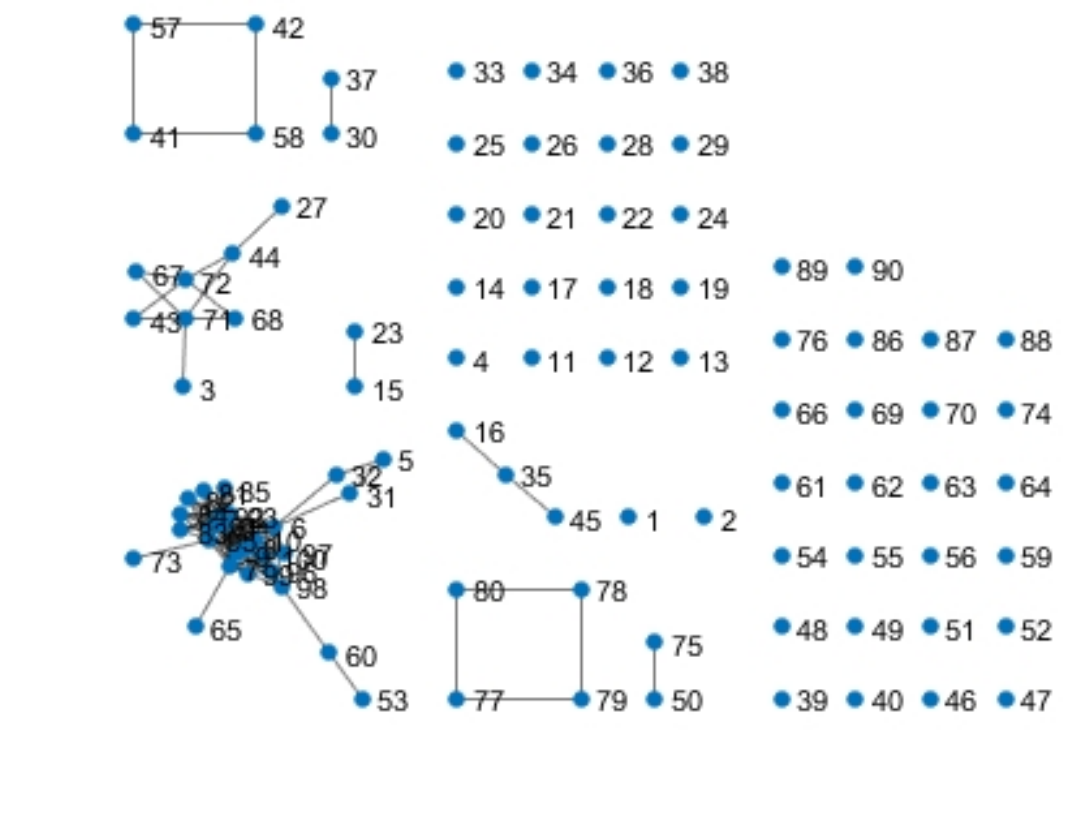}}
    \caption{SSONA.}  \end{subfigure}
  \hfill
  \begin{subfigure} {0.3\textwidth}
  {\includegraphics[width=\textwidth]{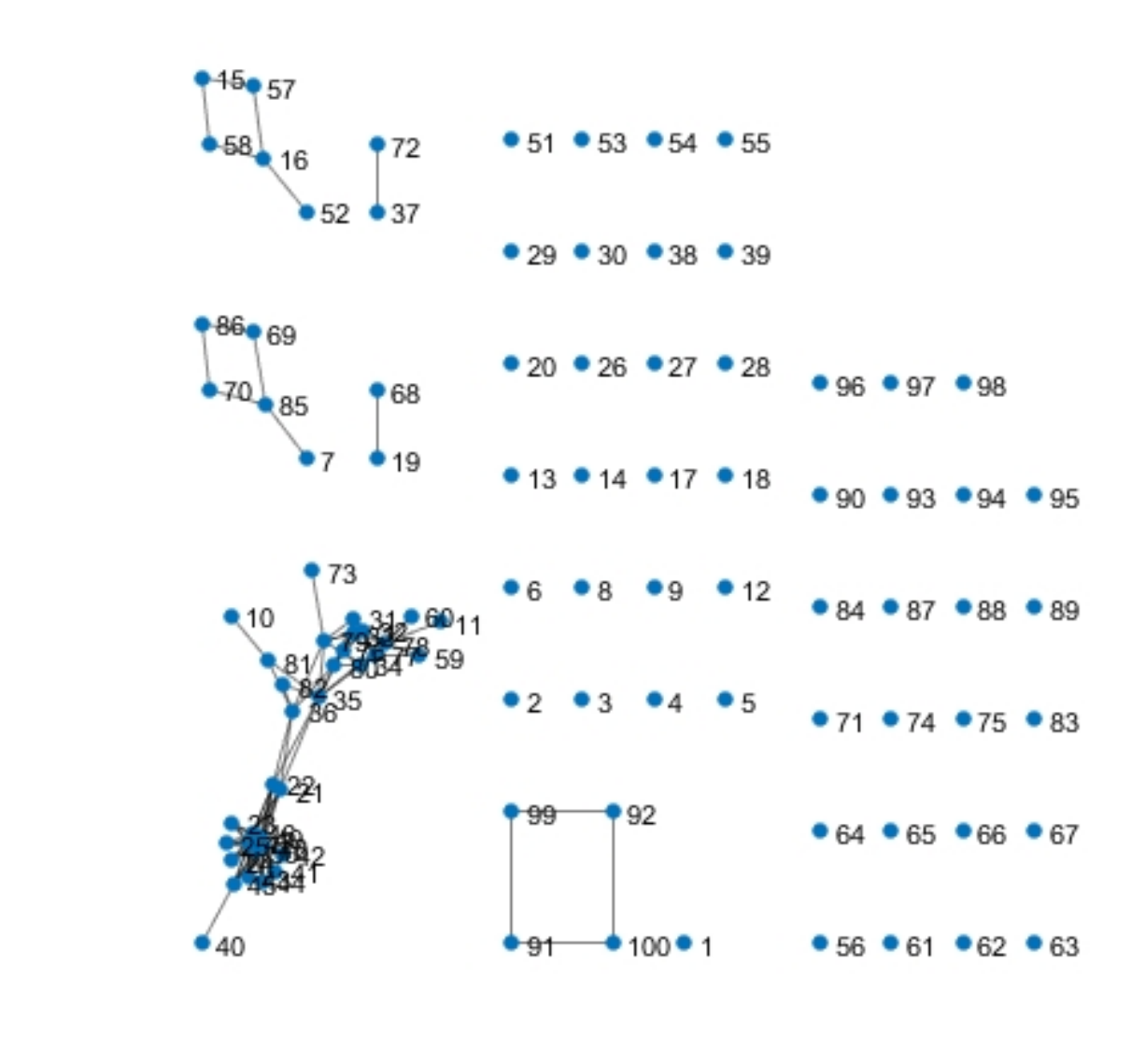}}
    \caption{Ground truth.}
  \end{subfigure}
 \begin{subfigure} {0.3\textwidth}
    {\includegraphics[width=\textwidth]{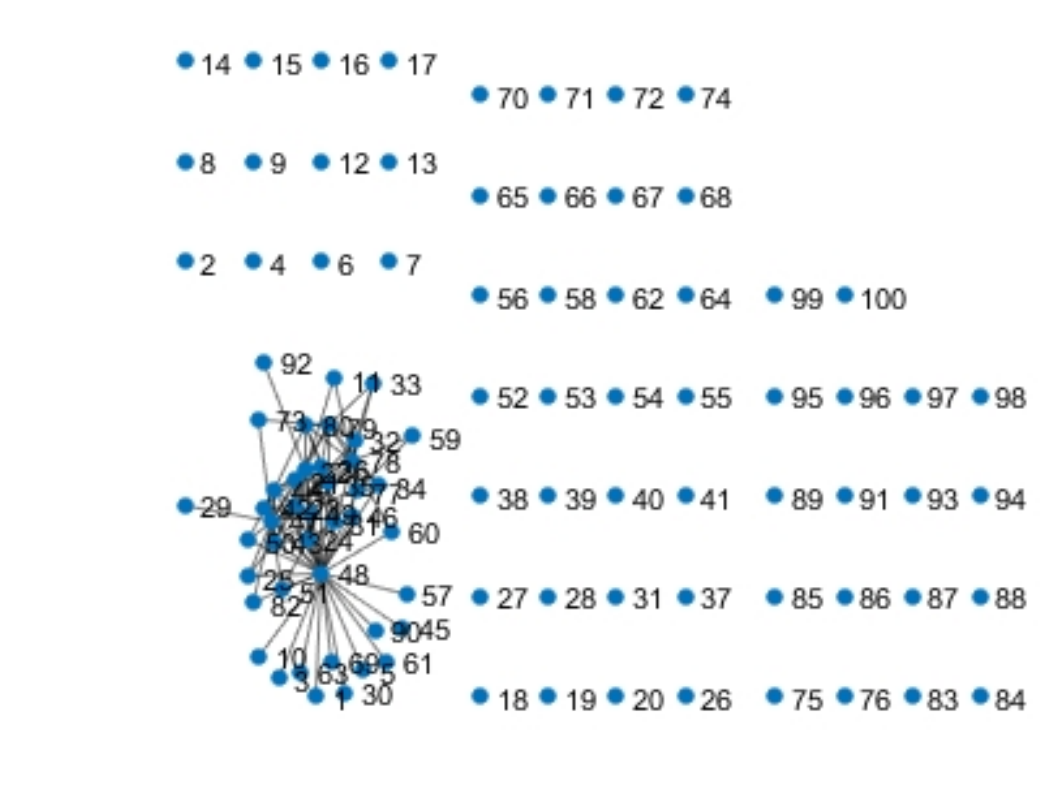}}
    \caption{Graphical lasso.}
  \end{subfigure}
  \hfill
  \begin{subfigure} {0.3\textwidth}
  {\includegraphics[width=\textwidth]{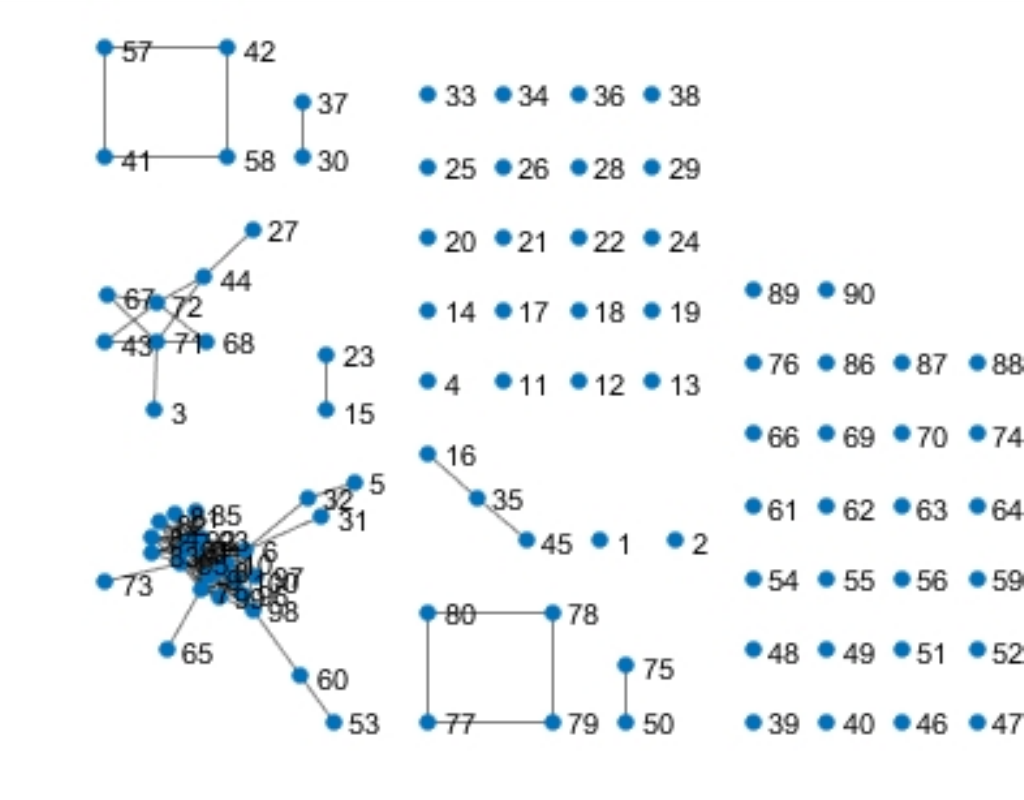}}
    \caption{SSONA.}
  \end{subfigure}
  \hfill
  \begin{subfigure} {0.3\textwidth}
  {\includegraphics[width=\textwidth]{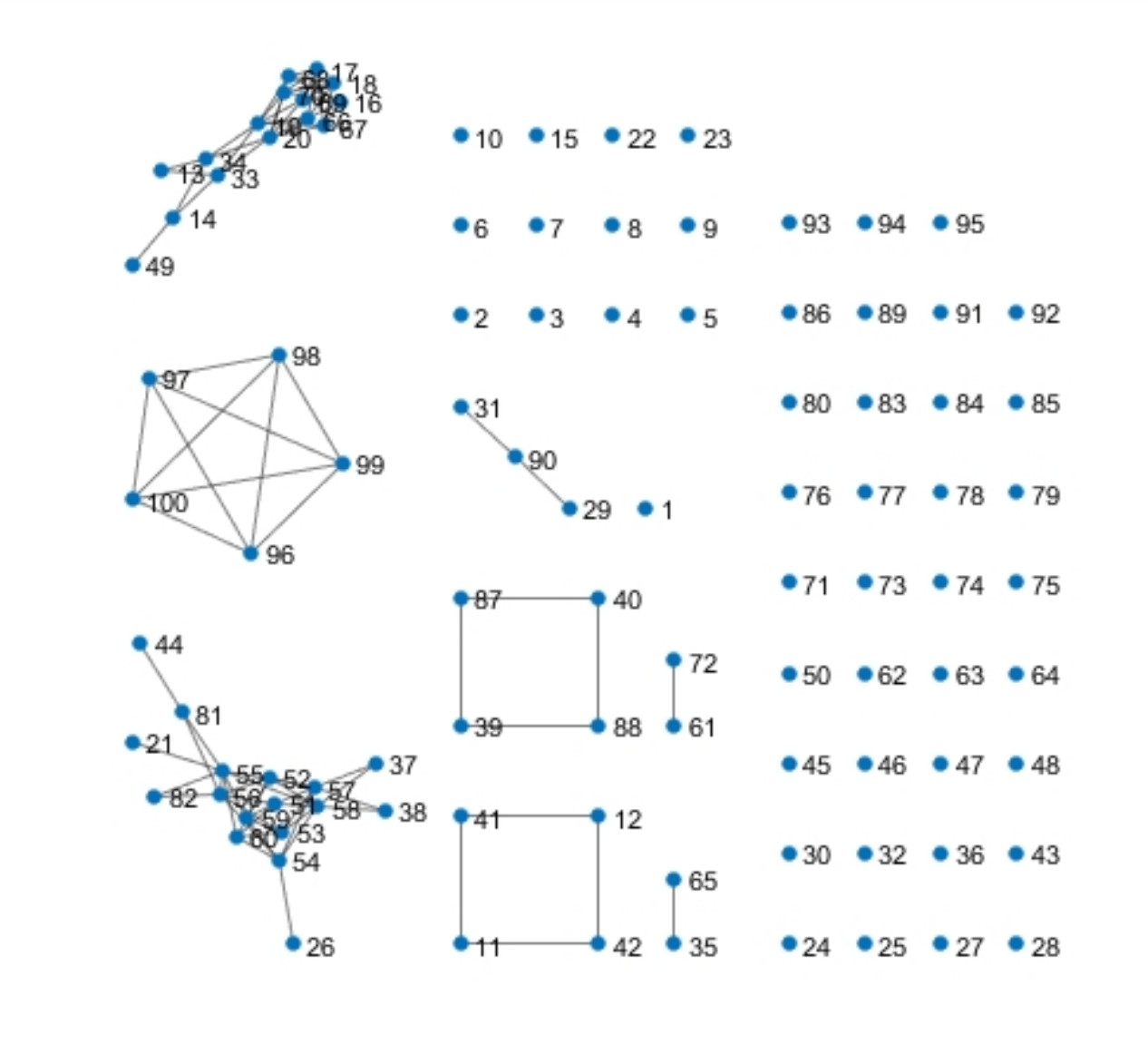}}
    \caption{Ground truth.} 
  \end{subfigure}
  \caption{ \small{Simulation for the Covariance graph model. Row~I: Results for $p=100$ and $m=200$. Row~II:
Results for $p=100$ and $m=100$.}} \label{figs2}
\end{figure}

\begin{figure}[!ht]
 \begin{subfigure} {0.3\textwidth}
    {\includegraphics[width=\textwidth]{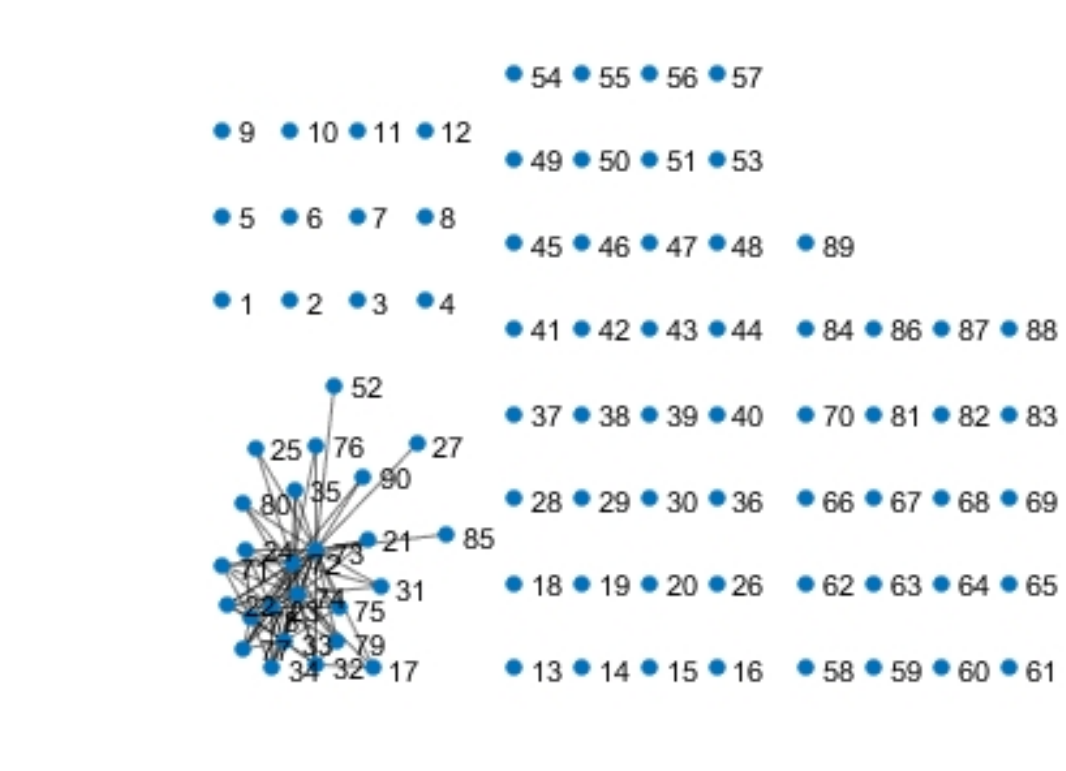}}
    \caption{Graphical lasso.}
  \end{subfigure}
  \hfill
  \begin{subfigure} {0.29\textwidth}
  {\includegraphics[width=\textwidth]{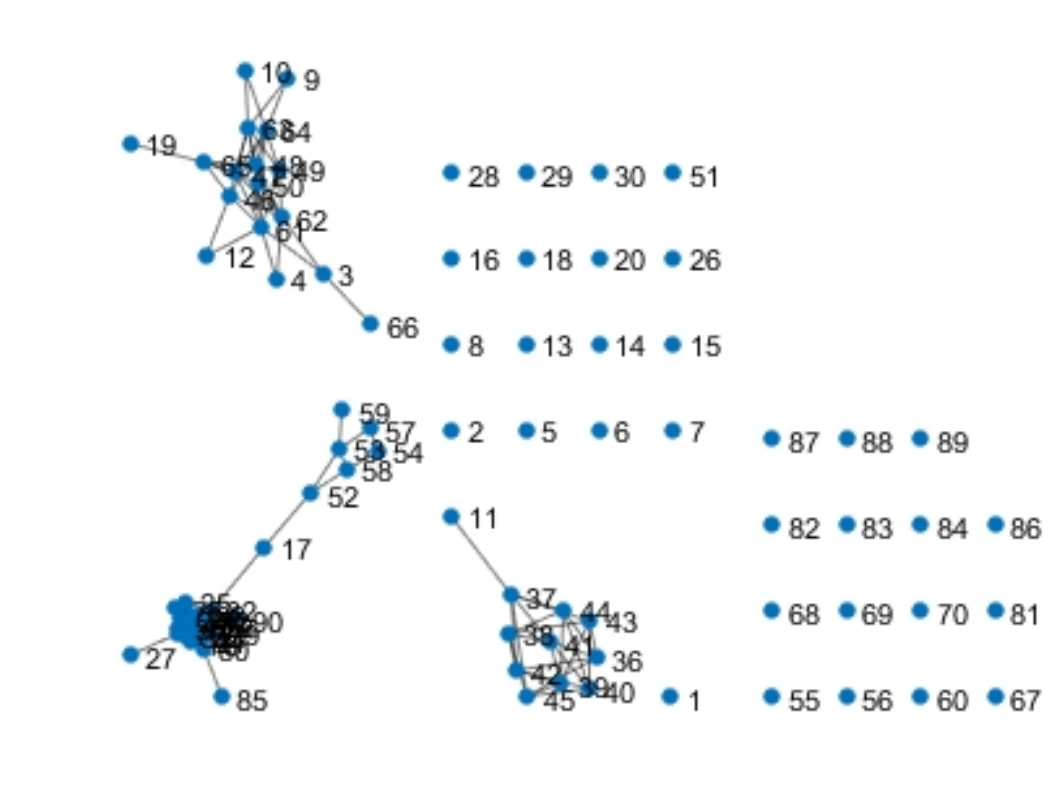}}
    \caption{SSONA.}  \end{subfigure}
  \hfill
  \begin{subfigure} {0.3\textwidth}
  {\includegraphics[width=\textwidth]{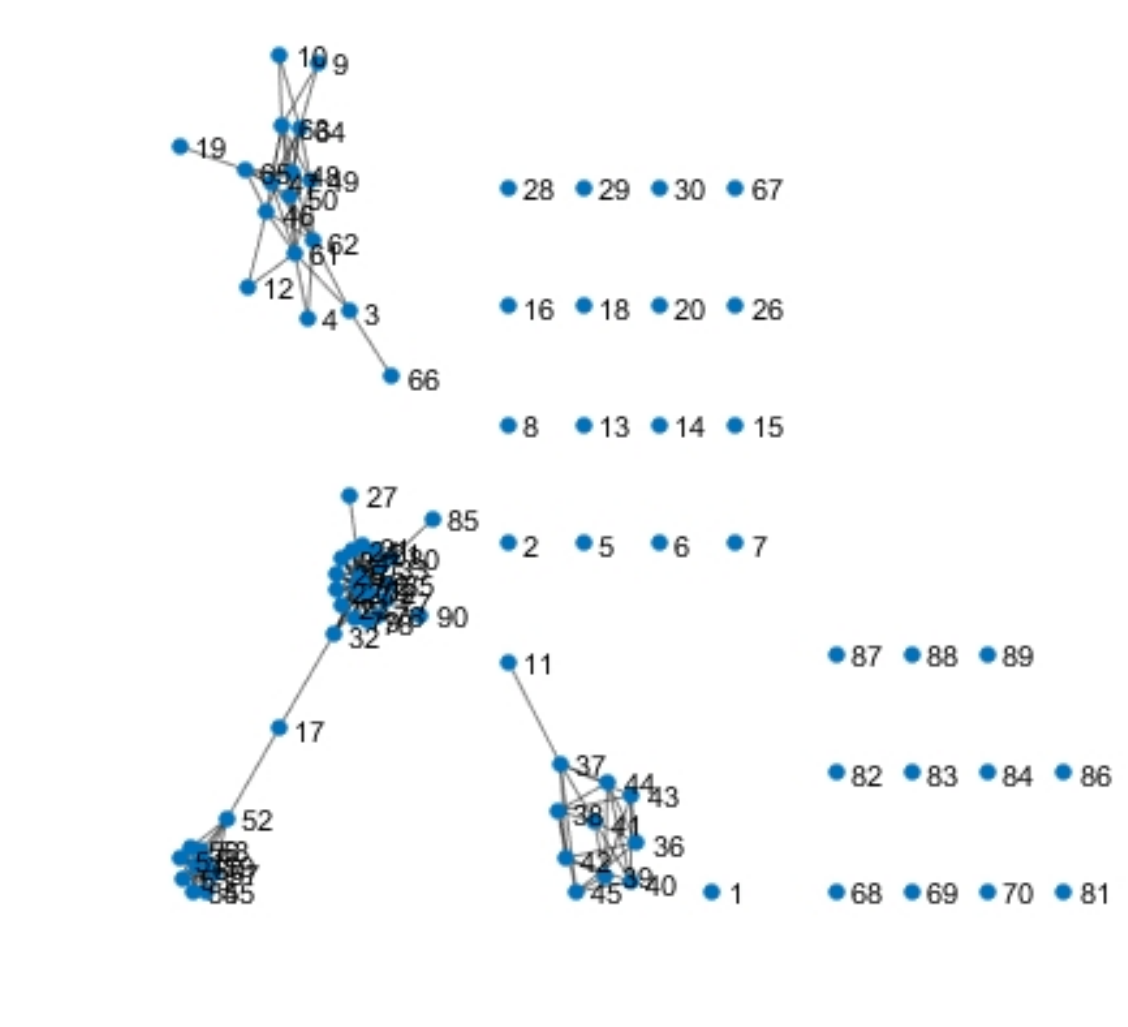}}
    \caption{Ground truth.}
  \end{subfigure}
 \begin{subfigure} {0.3\textwidth}
    {\includegraphics[width=\textwidth]{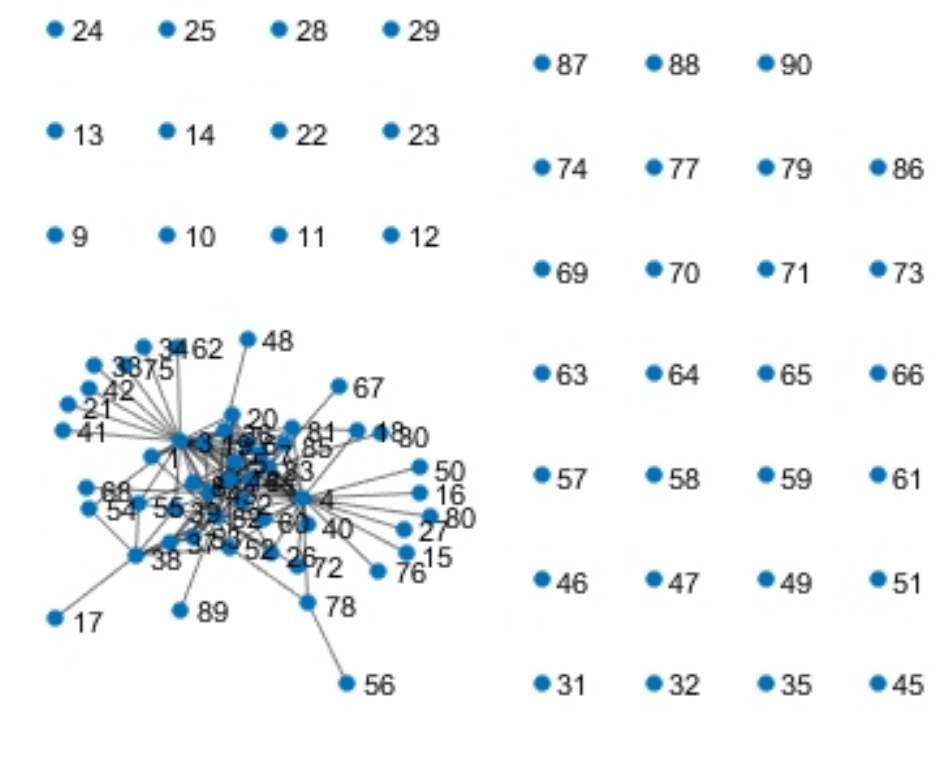}}
    \caption{Graphical lasso.}
  \end{subfigure}
  \hfill
  \begin{subfigure} {0.3\textwidth}
  {\includegraphics[width=\textwidth]{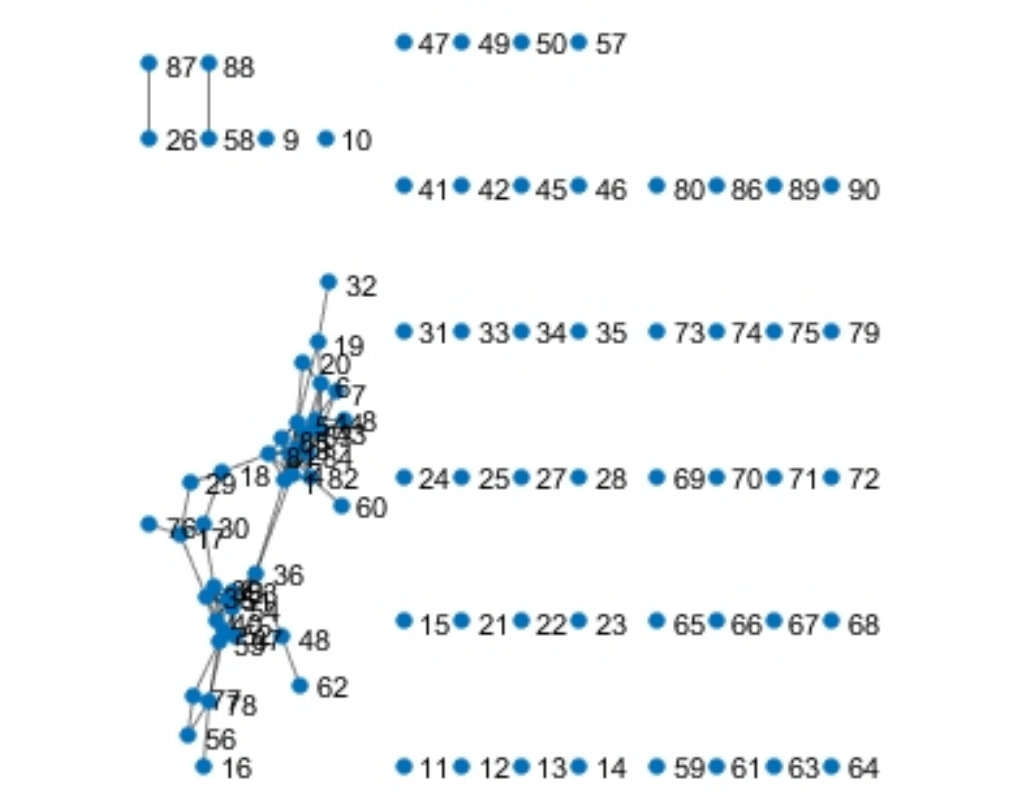}}
    \caption{SSONA.}
  \end{subfigure}
  \hfill
  \begin{subfigure} {0.3\textwidth}
  {\includegraphics[width=\textwidth]{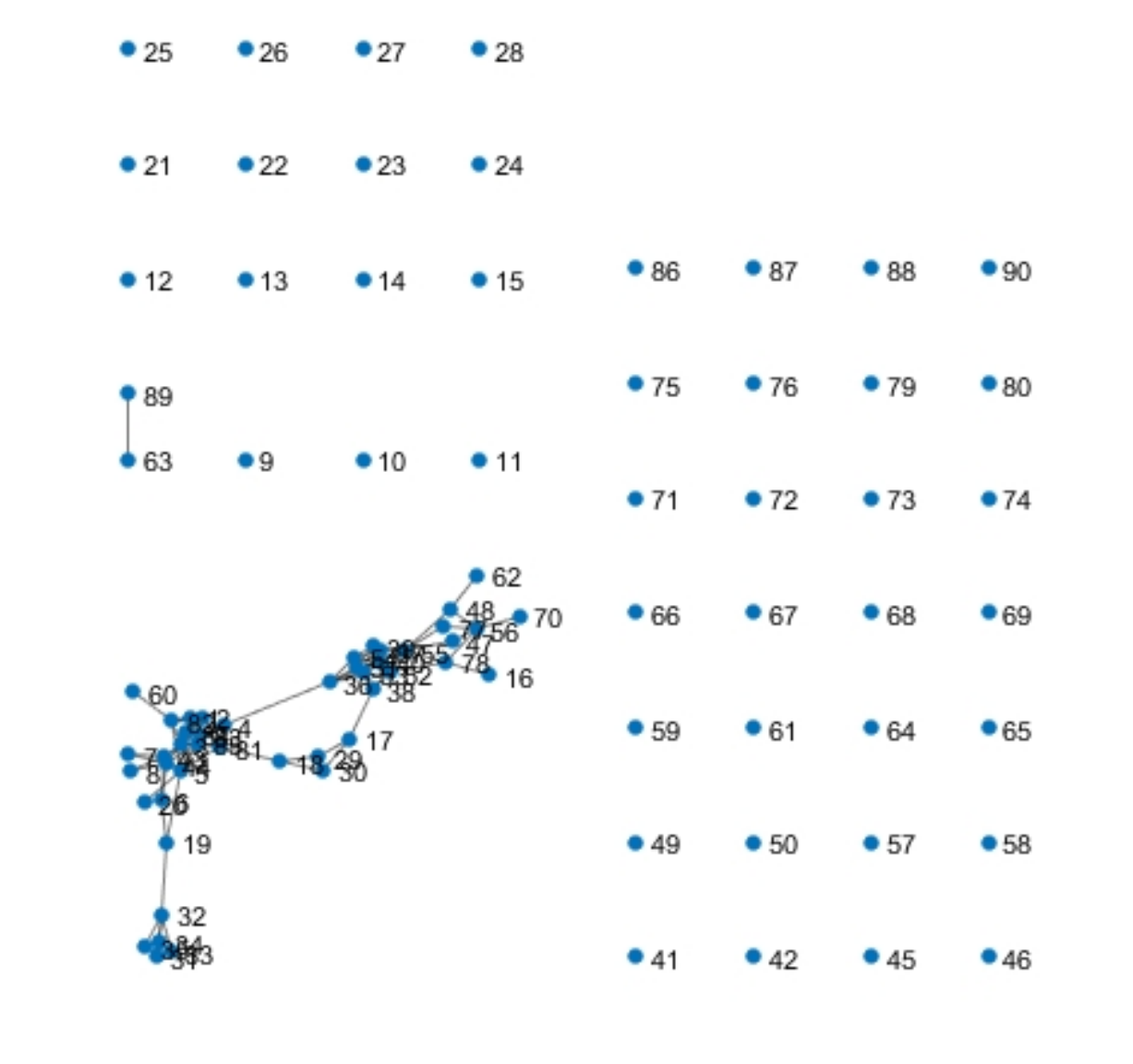}}
    \caption{Ground truth.} 
  \end{subfigure}
\caption{\small{Simulation for the Gaussian graphical model with 10 latent variables. Row I: Results for $p=100$ and $m=200$. Row II:
Results for $p=100$ and $m=100$.}}\label{figs3}
\end{figure}

\begin{figure}[!ht]
 \begin{subfigure} {0.3\textwidth}
    {\includegraphics[width=\textwidth]{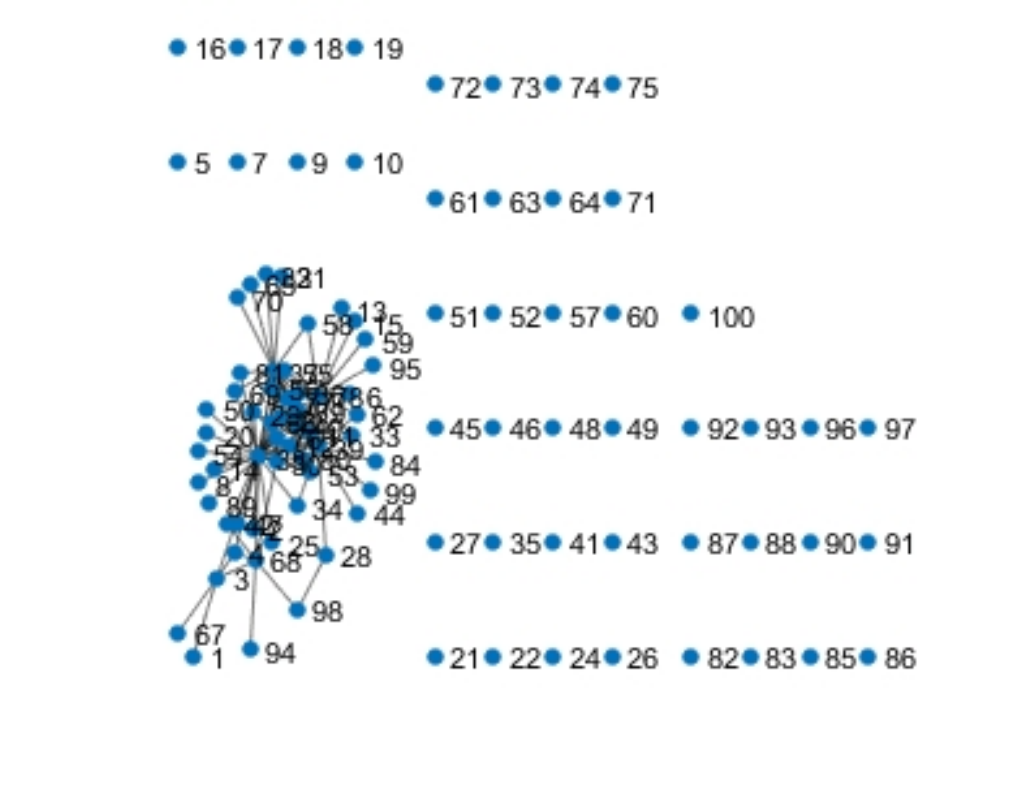}}
    \caption{Graphical lasso.}
  \end{subfigure}
  \hfill
  \begin{subfigure} {0.31\textwidth}
  {\includegraphics[width=\textwidth]{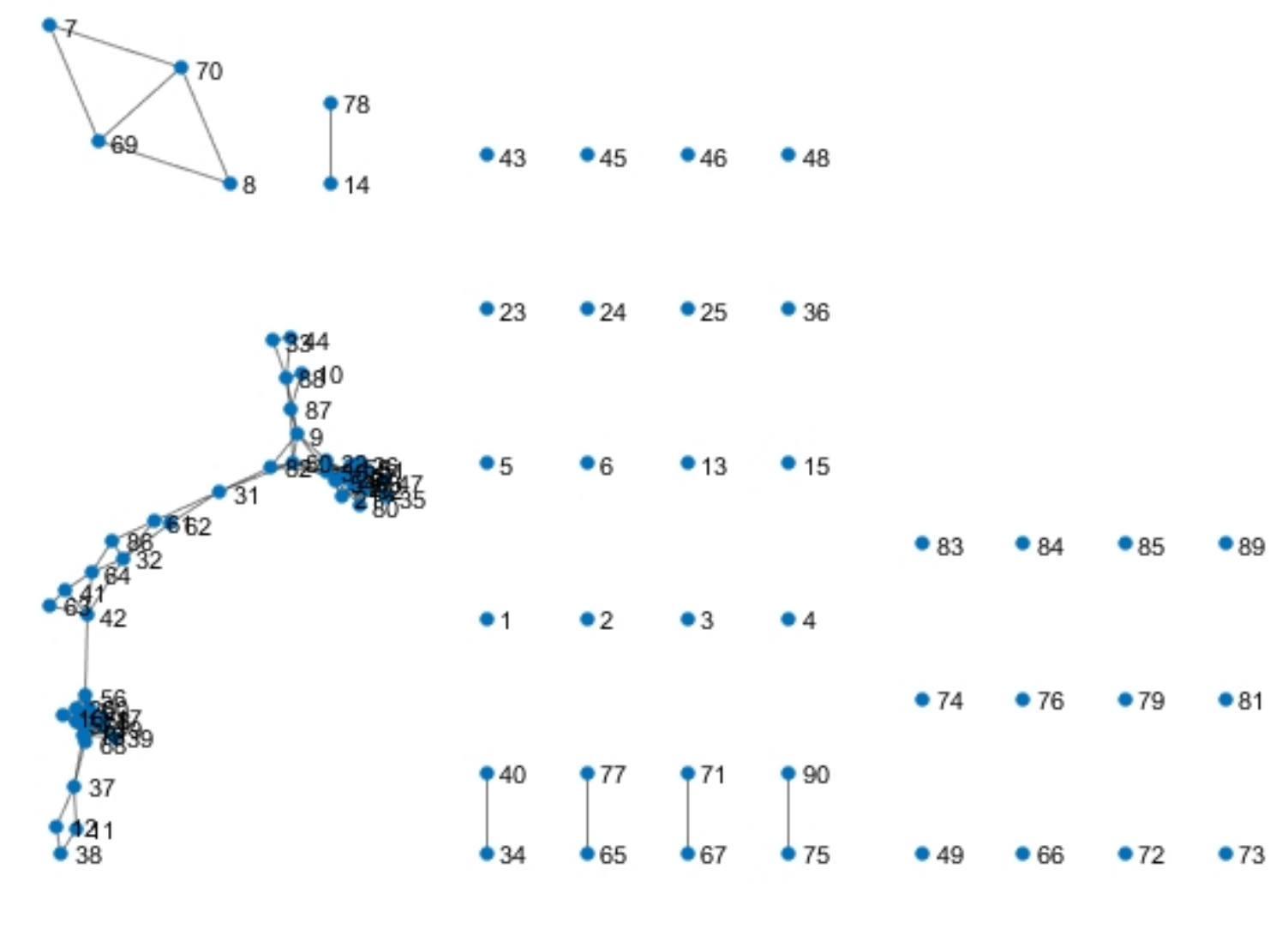}}
    \caption{SSONA.}
  \end{subfigure}
  \hfill
  \begin{subfigure} {0.26\textwidth}
  {\includegraphics[width=\textwidth]{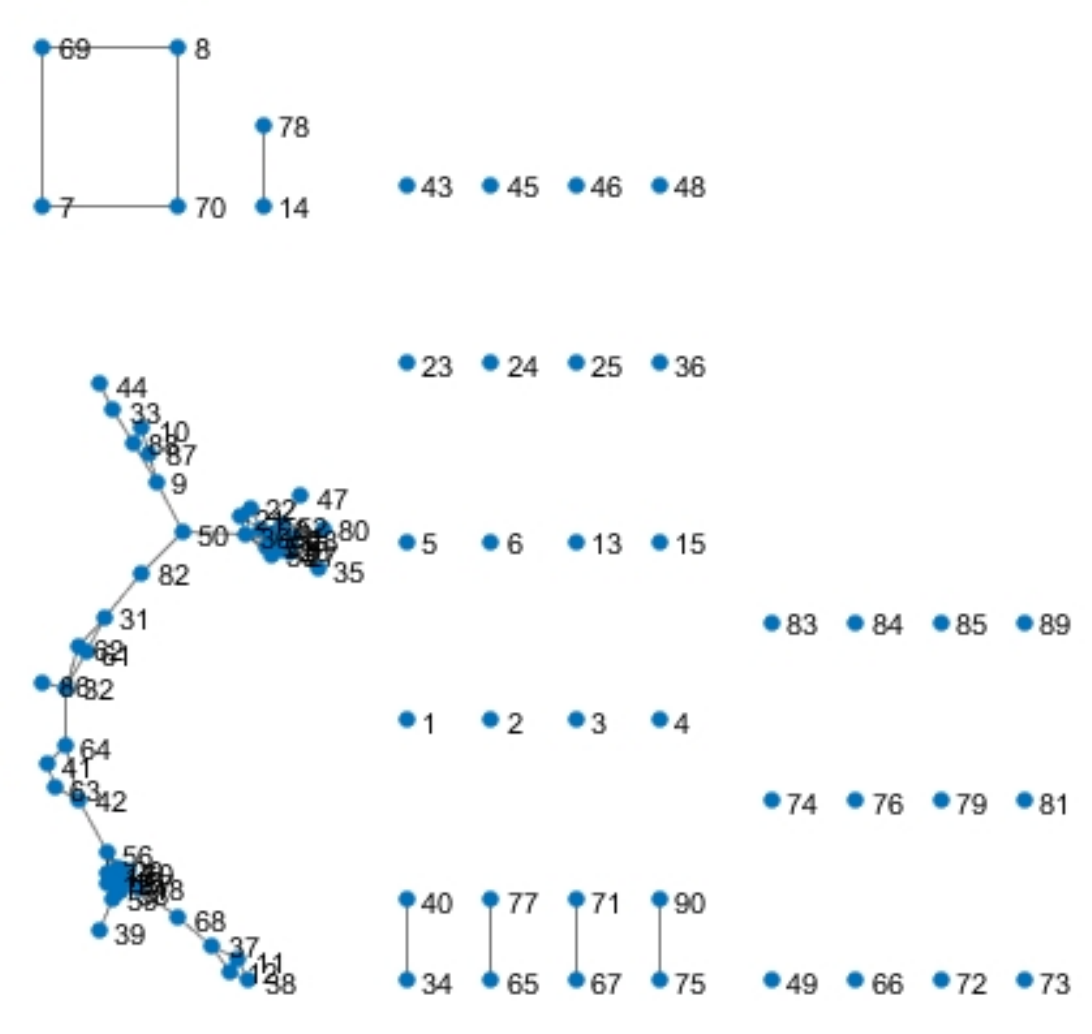}}
    \caption{Ground truth.}
  \end{subfigure}
 \begin{subfigure} {0.3\textwidth}
    {\includegraphics[width=\textwidth]{binP.pdf}}
    \caption{Graphical lasso.}
  \end{subfigure}
  \hfill
  \begin{subfigure} {0.31\textwidth}
  {\includegraphics[width=\textwidth] {binS.pdf}}
    \caption{SSONA.}
  \end{subfigure}
  \hfill
  \begin{subfigure} {0.26\textwidth}
  {\includegraphics[width=\textwidth]{binT.pdf}}
    \caption{Ground truth.} 
  \end{subfigure}
\caption{\small{Simulation for the binary Ising Markov random field. Row~I: Results for $p=100$ and $m=200$. Row~II:
Results for $p=100$ and $m=100$.}}\label{figs4}
\end{figure}

In this section, we present numerical results on structured graphical models to demonstrate the efficiency of SSONA. We compare the behavior of SSONA for a fixed value of $p=100$ with a lasso version of our algorithm. Results provided in Figures~ \ref{figs1}, \ref{figs2}, \ref{figs3} and \ref{figs4} indicate the efficiency of algorithm~\ref{alg:1} on structured graphical models. These results also show how the structure of the network returned by the two algorithms changes with growing $m$ (note that $\lambda_i$  and $\hat{\lambda}_i$ are kept fixed for each value of $m$). It can be easily seen from these figures (comparing Row~I and II) that SSONA is less sensitive to the number of samples and shows a better approximation of the network structure even for small sample sizes.

\subsection{Classification and clustering accuracy based on SSONA}\label{Exp1}

In this section, we evaluate the efficiency of SSONA on  real data sets in recovering complex structured sparsity patterns and subsequently evaluate them on a classification task. The two data sets deal with applications in cancer genomic and document classification.

\subsubsection{SSONA for Gene Selection Task}

Classification with a sparsity constraint has become a standard tool in applications involving Omics data, due to the large number of available features
and the small number of samples. The data set under study considers gene expression profiles of lung cancer tumors. Specifically, the data\footnote{
\texttt{http://www.broadinstitute.org/cgibin/cancer/publications/view/87}.
} consist of gene expression profiles of 12,626 genes for 197 lung tissue
samples, with 139 adenocarcinomas(AD),
21 squamous cell carcinomas(SQ), 20 carcinoids
(COID) and 17 normal lung tissue (NL). To distinguish lung adenocarcinomas from the normal lung tissues, we consider the diagnosis of lung cancer as a binary classification problem. Let the 17 normal lung comprise the positive class and the 139 lung adenocarcinomas the  negative class. Following the workflow in \citet{monti03}, we reserve the 1000 most significant genes after a preprocessing step. In the numerical experiment, we compare group lasso \citep{Yuan07}, group lasso with overlap \citep{Obo11} and SSONA according to the following two criteria: average classification accuracy and gene selection performance. The experiment is repeated ten times and the average accuracy and performance are depicted in Table~\ref{tab:clas}.
 \begin{table}[H]
\begin{center}
\begin{tabular}{ >{\arraybackslash}m{3.5in}  >{\centering\arraybackslash}m{.95in} >{\centering\arraybackslash}m{.95in} >{\centering\arraybackslash}m{.95in} >{\centering\arraybackslash}m{.95in} >{\centering\arraybackslash}m{.95in}}
\toprule[1.75pt]
{\bf Method} & {\bf Average classification accuracy } & {\bf Average number of genes selected} \\
\midrule
Group lasso \citep{Yuan07}  &  0.815(0.046) &  69.11(3.23) \\
Group lasso with overlap \citep{Obo11}  &   0.834(0.035) &  57.30(2.71)\\
SSONA (4 structured matrices) & 0.807(0.028) &  61.44(2.80)\\
SSONA (6 structured matrices) & 0.839(0.022) &  56.111(2.100)\\
\bottomrule[1.25pt]
\end {tabular}
\end{center}
\caption {Experimental results on lung cancer data over 10 replications (the standard deviations are reported in parentheses).
} \label{tab:clas}
\end{table}

As is shown in Table~\ref{tab:clas}, SSONA achieves higher classification accuracy than the group lasso and lower classification accuracy than the latent group lasso, although the performance of all three methods is very similar and within the variability induced by the replicates.
However, our SSON based lasso does not require a priori knowledge of group structures, which is a prerequisite for the other two methods. One can easily improve the the classification accuracy and gene selection performance of SSONA by adding more structured matrices.

In our experiments, SSONA selects the least number of  genes and achieves the smallest standard deviation of average number of genes {\em without any priori knowledge}. Due to the different number of randomly selected genes, the average number of gene sometimes will be a non-integer.

\subsubsection{SSONA for Document Classification Task}

The next example involves a data set
 \footnote{\texttt{http://qwone.com/jason/
20Newsgroups/}} containing 1427 documents with a corpus of size 17785 words. We randomly partition the data into 999 training, 214 validation and 214 test examples, corresponding to a 70/15/15 split \citep{Rao16}. We first train a Latent Dirichlet Allocation based topics  model \citep{blei03} to assign the words to 100 "topics". These correspond to our groups, and since a single word can be assigned to multiple topics, the groups overlap. We then train a  lasso
logistic model using as outcome variable indicating whether the document discusses atheism or not , together with an overlapping group lasso and a SSON
based lasso model where the tuning parameters are selected based on cross validation. Table \ref{tab222} shows that the variants of the SSON yield almost the same misclassification rate compared to the other two methods,  while it does not require a priori knowledge of group structures.

\begin{table}[H]
\begin{center}
\begin{tabular}
{>{\arraybackslash}m{3.5in}  >{\centering\arraybackslash}
>{\centering\arraybackslash}m{2in}
>{\centering\arraybackslash}
>{\centering\arraybackslash}m{2in} >{\centering\arraybackslash}m{2in}}
\toprule[2pt]
{\bf Method} & {\bf Misclassification Rate} \\
\midrule
Group lasso \citep{Yuan07}  &  0.445 \\
Group lasso with overlap \citep{Obo11}  &  0.390\\
SSONA (5 structured matrices) & 0.435\\
SSONA (6 structured matrices) & 0.421\\
SSONA (7 structured matrices) & 0.401\\
\bottomrule[2pt]
\end {tabular}
\caption{Misclassification rate on the test set for document classification.} \label{tab222}
\end{center}
\end{table}

\subsubsection{SSONA for structured subspace clustering}

Our last example focuses on data clustering. The data come from multiple low-dimensional linear or affine subspaces embedded in a high-dimensional space. Our method is based on \eqref{eq:7}, wherein each point in a union of subspaces has a representation with respect to a dictionary formed by all other data points. In general, finding such a
representation is NP hard.  We apply our subspace clustering algorithm to a structured data in the presence of noise. The segmentation of the data is obtained by applying SSONA to the adjacency matrix built from the data. Our method can handle noise and missing data and is effective to detect the clusters.

\begin{figure}[!ht]
\captionsetup[subfigure]{labelformat=empty}
\begin{subfigure} {0.4\textwidth}
    {\includegraphics[width=\textwidth]{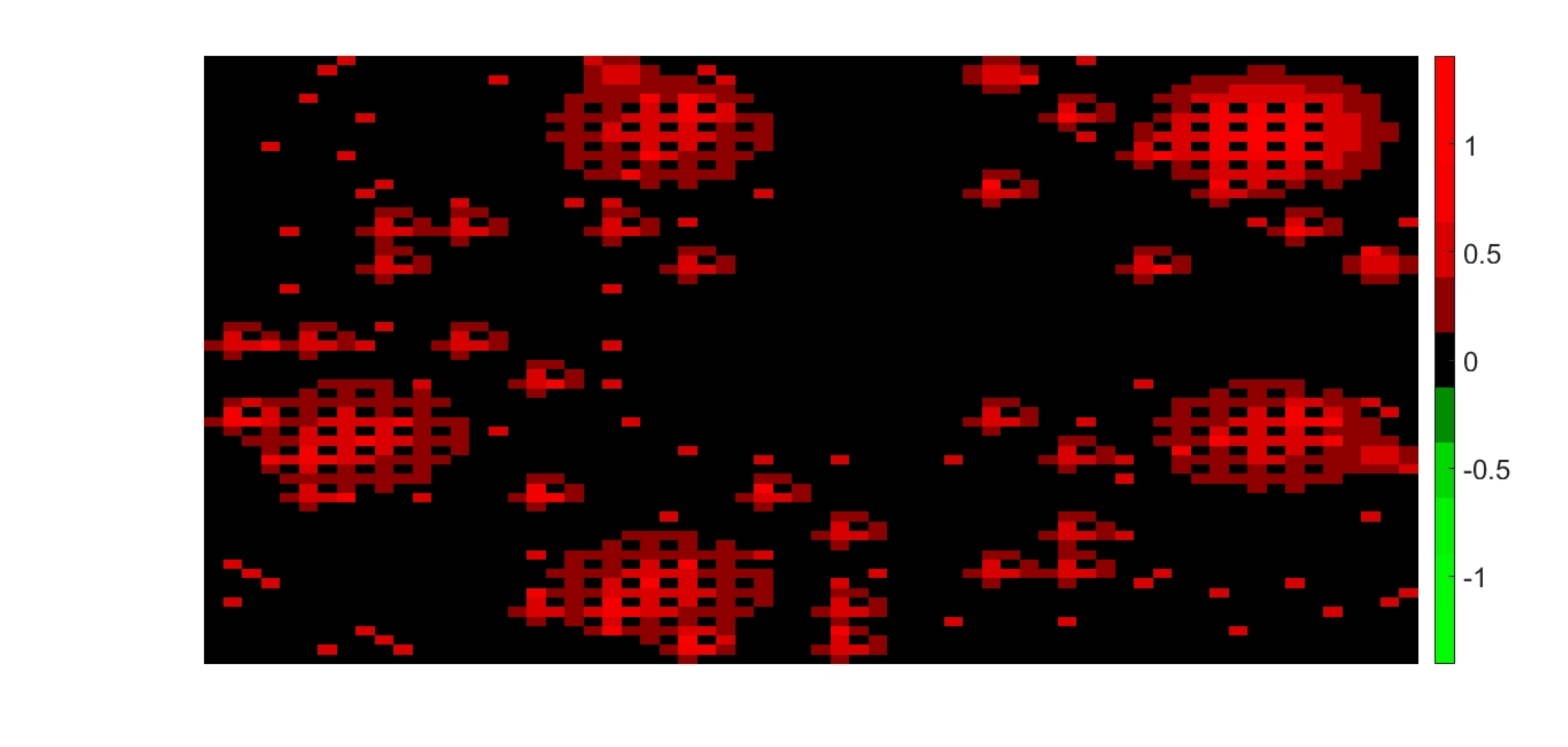}}
    \caption{Ground Truth.}
  \end{subfigure}
  \hfill
  \begin{subfigure} {0.4\textwidth}
  {\includegraphics[width=\textwidth]{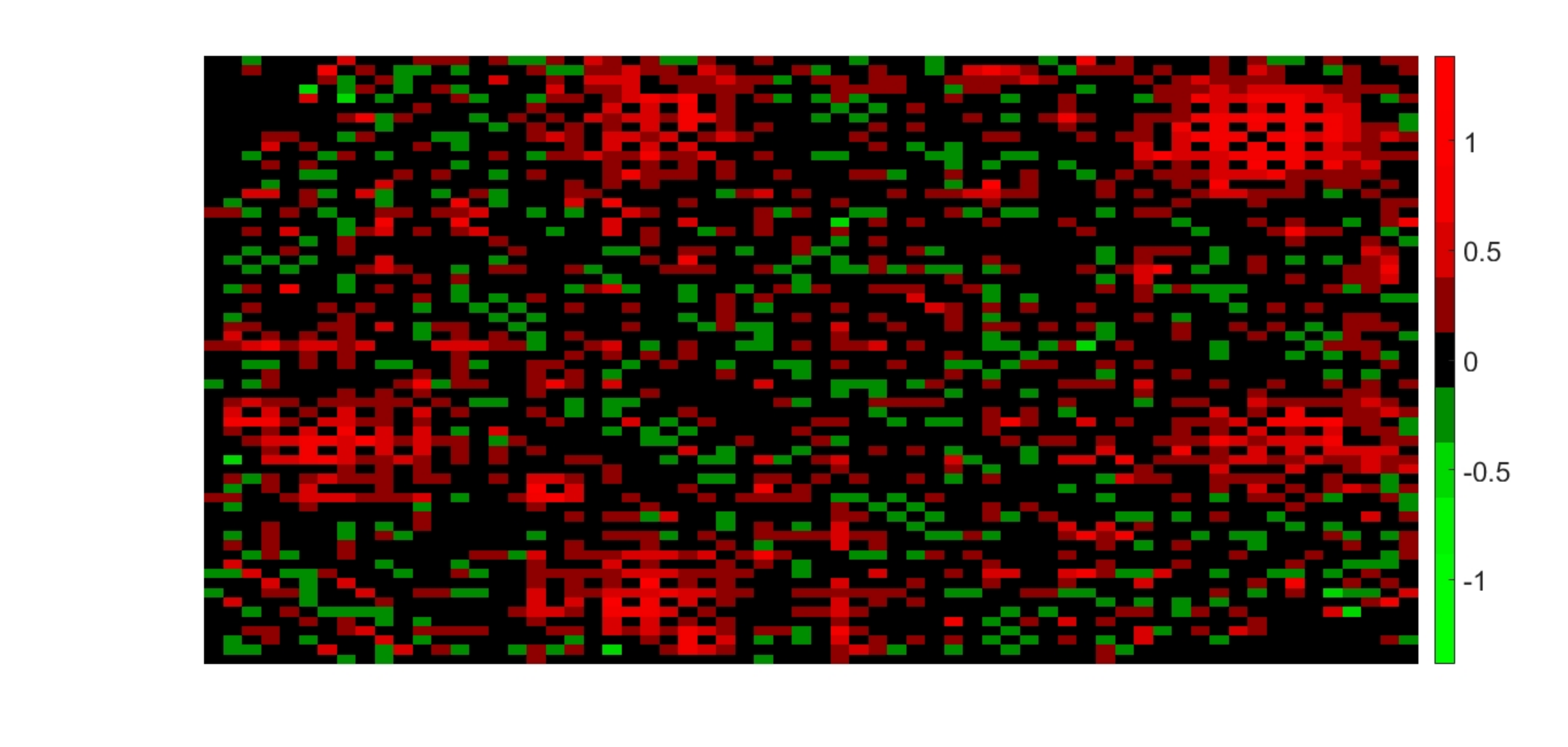}}
    \caption{Ground Truth+ Noise.}
    \end{subfigure}
  \hfill
  \begin{subfigure} {0.4\textwidth}
  {\includegraphics[width=\textwidth]{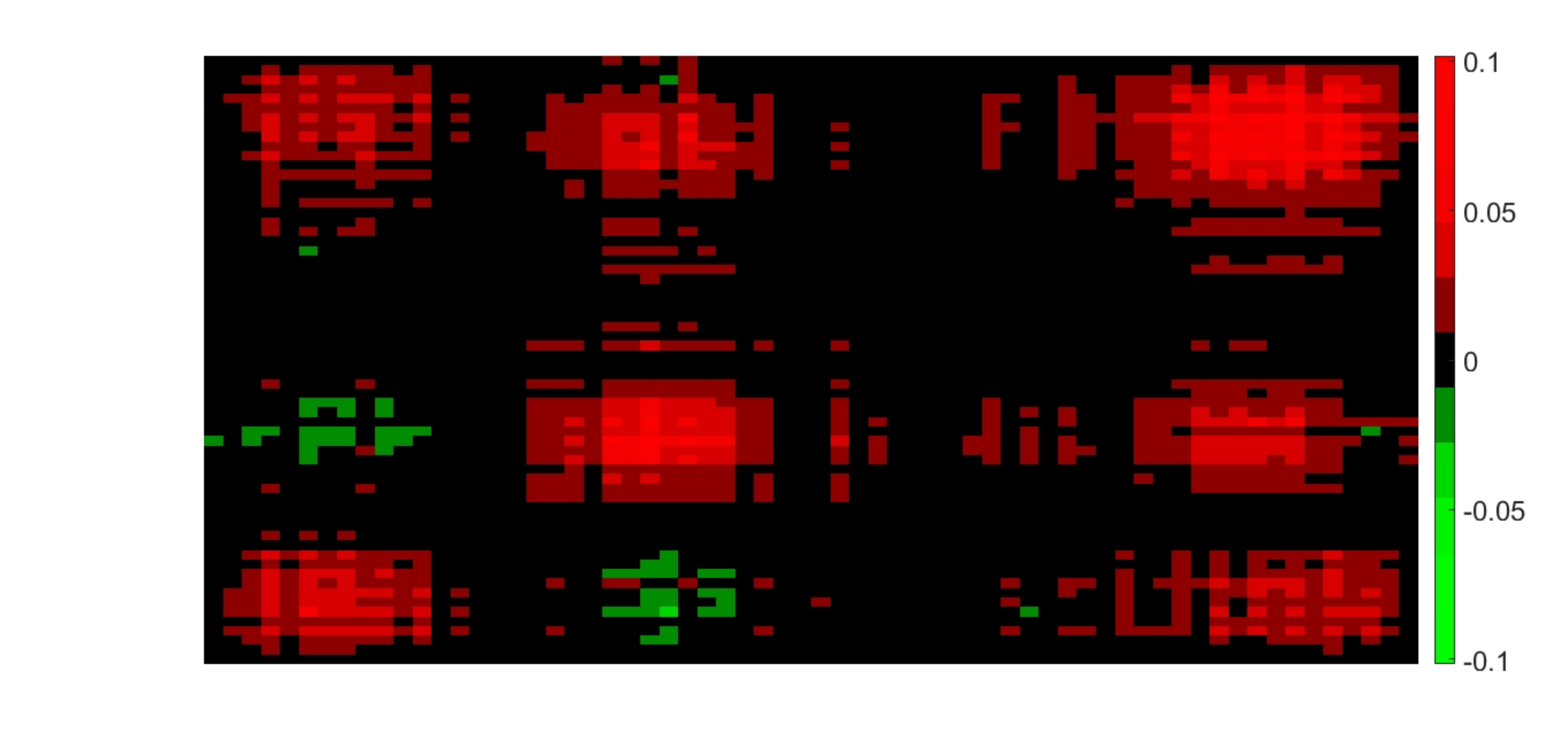}}
    \caption{LRR \citep{liu13}.}
  \end{subfigure}
 \hfill
 \begin{subfigure} {0.4\textwidth}
    {\includegraphics[width=\textwidth]{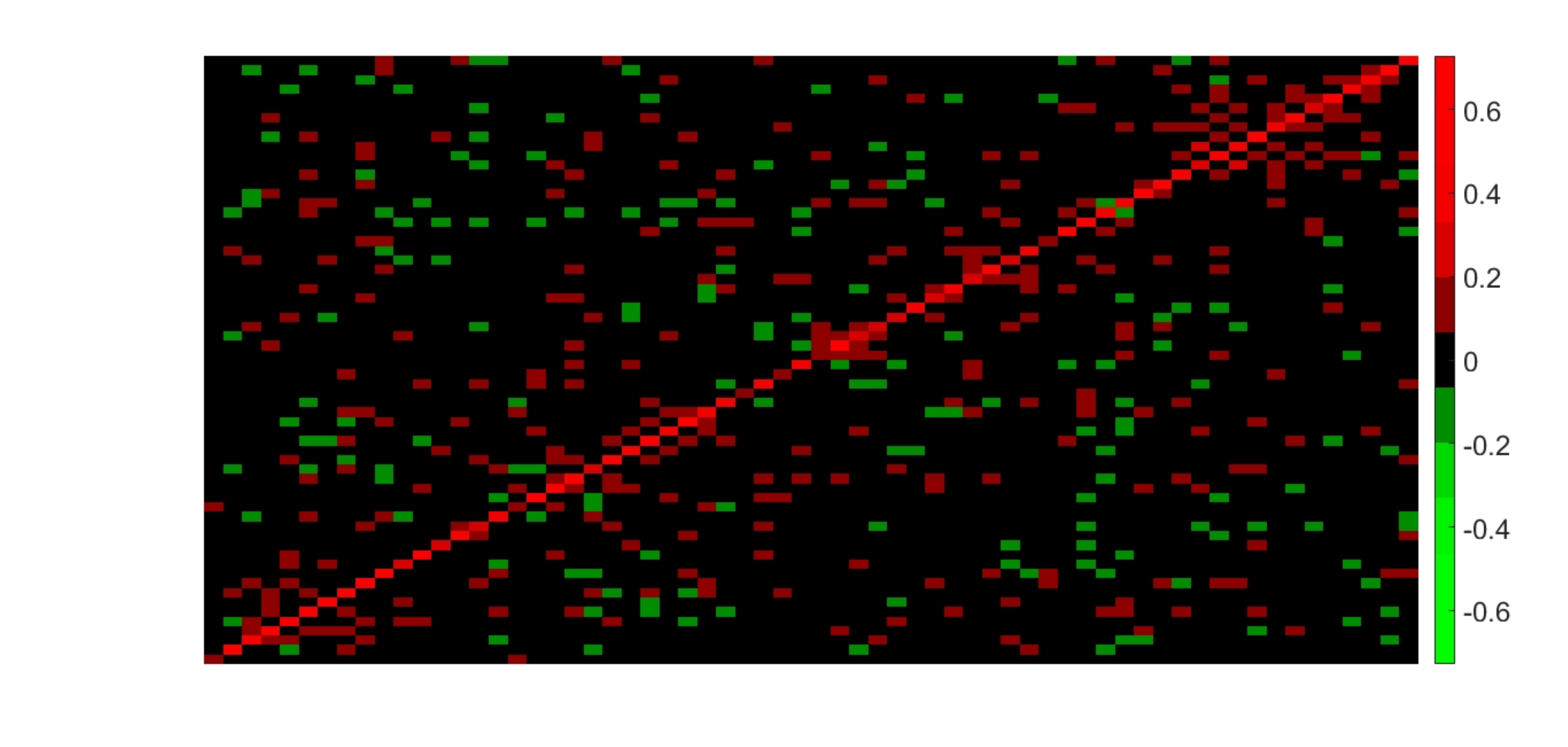}}
    \caption{SSC \citep{elham09}.}
  \end{subfigure}
  \hfill
  \begin{center}
  \begin{subfigure} {.5\textwidth}
  {\includegraphics[width=\textwidth] {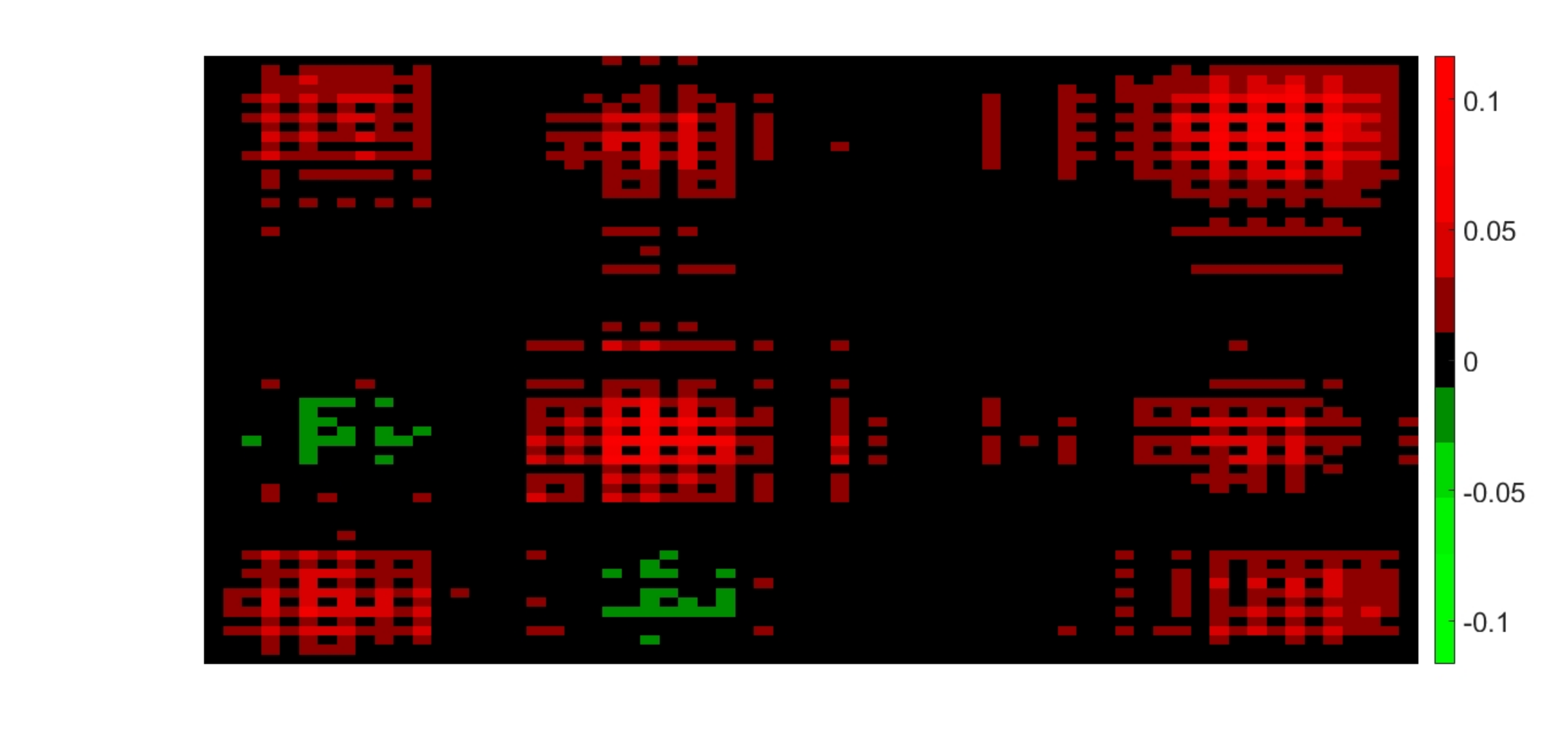}}
    \caption{SSONA.}
  \end{subfigure}
  \end{center}
\caption{Heatmap of different algorithms for detecting clusters in data.}
\label{figsub}
\end{figure}

Figure~\ref{figsub} shows that our approach significantly outperforms state-of-the-art methods.

\subsection{Application to real data sets}\label{realdata}

Next, we use the SSON framework to analyze three data sets from molecular and social science domains. Although there is no known ground truth, the proposed framework recovers interesting patterns and highly interpretable structures.

\paragraph{Analysis of connectivity in the financial sector.}

We applied the SSON methodology to analyze connectivity in the financial sector. We use monthly stock returns data from August, 2001 to July, 2016 for three financial sectors, namely banks (BA), primary broker/dealers (PB), and insurance companies (INS). The data are obtained from the University of Chicago's Center for Research in Security Prices database (CRSP).

\begin{figure}
\centering
   \includegraphics[width=1\textwidth]{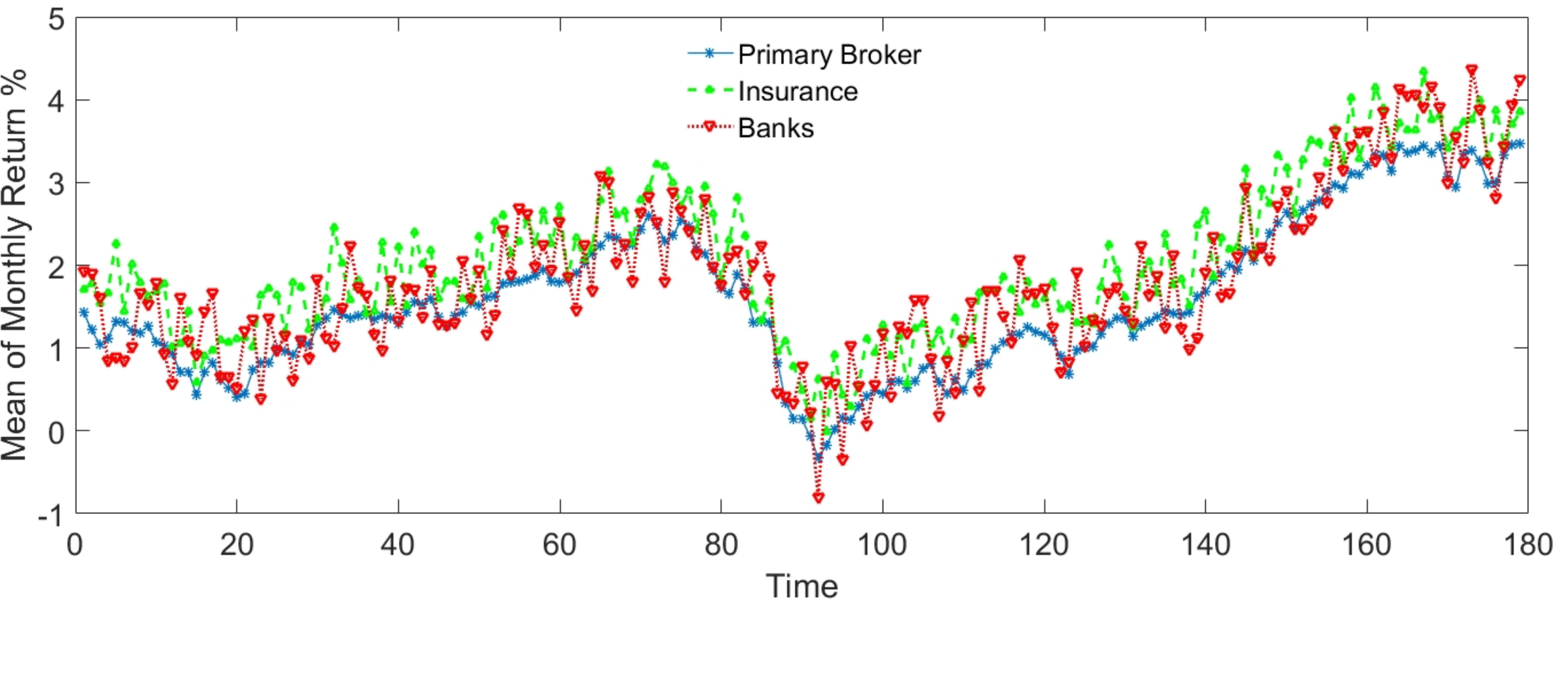}
 \caption{Average monthly return of firms in the three sectors- Bank, primary broker-dealer and insurance firms, in different 3-year rolling windows during 180 months. The figure shows diminished performance during the 2007-2009 crisis (time step : 80-100) and also clearly captures the strong recovery
of stock performance starting in 2009.}
 \label{mean}
\end{figure}

Our final sample covers 75 different institutions spanning a 16-year period. Figure~\ref{mean} shows the mean (in \%) of monthly stock returns across different sectors in each 3-year long rolling windows. As expected, the average returns are significantly lower during the financial 2007-2009 crisis period, compared to any other period in our sample. Indeed, looking across the sectors, all three sectors experienced diminished performance during the 2007-2009 crisis. Further, the almost linear ramp-up following 2009 clearly captures the recovery of financial stocks and the broader market.

\begin{figure}[h]
 \begin{subfigure} {0.85 \textwidth}
 \includegraphics[width=.85\textwidth]{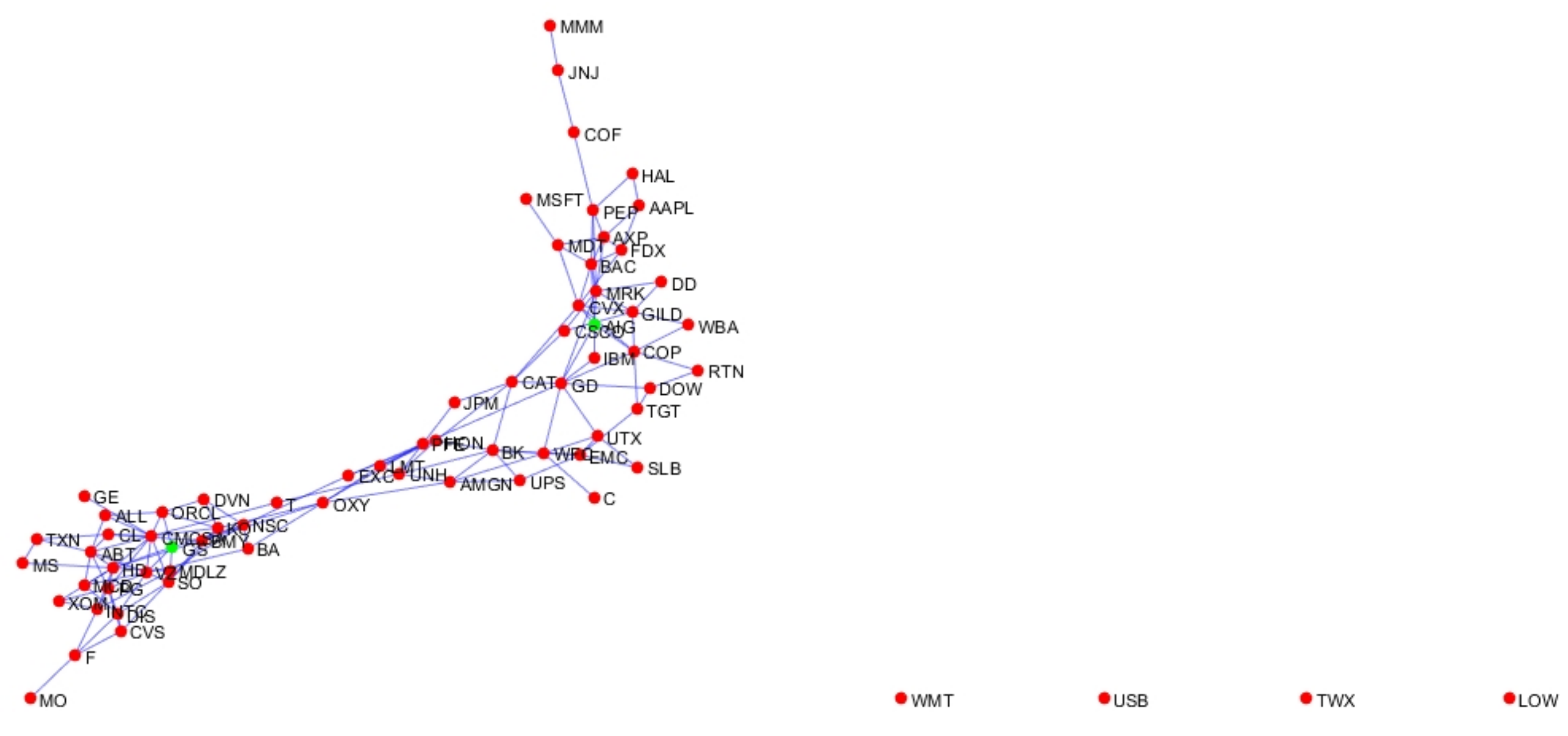}
    \caption{SSONA VAR}
  \end{subfigure}
  \hfill
  \begin{subfigure} {0.85 \textwidth}
\includegraphics[width=.85\textwidth]{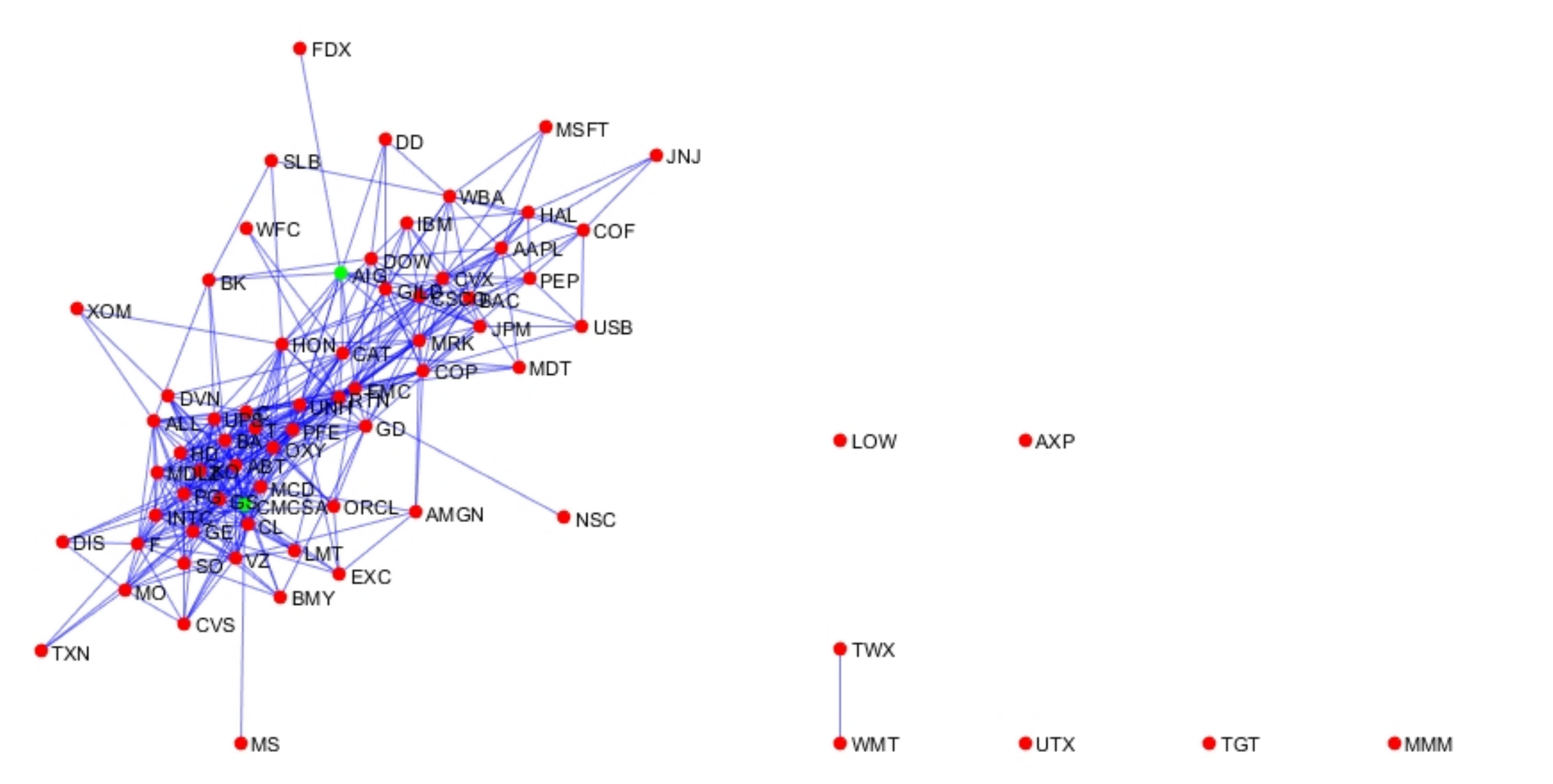}
    \caption{Lasso VAR}
    \end{subfigure}
 \caption{Networks estimate by SSONA and Lasso VAR during crisis period of Jan 2007 to Oct 2009.} \label{crisnet}
\end{figure}

Next, we estimate a measure of network connectivity for a sample of the 71 components of the SP100 index that were present during the entire 2001-16 period under consideration.
Figure \ref{crisnet} depicts the network estimates of the transition (lead-lag) matrices using straight lasso VAR and SSONA based VAR for the January 2007 to Oct 2009 period. It can be seen that the lasso VAR estimates produce a more highly connected network, while the SSONA ones identify
two more connected components. Both methods highlight the key role played by AIG and GS (Goldman Sachs), but the SSONA based network indicates that one dense connected component is centered around the former, while the other dense connected component around the latter.
In summary, both methods capture the main connectivity patterns during the crisis period, but SSONA provides a more nuanced picture.

\paragraph{US House voting data set.}
We applied SSONA to describe the relationships amongst House Representatives in the U.S. Congress during the 2005-2006 period (109th Congress). The variables correspond to the 435 representatives, and the observations to the 1210 votes that the House deliberated and voted on during that period, which include bills, resolutions, motions, debates and roll call votes. The assumption of our model is that bills are i.i.d. sample from the same underlying Ising model. The votes are recorded as "yes" (encoded as "1") and "no" (encoded as "0"). Missing observations were replaced with the majority vote of the House member’s party on that particular vote.
\begin{figure}[!ht]
\captionsetup[subfigure]{labelformat=empty}
 \begin{subfigure} {0.3\textwidth}
 \includegraphics[width=\textwidth]{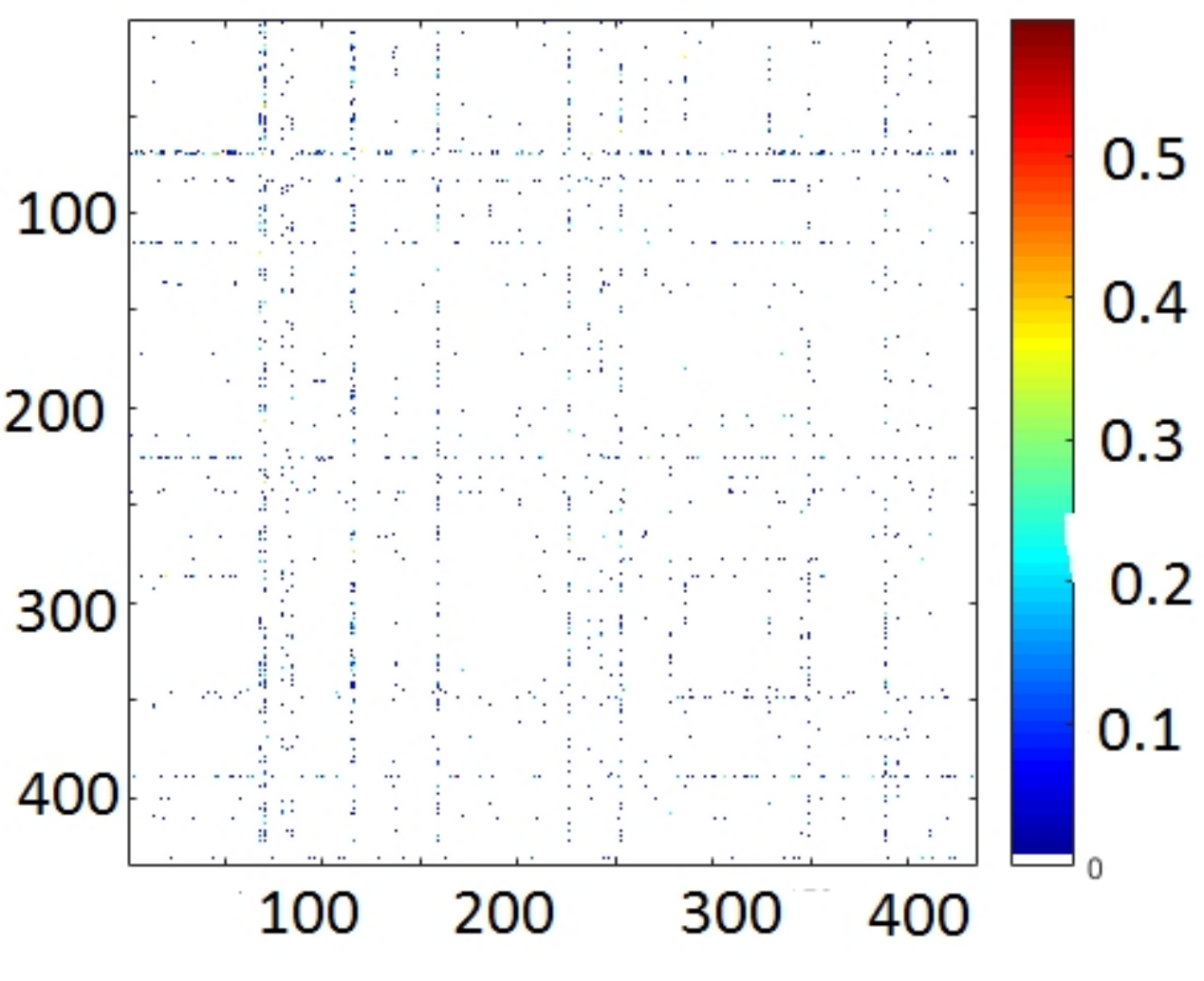}
    \caption{$Z_1+Z_1^\top$}
  \end{subfigure}
  \hfill
  \begin{subfigure} {0.3\textwidth}
   \includegraphics[width=\textwidth]{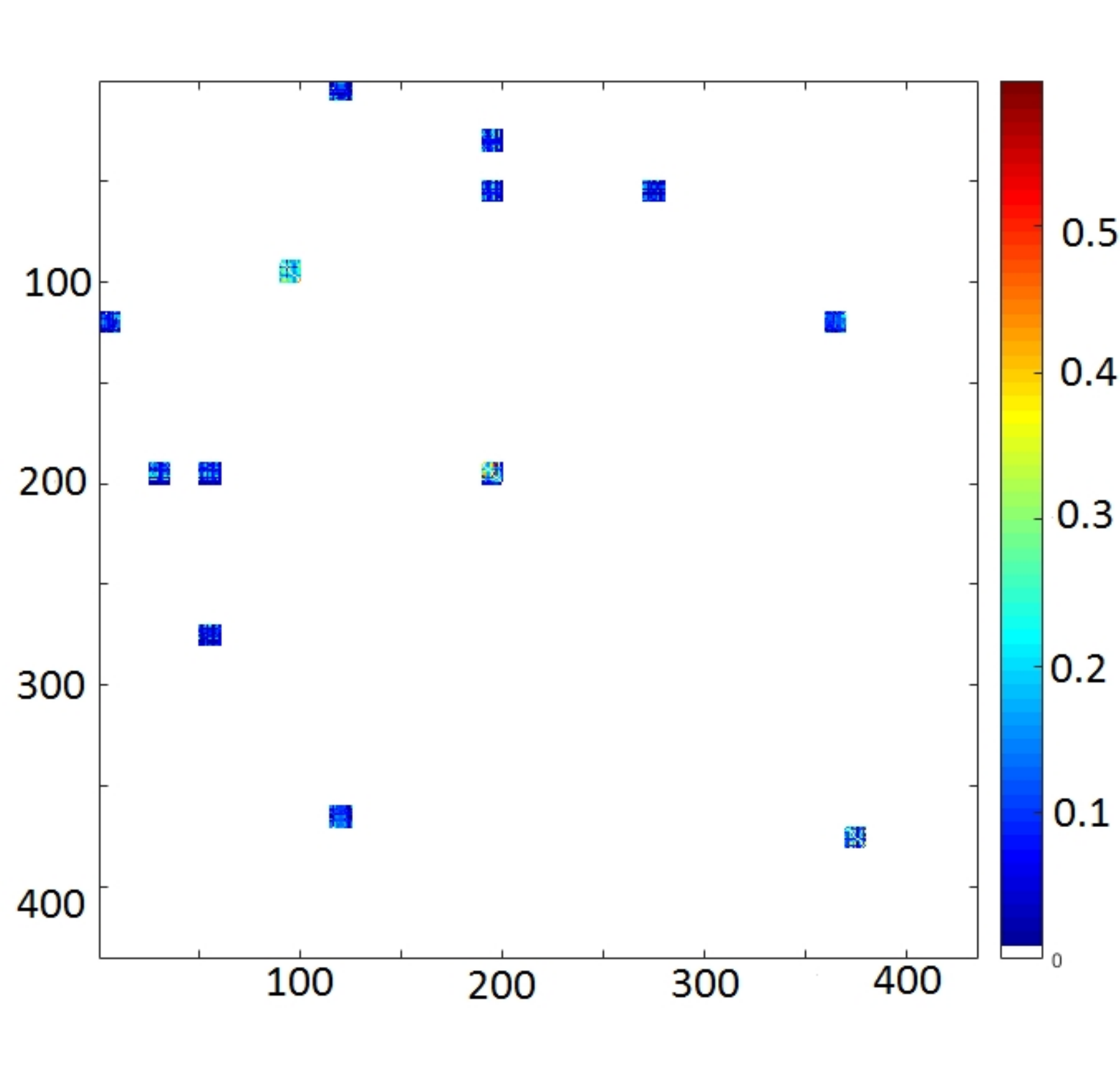}
    \caption{$Z_2+Z_2^\top$}
    \end{subfigure}
  \hfill
  \begin{subfigure} {0.3\textwidth}
 \includegraphics[width=\textwidth]{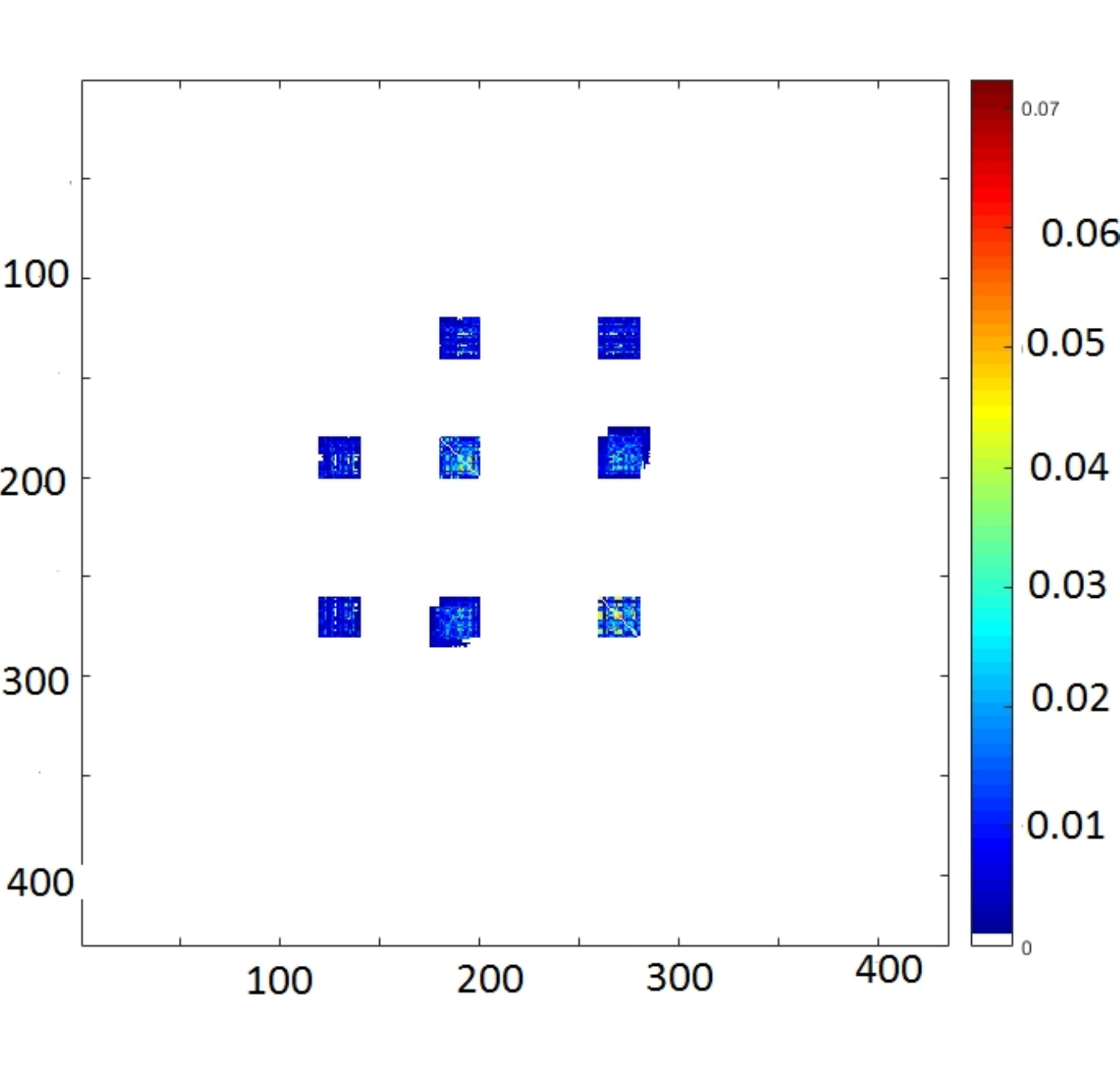}
    \caption{$Z_3+Z_3^\top$}
  \end{subfigure}
 \begin{subfigure} {0.3\textwidth}
 \includegraphics[width=\textwidth]{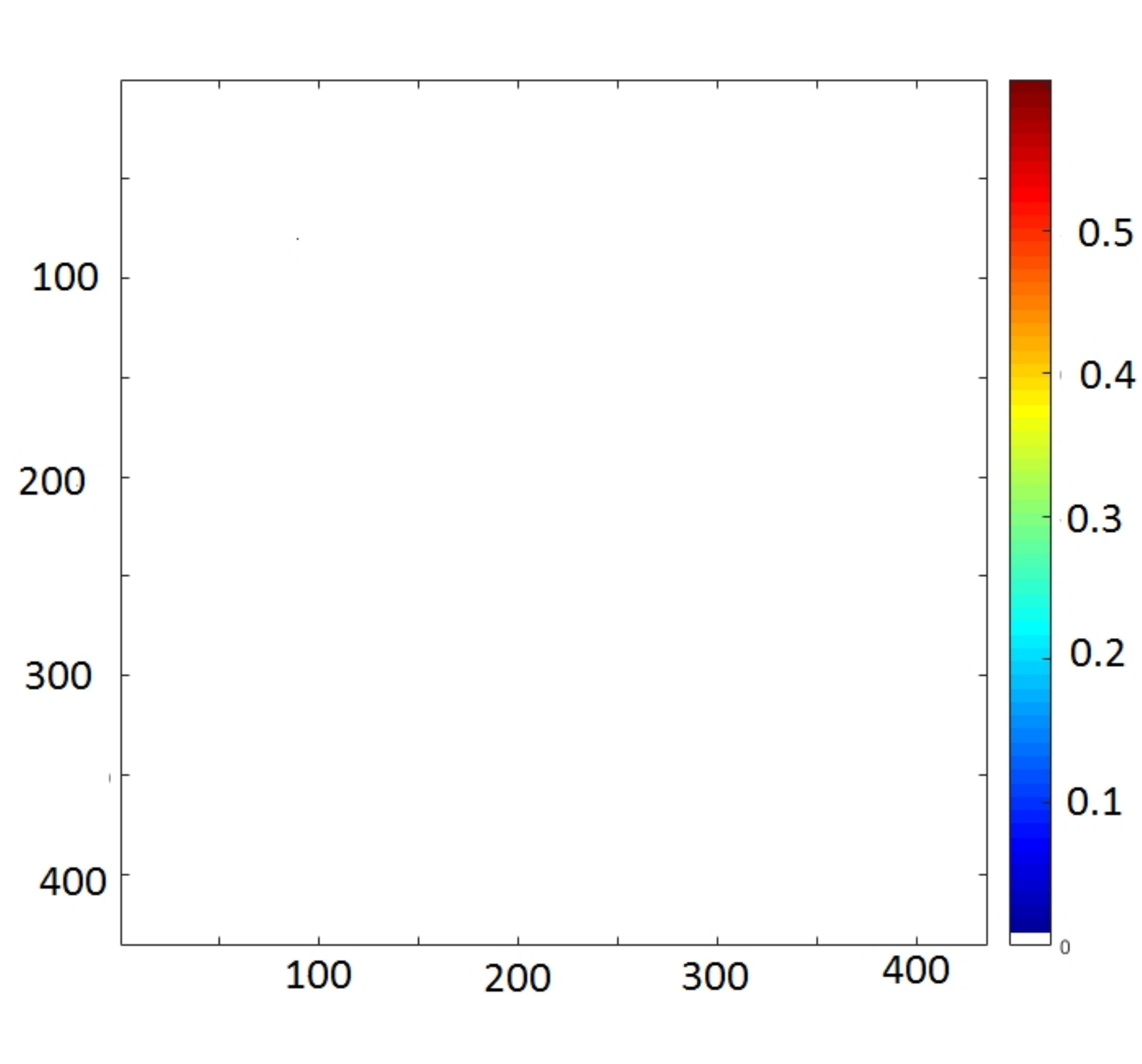}
    \caption{$Z_4+Z_4^\top$}
  \end{subfigure}
  \hfill
  \begin{subfigure} {0.3\textwidth}
 \includegraphics[width=\textwidth]{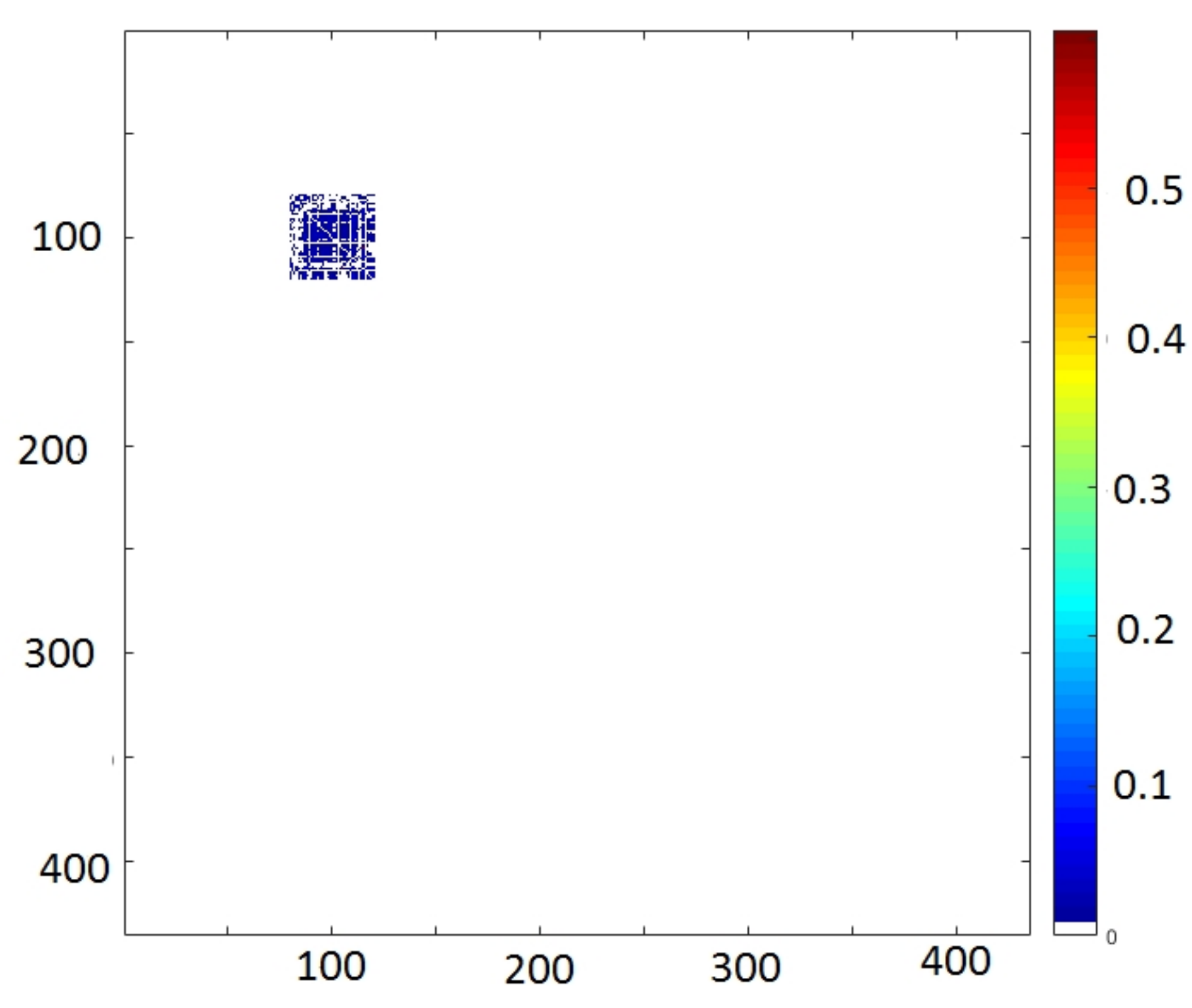}
    \caption{$Z_5+Z_5^\top$}
  \end{subfigure}
    \hfill
  \begin{subfigure} {0.3\textwidth}
  \includegraphics[width=\textwidth]{vottt.pdf}
    \caption{$Z_6+Z_6^\top$}
  \end{subfigure}
\caption{Heatmap of the structured precision matrix $\Theta$ decomposed into $Z_1+Z_1^\top +\dots+Z_6+Z_6^\top$ in the House voting data, estimated by SSONA.} \label{figvot1}
\end{figure}
Following \citet{Guo15}, we used a bootstrap procedure with the proposed SSONA estimator to evaluate the confidence of the estimated edges. Specifically, we estimated the network for multiple bootstrap samples of the same size, and only retained the edges that appeared more that $\omega$ percent of the time. The goal of the analysis is to understand the type of relationships that existed among the House members in the 109th Congress. In particular, we wish to identify and interpret the presence of densely connected components, as well of sparse components. The heatmap of the adjacency matrix of the estimated network by using SSONA is depicted in Figure \ref{figvot1}.
\begin{figure}[!ht]
  \centering
    \includegraphics[width=.4\textwidth]{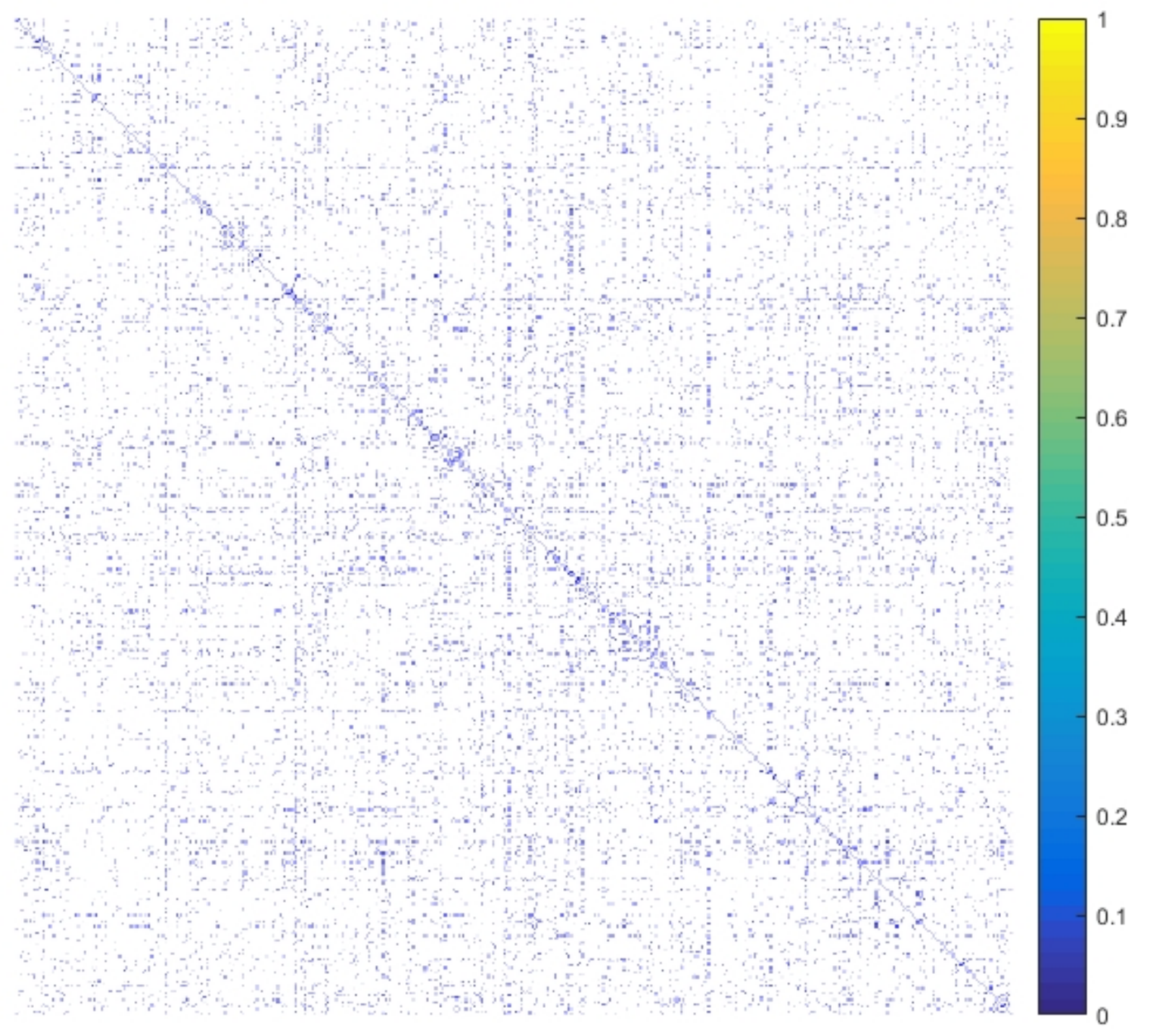}
\caption{Heatmap of the inverse covariance matrix in the voting record of the U.S. House of Representatives, estimated by the graphical lasso method \citep{Friedman07}.}
\label{figvote00}
\end{figure}
It can be easily seen that there exist densely connected components in the network, a fact that the glasso algorithm \citep{Friedman07} fails to recover (see, Figure~\ref{figvote00}).

\begin{figure}[!ht]
\captionsetup[subfigure]{labelformat=empty}
 \begin{subfigure} {0.3\textwidth}
   \includegraphics[width=\textwidth]{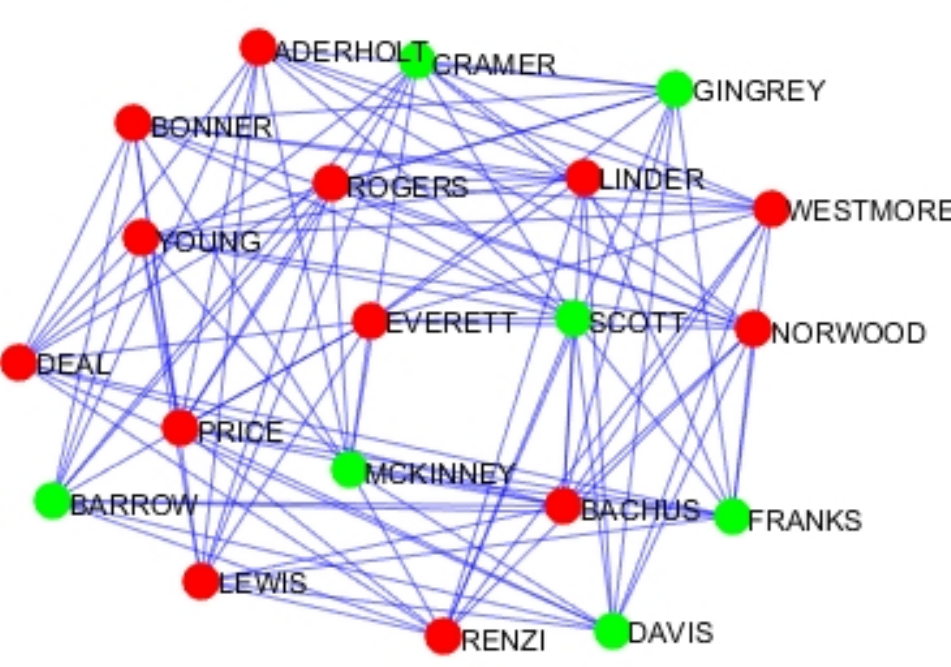}
    \caption{}
  \end{subfigure}
  \hfill
  \begin{subfigure} {0.3\textwidth}
    \includegraphics[width=\textwidth]{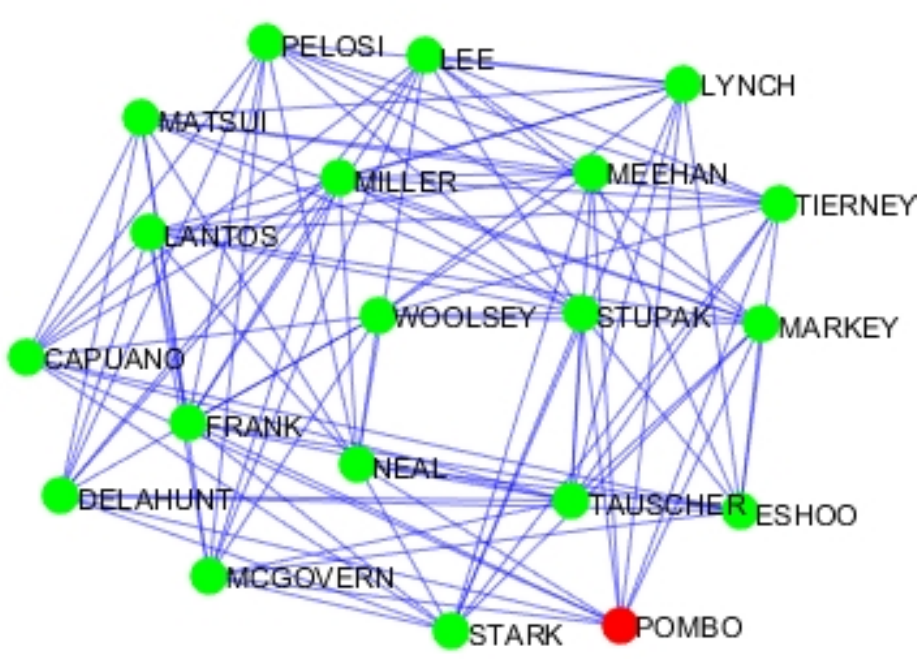}
    \caption{}
    \end{subfigure}
  \hfill
  \begin{subfigure} {0.3\textwidth}
    \includegraphics[width=\textwidth]{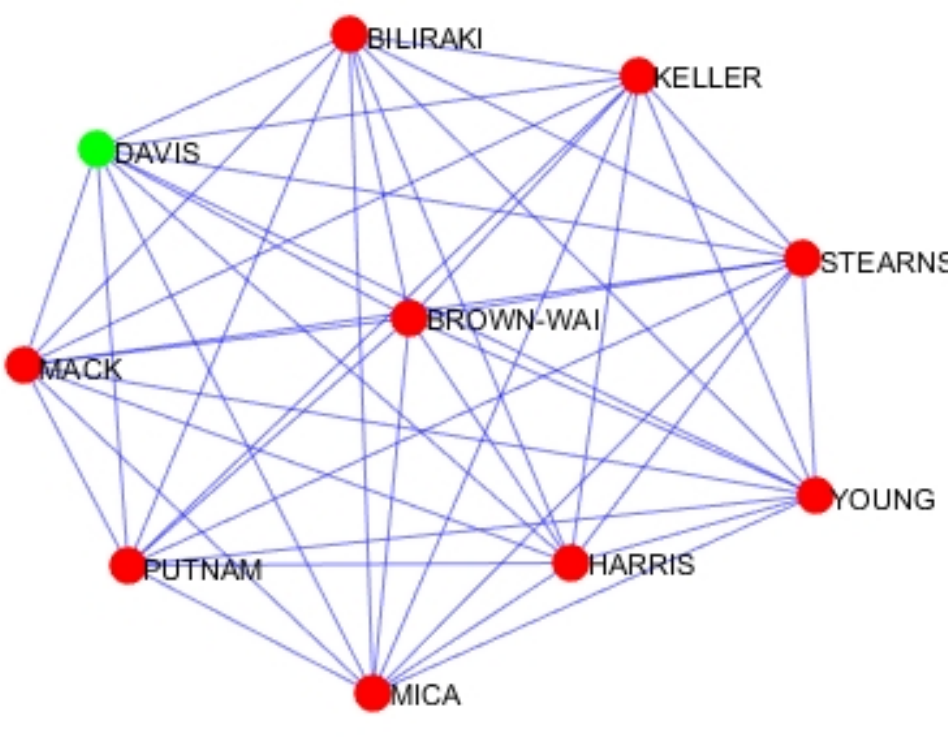}
    \caption{}
  \end{subfigure}
 \begin{subfigure} {0.3\textwidth}
   \includegraphics[width=\textwidth]{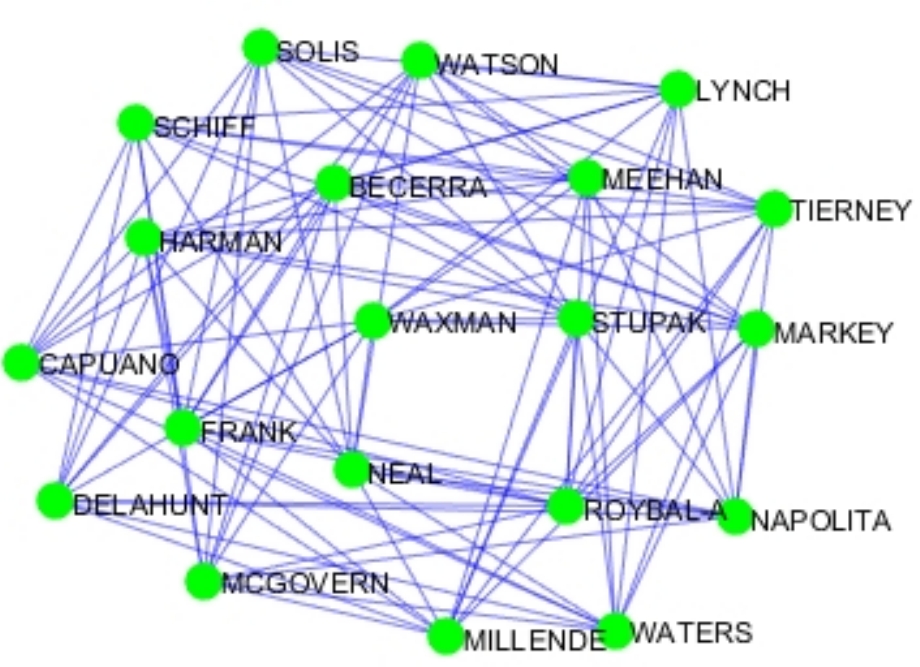}
    \caption{}
  \end{subfigure}
  \hfill
  \begin{subfigure} {0.3\textwidth}
    \includegraphics[width=\textwidth]{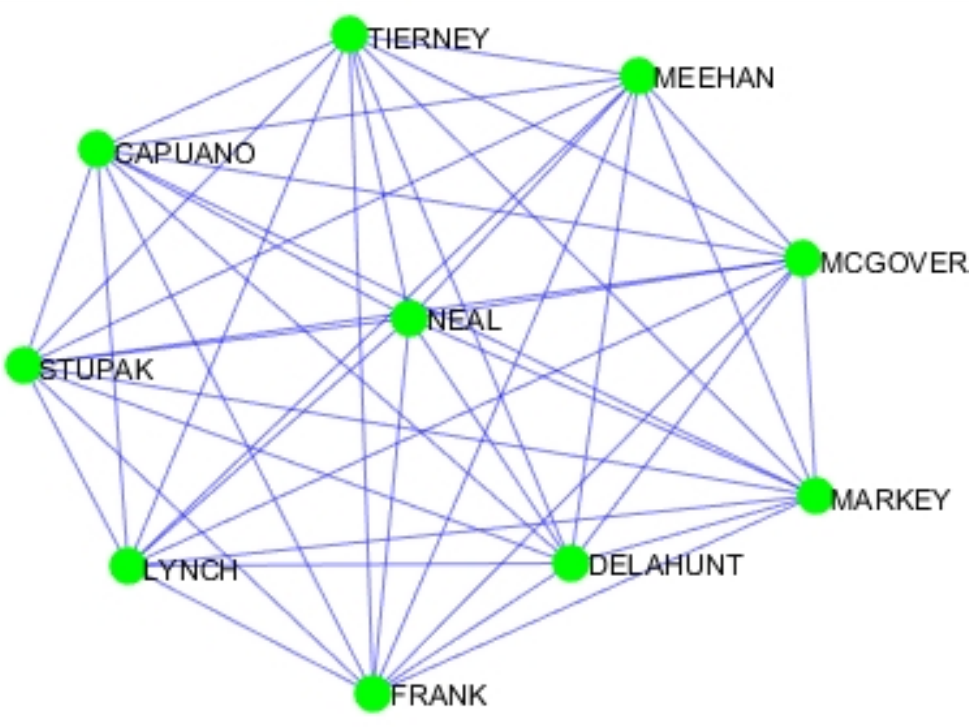}
    \caption{}
  \end{subfigure}
    \hfill
  \begin{subfigure} {0.3\textwidth}
    \includegraphics[width=\textwidth]{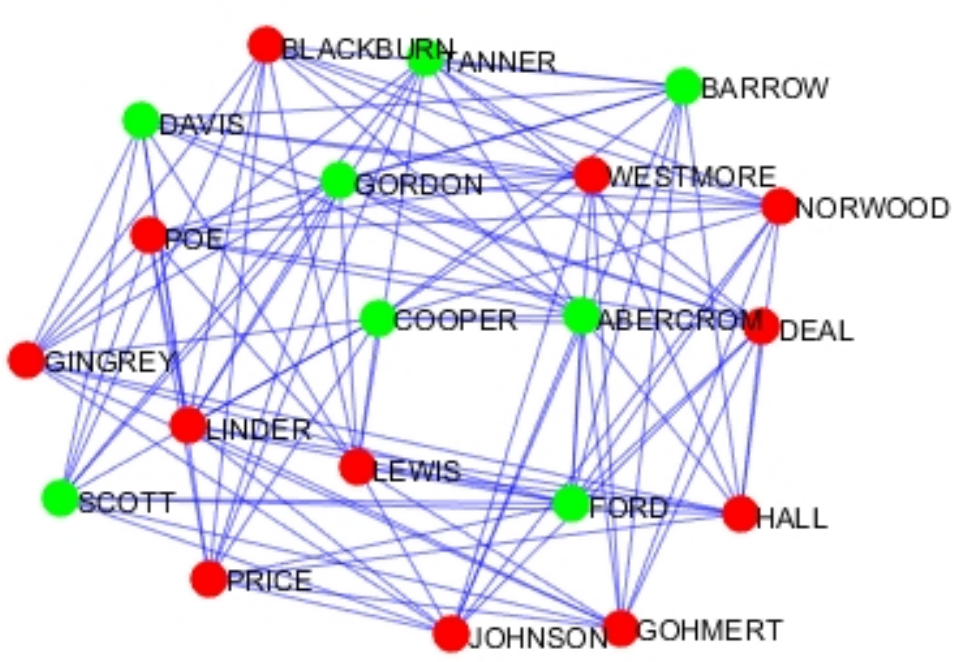}
    \caption{}
  \end{subfigure}
  \hfill
  \begin{subfigure} {0.3\textwidth}
    \includegraphics[width=\textwidth]{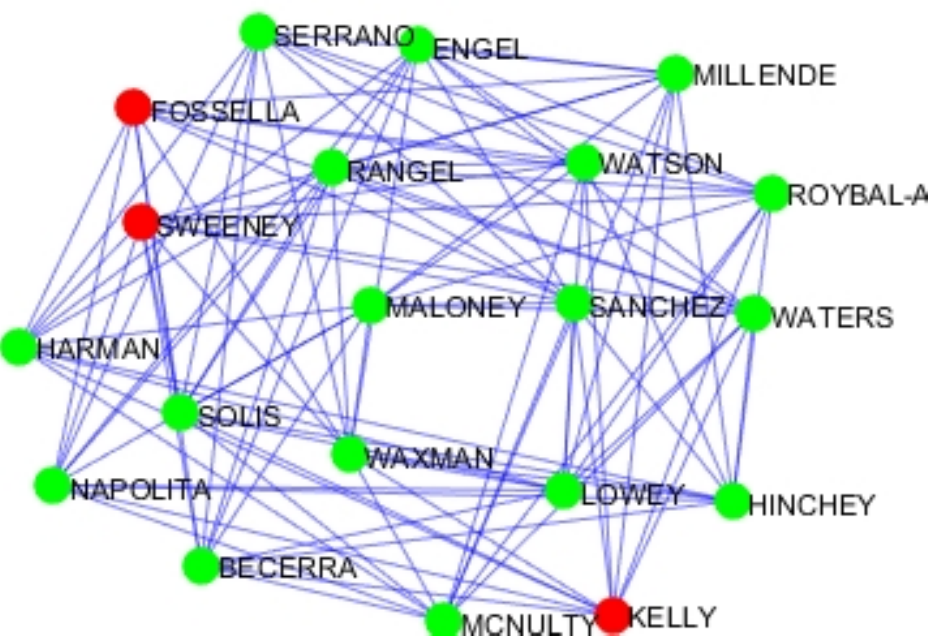}
    \caption{} 
  \end{subfigure}
      \hfill
  \begin{subfigure} {0.55\textwidth}
    \includegraphics[width=\textwidth]{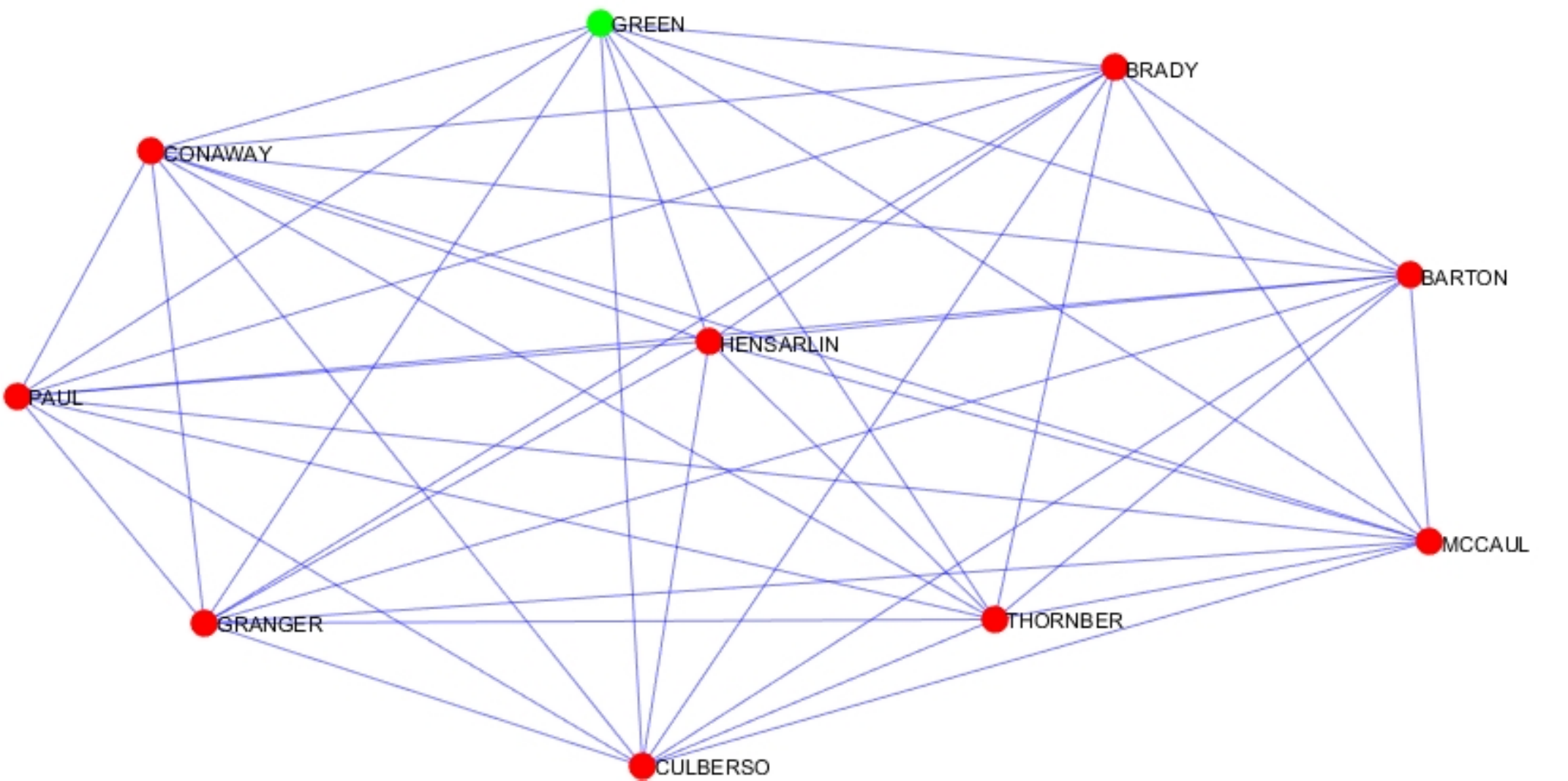}
    \caption{}
  \end{subfigure}
  \caption{Dense subgraphs identified by SSONA for the House voting data with an inclusion
cutoff value of 0.6. Subfigures correspond to a densely connected area in Figure~\ref{figvot1} for the symmetric structured matrix $Z_2+Z_2^\top$. The nodes represent House members, with red and green colored nodes corresponding to Republicans and  Democrats, respectively. A blue line corresponds to an edge between two nodes.}
\label{figvot12}
\end{figure}

\begin{figure}[!ht]
\captionsetup[subfigure]{labelformat=empty}
 \begin{subfigure} {0.45\textwidth}
     \includegraphics[width=\textwidth]{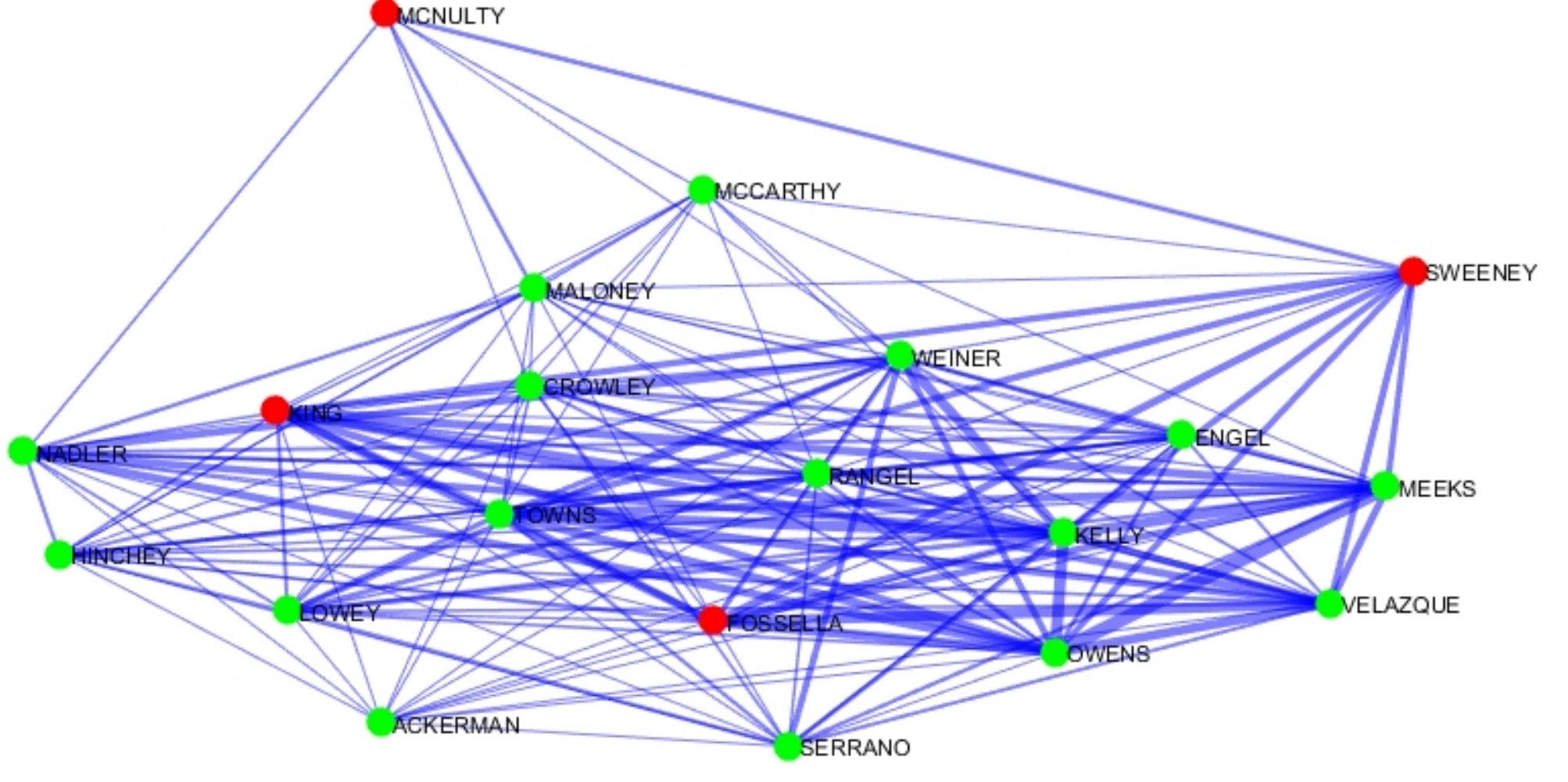}
    \caption{}
  \end{subfigure}
  \hfill
  \begin{subfigure} {0.45\textwidth}
       \includegraphics[width=\textwidth]{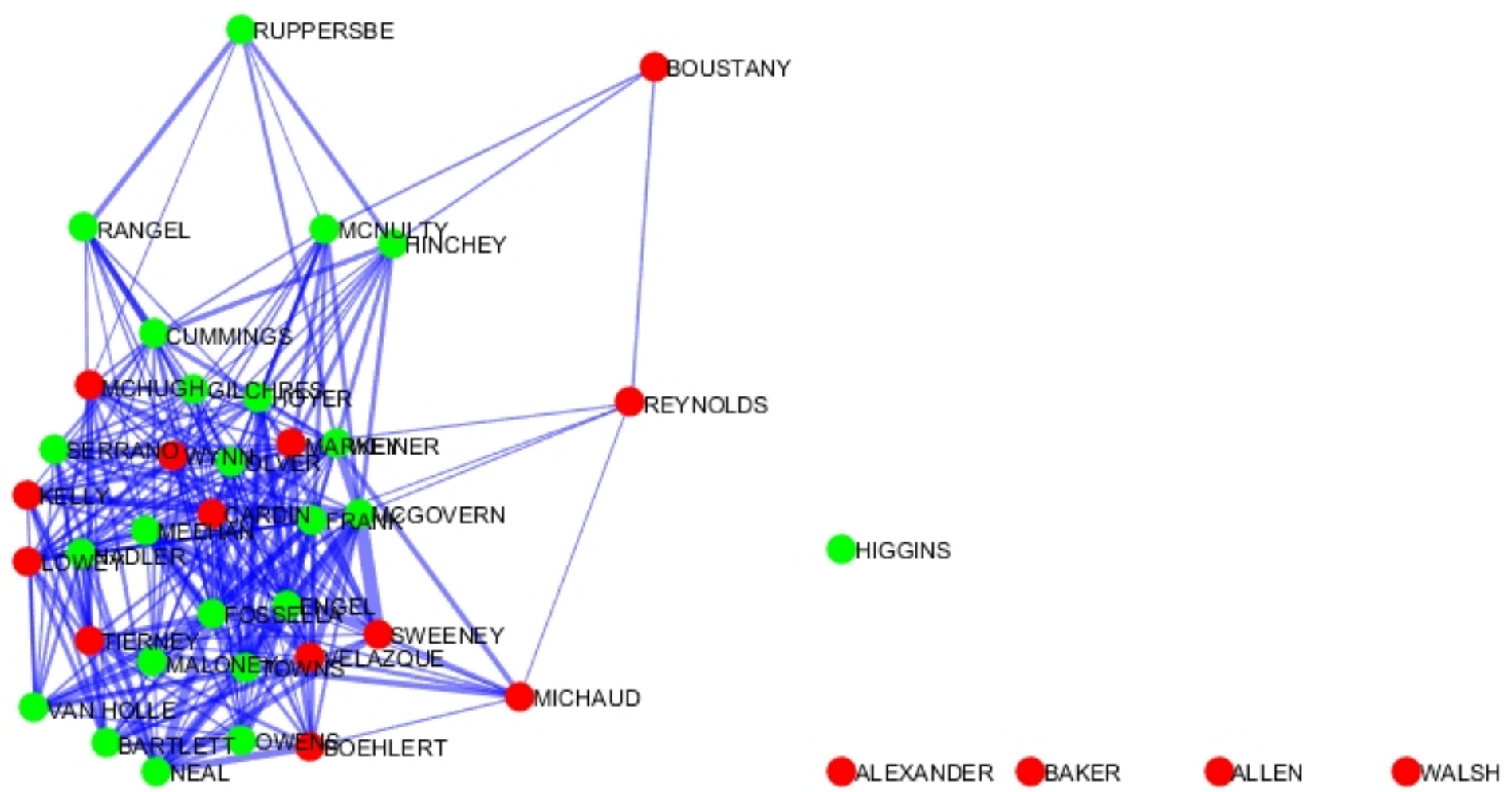}
    \caption{}
    \end{subfigure}
  \hfill
  \begin{subfigure} {0.45\textwidth}
  \includegraphics[width=\textwidth]{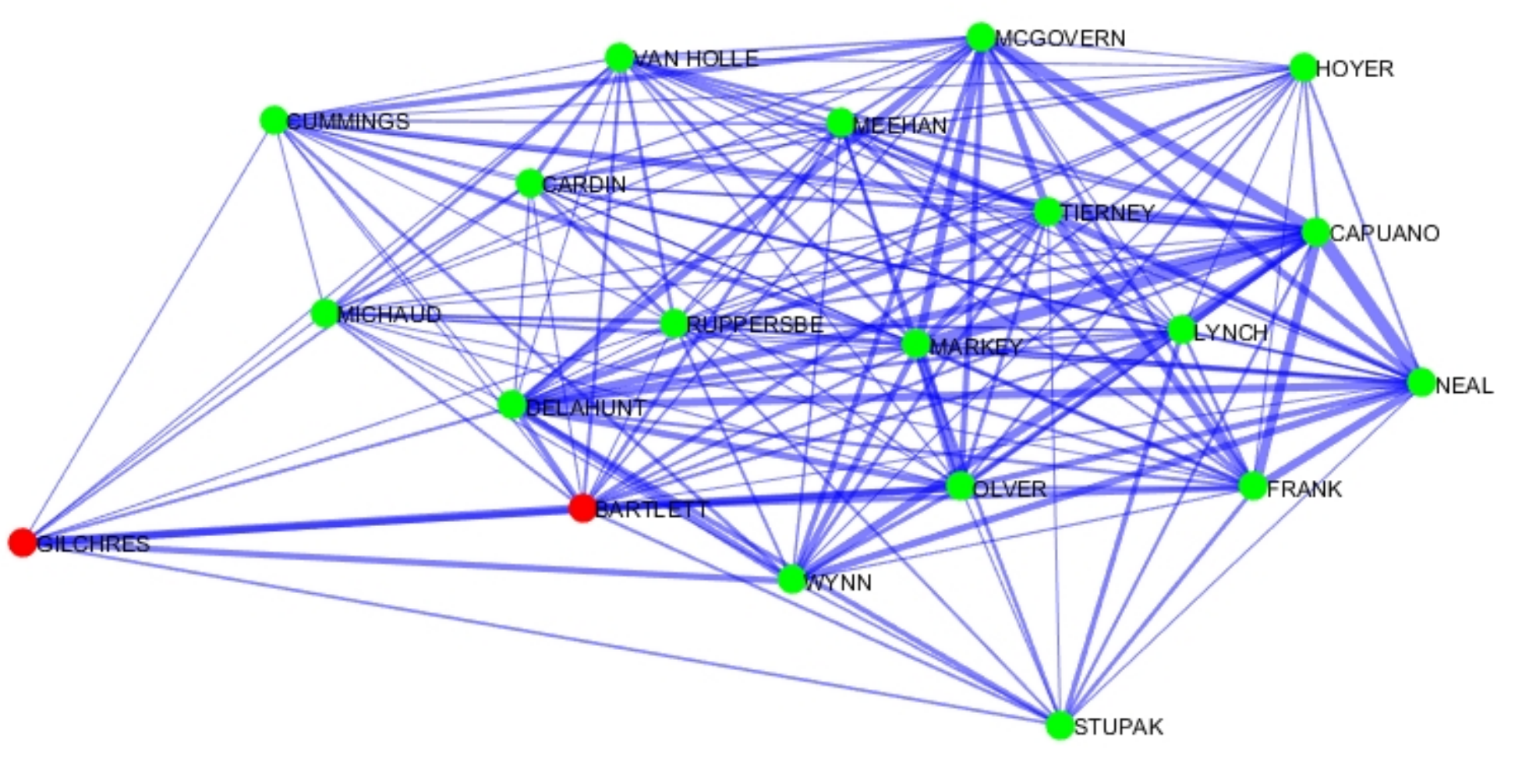}
    \caption{}
  \end{subfigure}
 \begin{subfigure} {0.45\textwidth}
      \includegraphics[width=\textwidth]{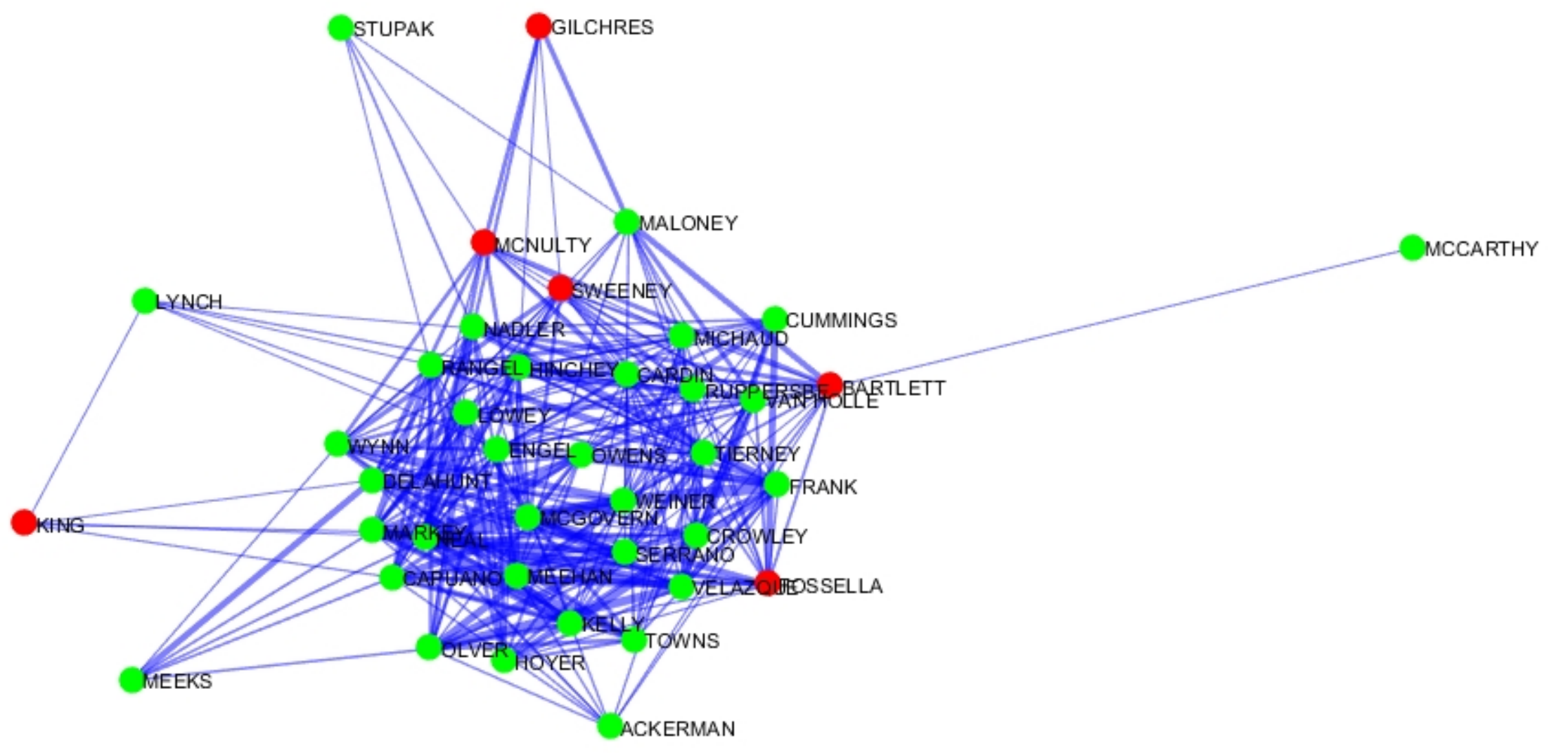}
    \caption{}
  \end{subfigure}
  \hfill
  \begin{subfigure} {0.45\textwidth}
       \includegraphics[width=\textwidth]{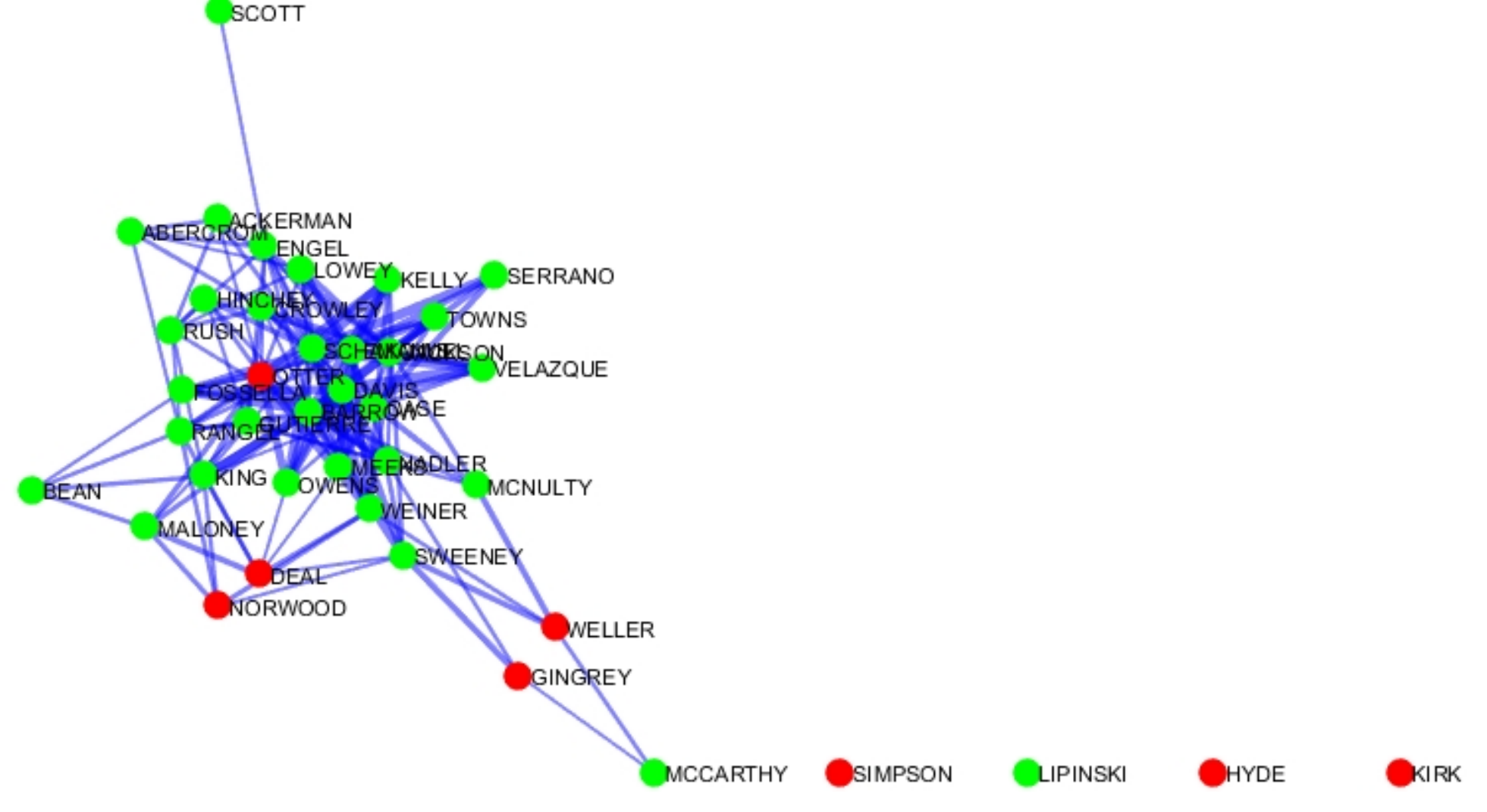}
    \caption{}
  \end{subfigure}
    \hfill
  \begin{subfigure} {0.45\textwidth}
       \includegraphics[width=\textwidth]{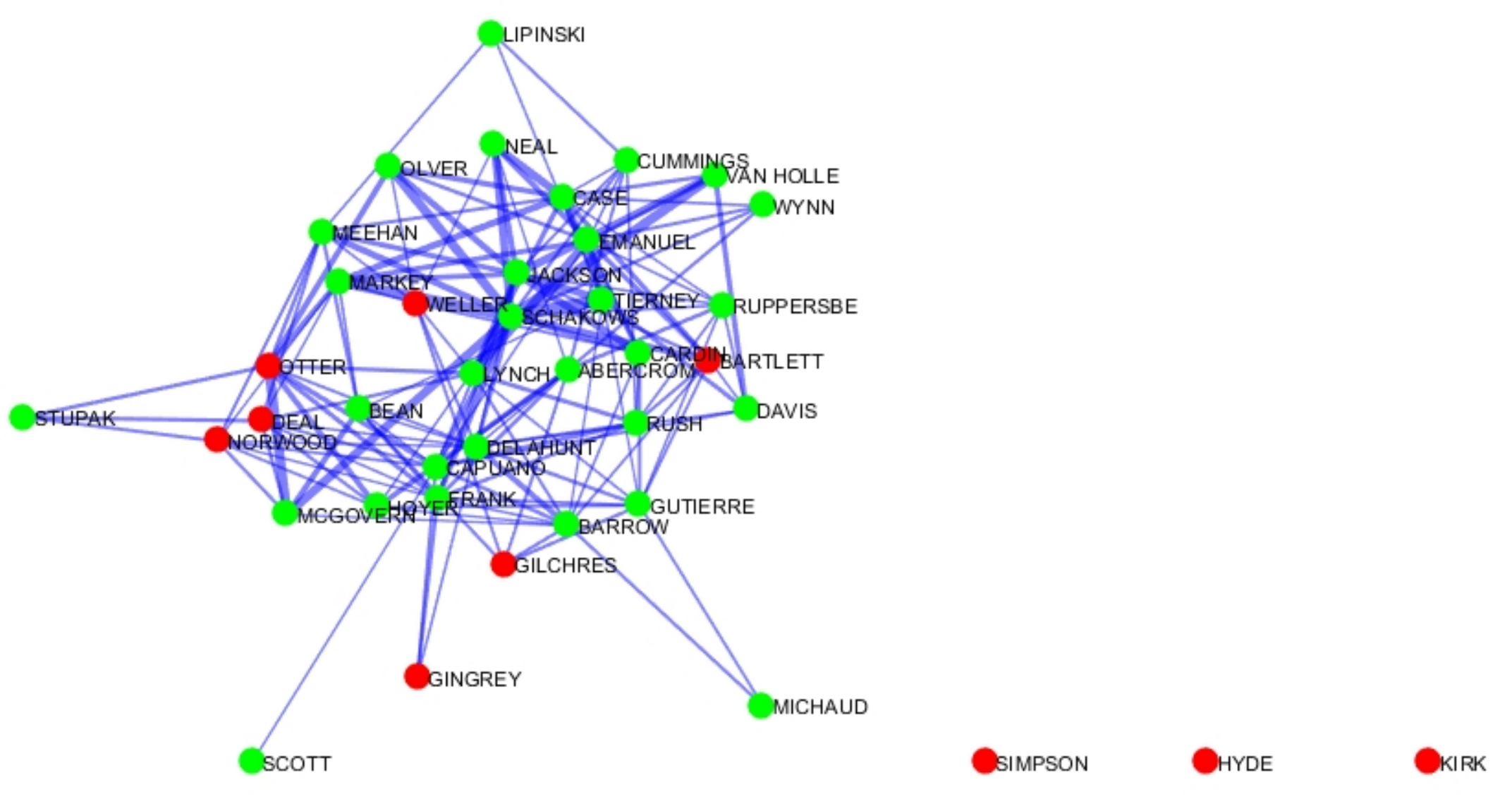}
    \caption{}
  \end{subfigure}
  \caption{Dense subgraphs identified by SSONA for the House voting data with an inclusion
cutoff value of 0.6. Subfigures correspond to a densely connected area in Figure~\ref{figvot1} for the symmetric structured matrix $Z_3+Z_3^\top$. The nodes represent House members, with red and green colored nodes corresponding to Republicans and  Democrats, respectively. A blue line corresponds to an edge between two nodes.}
\label{figvot13}
\end{figure}

\begin{figure}[!ht]
  \centering
      \includegraphics[width=.8\textwidth]{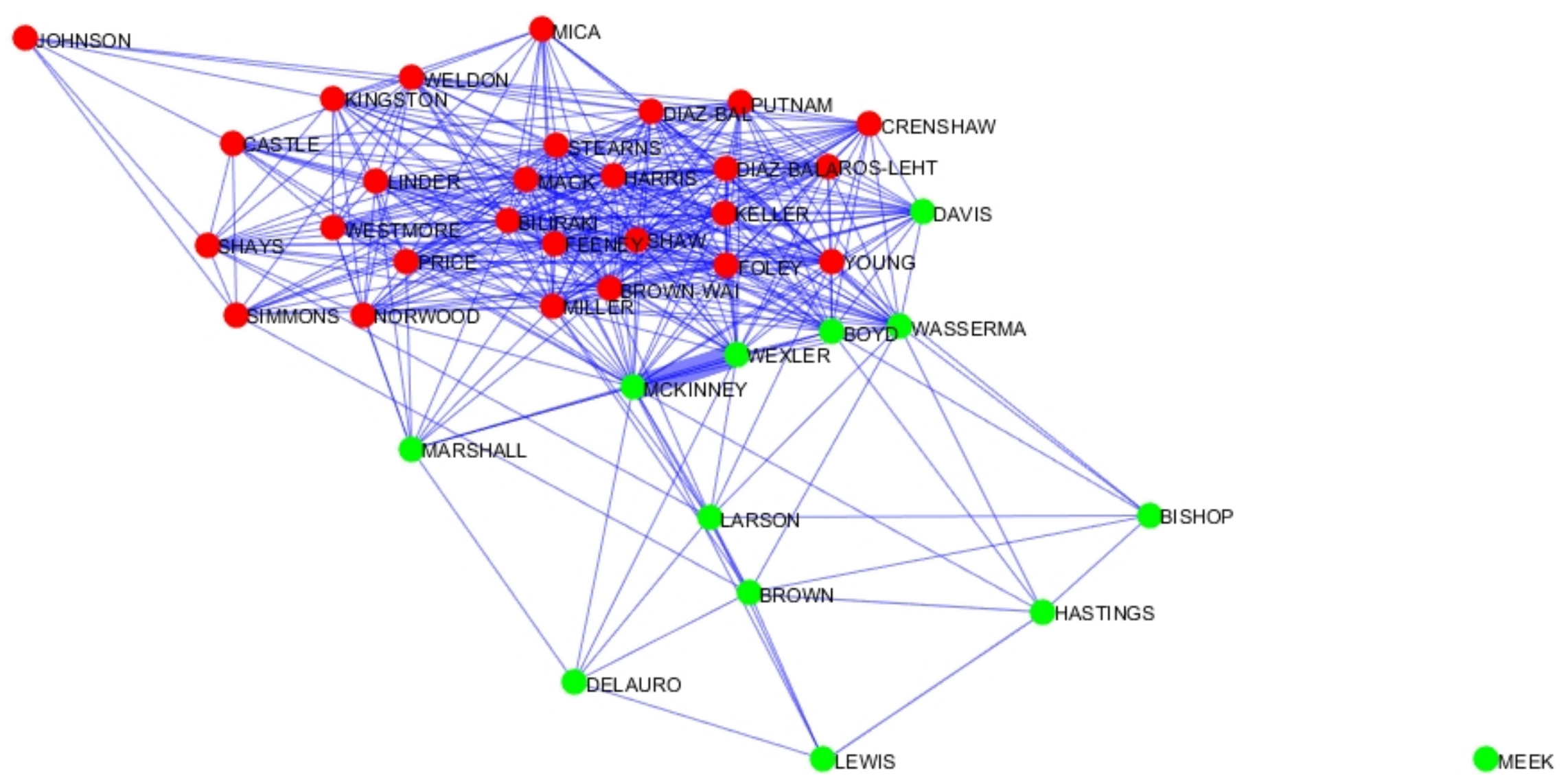}
      \caption{Dense subgraph identified by SSONA for the House voting data with an inclusion
cutoff value of 0.6. Subfigures corresponds to a densely connected area in Figure \ref{figvot1} for the symmetric structured matrix $Z_5+Z_5^\top$. The nodes represent House members, with red and blue node colors corresponding to Republicans and  Democrats, respectively. A blue line corresponds to an edge between two nodes.}\label{figvot14}
\end{figure}

The network representation of subgraphs, with a cut-off value of 0.6, is given in Figures~\ref{figvot12}, \ref{figvot13} and \ref{figvot14}. We only plot the edges associated with the subgraphs to enhance the visual reading of densely correlated areas. An interesting result of applying SSONA on this data set is the clear separation between members of the Democratic and Republican parties, as expected (see, Figures~\ref{figvot12}, \ref{figvot13} and \ref{figvot14}). Moreover, voting relationships within the two parties exhibit a clustering structure, which a closer inspection of the votes and subsequent analysis showed was mainly driven by the position of the House member on the ideological/political spectrum.

Other interesting patterns emerging from the analysis is that SSONA recovers members of opposite parties as a sparse component in each subgraph (see, Figures~\ref{figvot12}, \ref{figvot13} and \ref{figvot14}). For instance, Figure~\ref{figvot13} shows that Republican members such as Simpson, Kirk and Hyde are sparsely connected in a clustered group of Democratic members. This is possibly due to the overall centrist record of Kirk and alignment of Hyde and Simpson on selected issues. Similarly, Figure~\ref{figvot13} indicates that Democratic members Bishop, Hastings and Meek are approximately sparsely connected to a subgraph of Republican members. Bishop from Georgia has compiled a fairly conservative voting record. The same conclusion can be derived from  Figure~\ref{figvot13}. Indeed, Figures~\ref{figvot12}, \ref{figvot13} and \ref{figvot14} reveals that there are strong positive associations between members of the same party and negative associations between members of opposite parties. Obviously, at the higher cutoff value the dependence structure between members of opposite parties becomes sparser.

Other patterns of interest include a strong dependence between members of two opposite parties in selected subgraphs when the members come from the same state, as is the case for New York state members Jerrold Nadler (D), Anthony D. Weiner (D), Ed Towns (D),  Major Owens (D), Nydia Velázquez (D), Vito Fossella (R), Carolyn B. Maloney (D), Charles B. Rangel (D), José Serrano (D), Eliot L. Engel (D), Nita Lowey (D), Sue W. Kelly (R), John E. Sweeney (R), Michael R. McNulty (D), Maurice Hinchey (D), John M. McHugh (R), Sherwood Boehlert (R), Jim Walsh (R), Tom Reynolds (R), Brian Higgins (D) -see Figure~\ref{figvot13}.
However, in this instance, there is also a cluster of positive associations between Democrats.

In summary, SSONA provides deeper insights into relationships between House members, going beyond the obvious separation into two parties,
according to their voting record.

\paragraph{Analysis of a breast cancer data set.}
We applied SSONA to a data set containing 800~gene expression measurements from large epithelial cells obtained from 255~patients with breast cancer.
\begin{figure}[!ht]
\captionsetup[subfigure]{labelformat=empty}
 \begin{subfigure} {0.3\textwidth}
 \includegraphics[width=\textwidth]{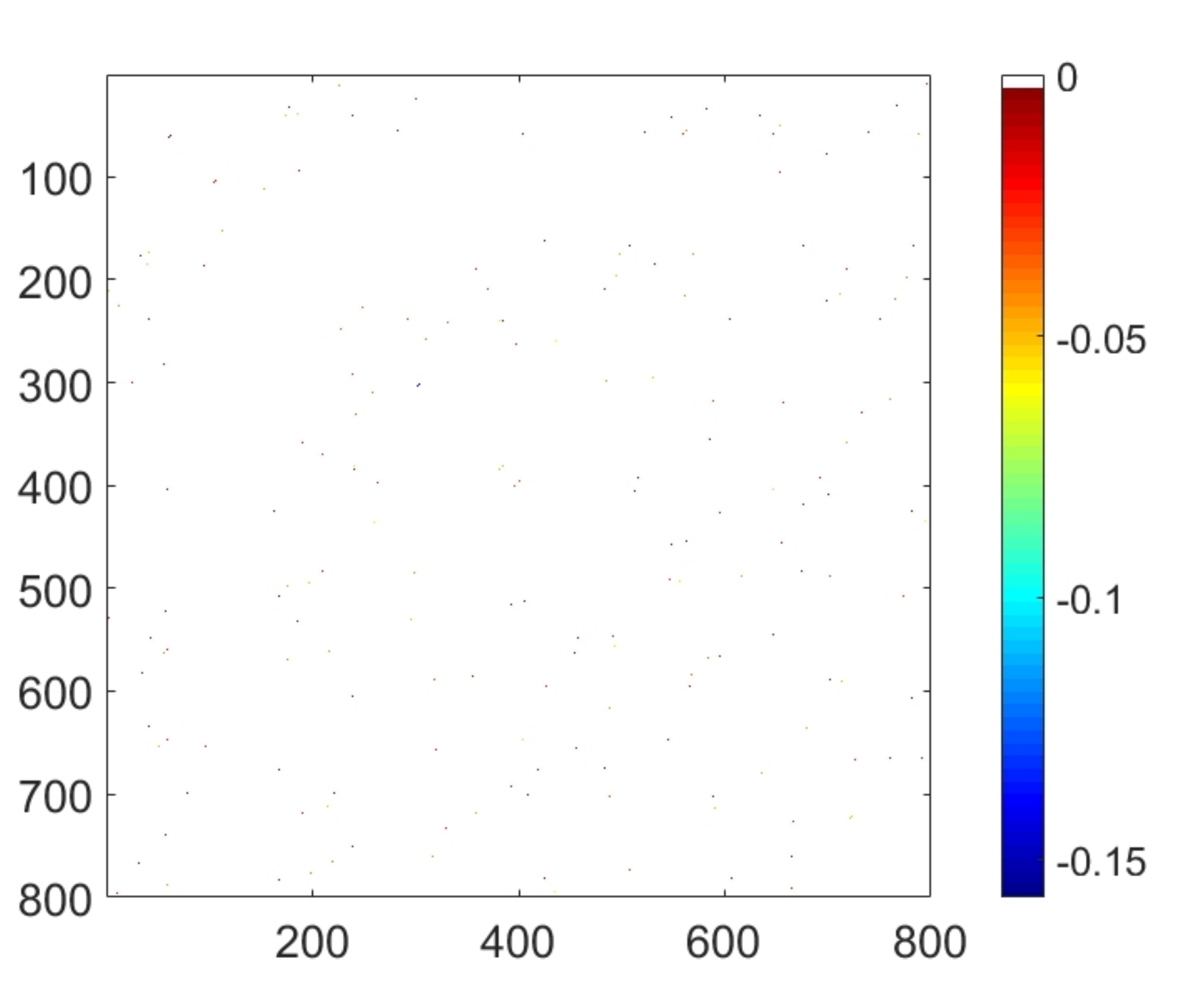}
    \caption{$Z_1+Z_1^\top$}
  \end{subfigure}
  \hfill
  \begin{subfigure} {0.3\textwidth}
\includegraphics[width=\textwidth]{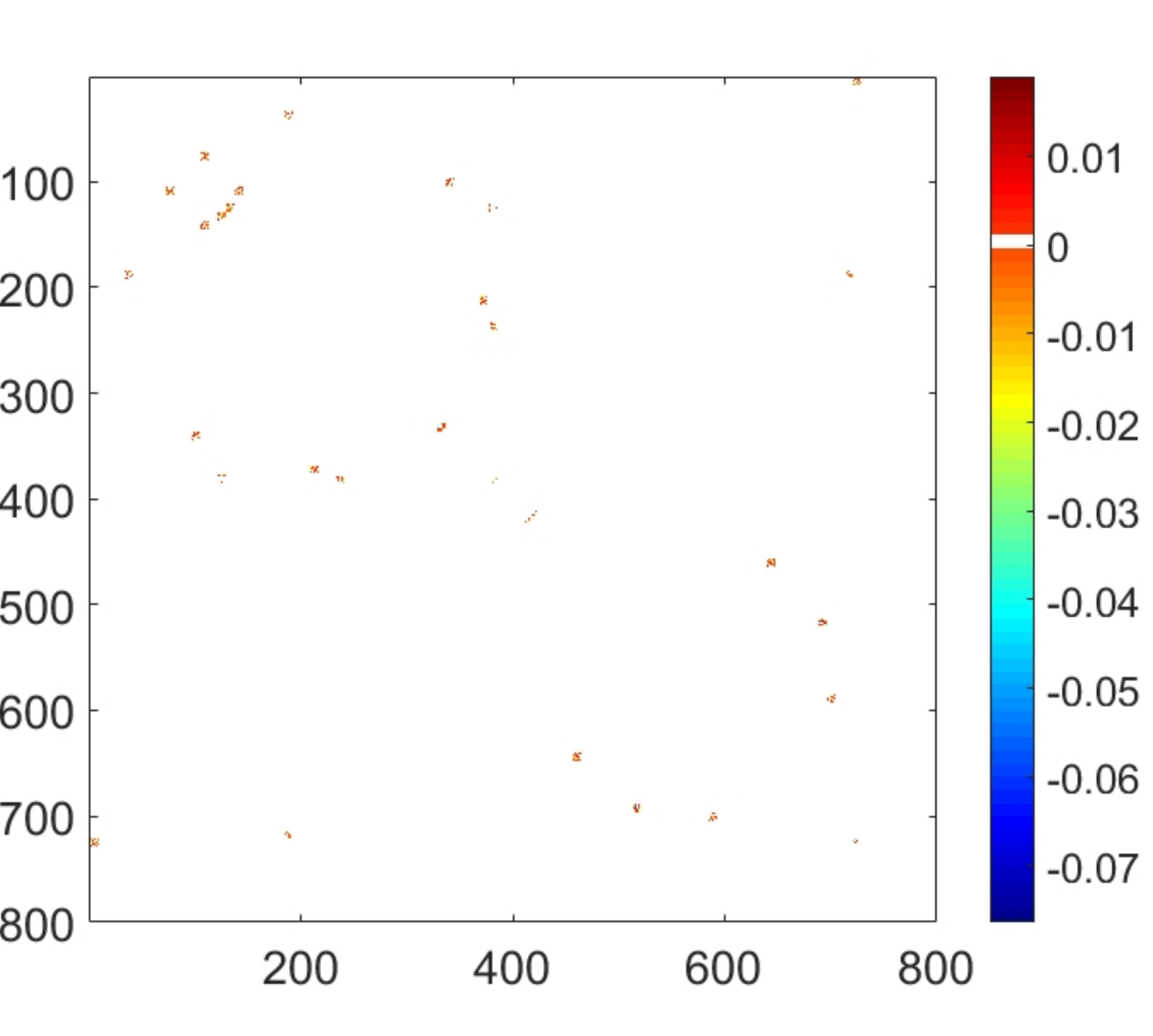}
    \caption{$Z_2+Z_2^\top$}
    \end{subfigure}
  \hfill
  \begin{subfigure} {0.3\textwidth}
 \includegraphics[width=\textwidth]{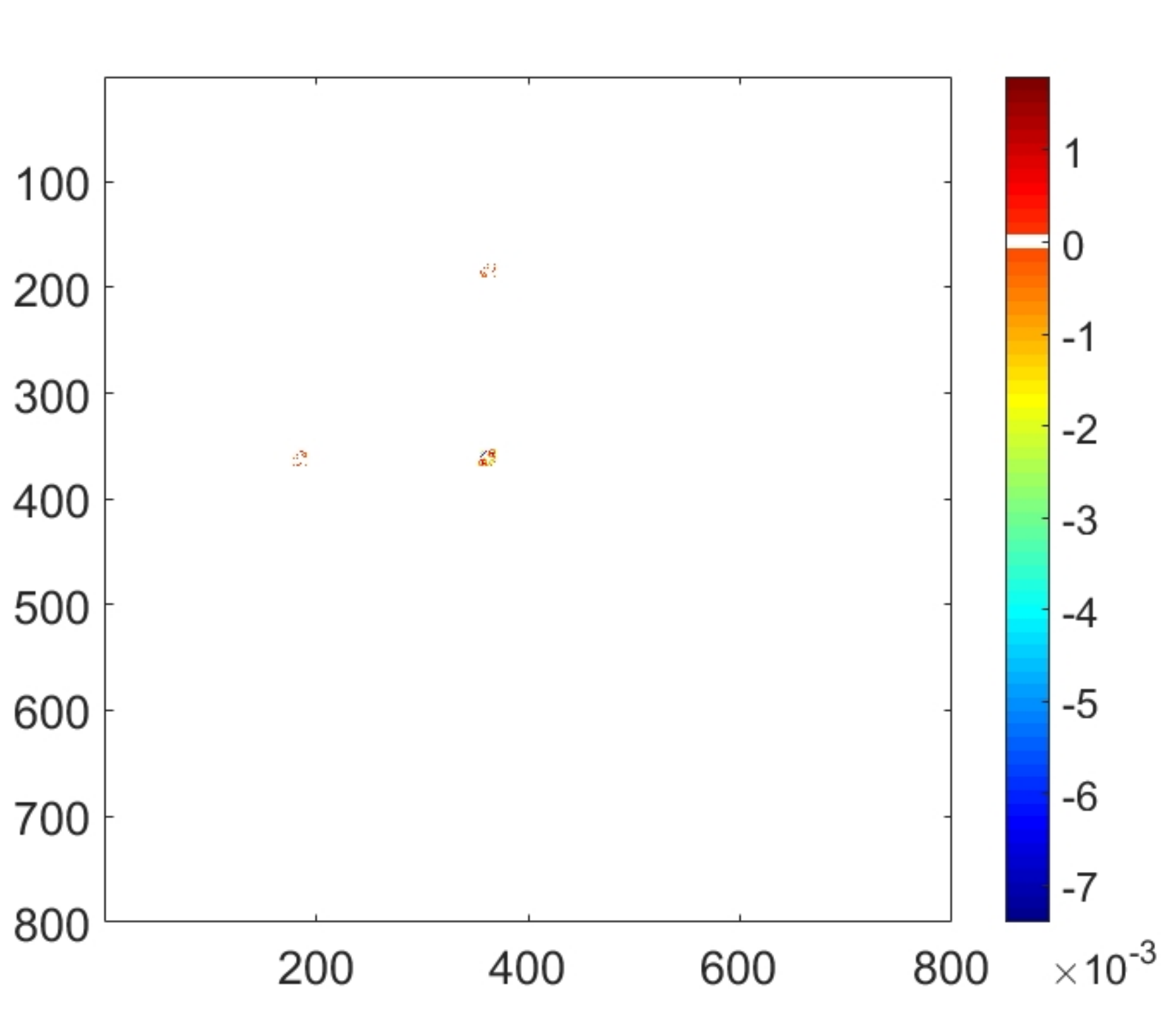}
    \caption{$Z_3+Z_3^\top$}
  \end{subfigure}
 \begin{subfigure} {0.3\textwidth}
 \includegraphics[width=\textwidth]{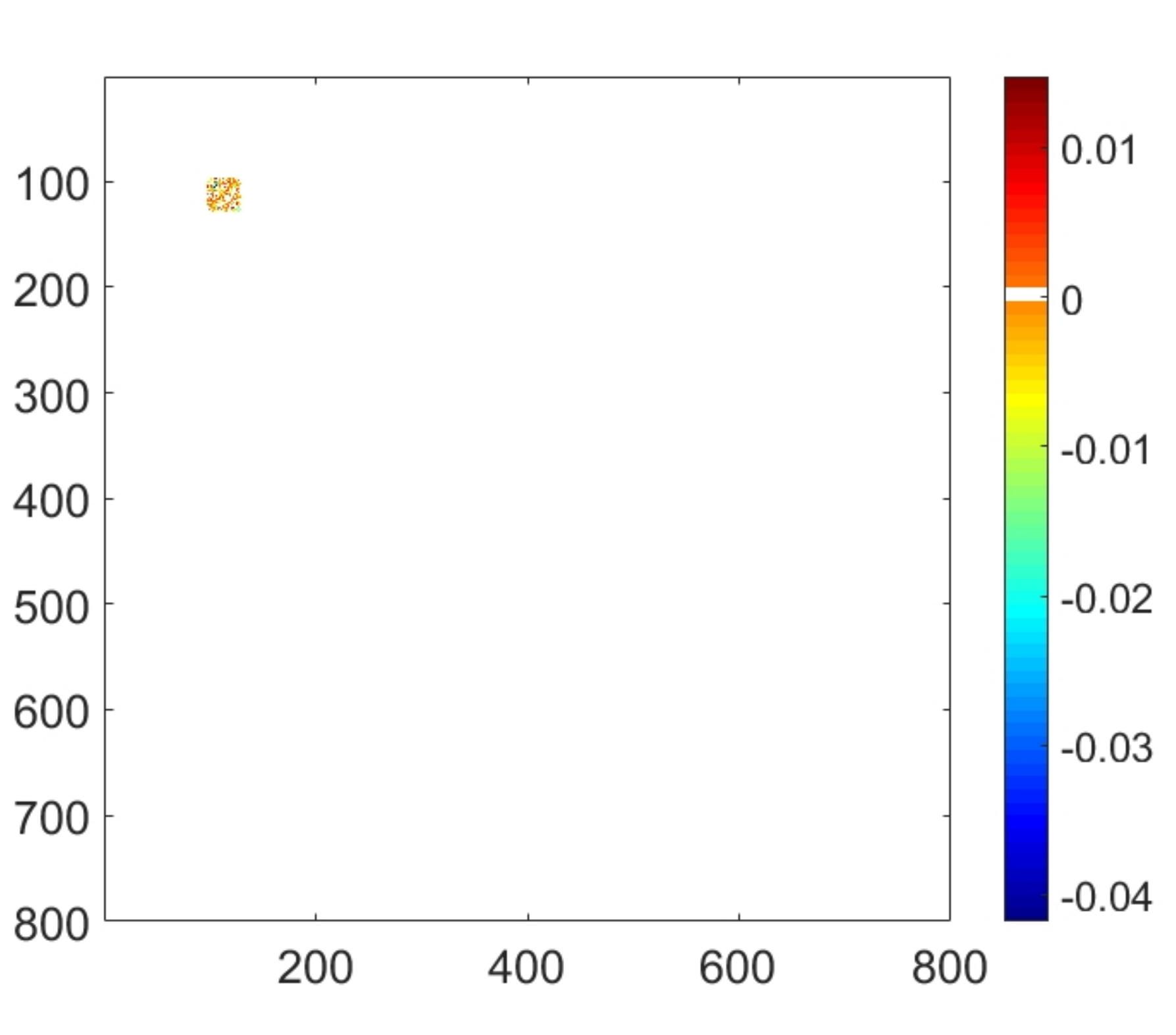}
    \caption{$Z_4+Z_4^\top$}
  \end{subfigure}
  \hfill
  \begin{subfigure} {0.3\textwidth}
  \includegraphics[width=\textwidth]{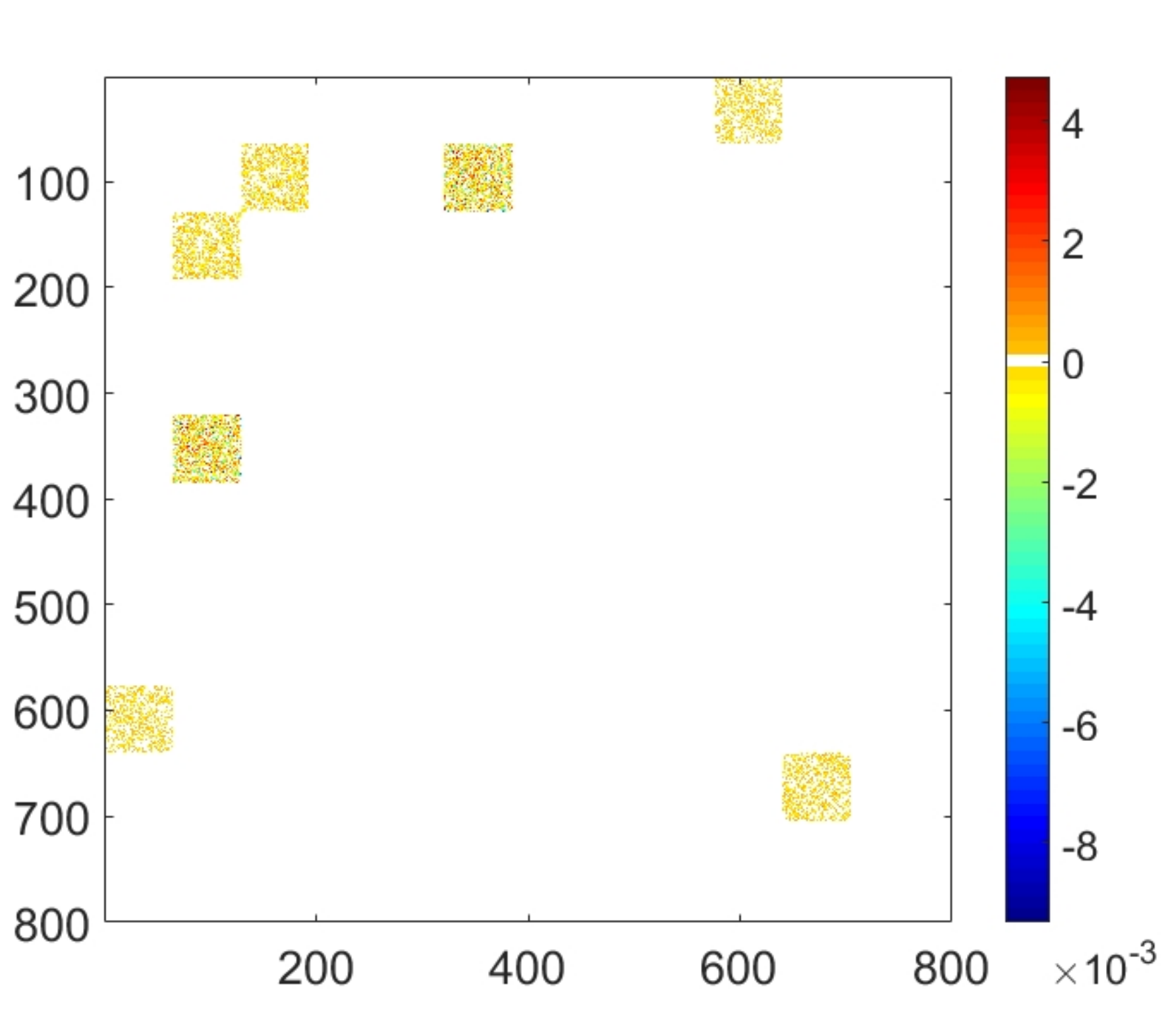}
    \caption{$Z_5+Z_5^\top$}
  \end{subfigure}
  \hfill
  \begin{subfigure} {0.3\textwidth}
  \includegraphics[width=\textwidth]{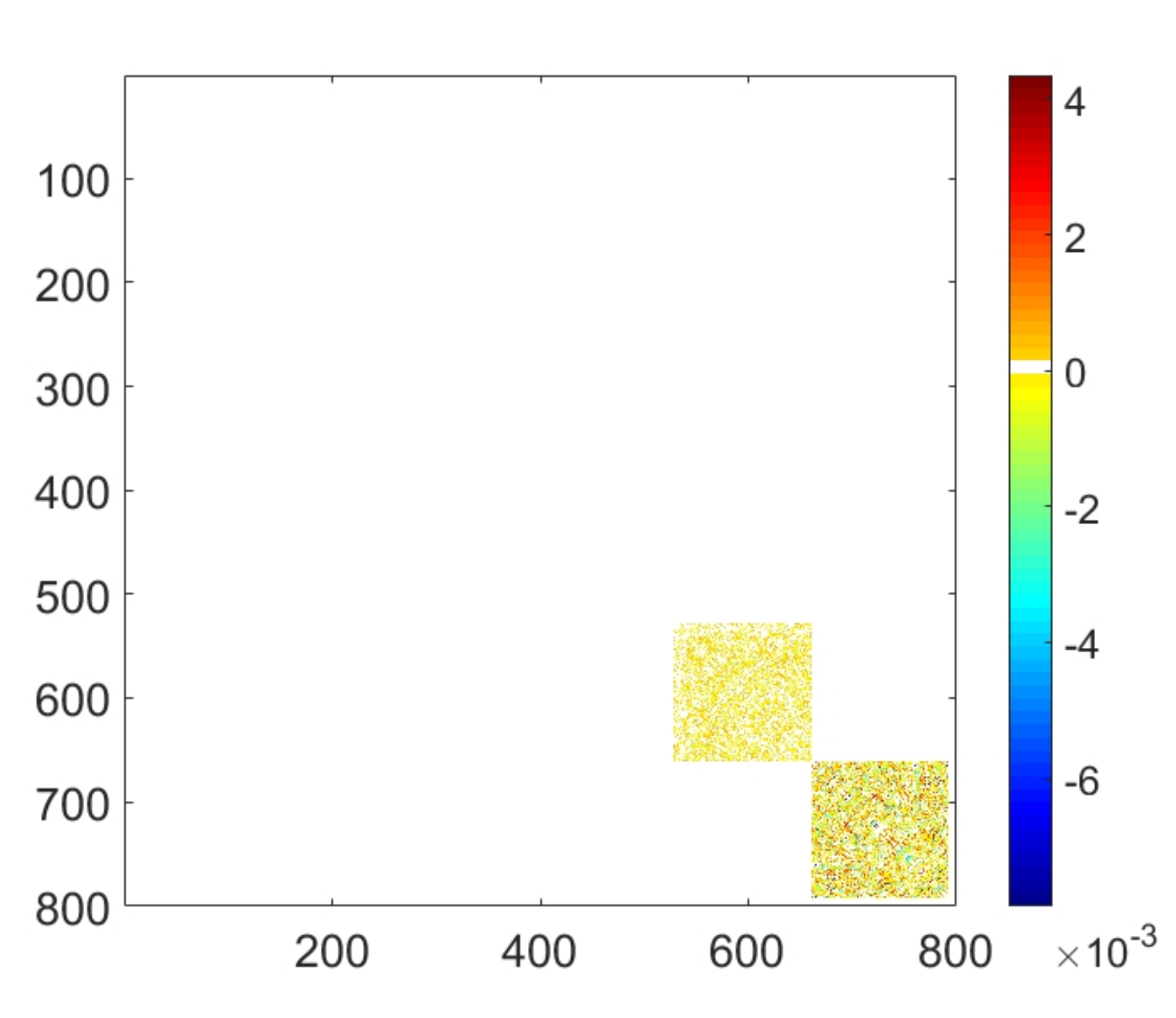}
    \caption{$Z_6+Z_6^\top$} 
  \end{subfigure}
\caption{Heat map of the structured precision matrix $\Theta$ decomposed into $Z_1+Z_1^\top+\dots +Z_6+Z_6^\top$ in the breast cancer data set, estimated by SSONA.} \label{figbreast}
\end{figure}
The goal is to capture regulatory interactions amongst the genes, as well as to identify genes that tend to have interactions with other genes in a group and hence act as master regulators, thus providing insights into the molecular circuitry of the disease. Figure \ref{figbreast} depicts the heat map of the estimated adjacency matrix for the breast cancer data set. As it is clear in Figure \ref{figbreast}, $Z_2+Z_2^\top,\dots,Z_5+Z_5^\top$ and $Z_6+Z_6^\top$ show that selected genes are densely connected, which is not the case when employing the the graphical lasso algorithm (see, Figure \ref{breaspa}). Therefore, SSONA can provide an intuitive explanation of the relationships among the genes in the breast cancer data set (see, Figure \ref{impgen} and \ref{areas} for two examples). These genes connectivity in the tumor samples may indicate a relationship that is common to an important subset of cancers. Many other genes belong to this network, each indicating a potentially interesting interaction in cancer biology. We omit the full list of densely connected genes in our estimated network and provide a complete list in the on-line supplementary materials available in the first author's homepage.
\begin{figure}
\centering
  \includegraphics[width=.5\textwidth]{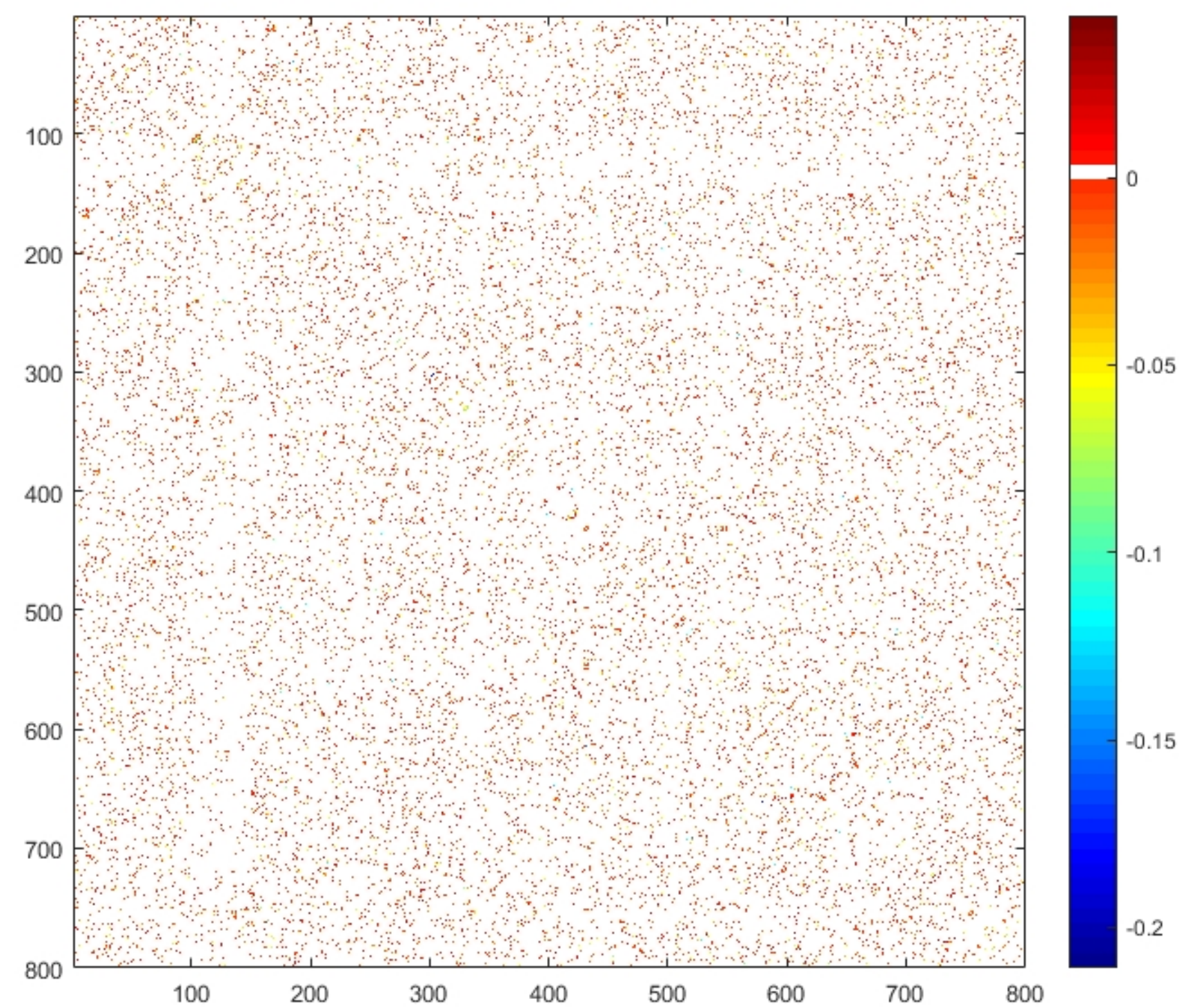}
 \caption{Heatmap of the inverse covariance matrix in the breast cancer data sets, estimated from graphical lasso \citep{Friedman07}.}  \label{breaspa}
\end{figure}

\begin{figure}[!ht]
  \centering
  \includegraphics[width=.75\textwidth]{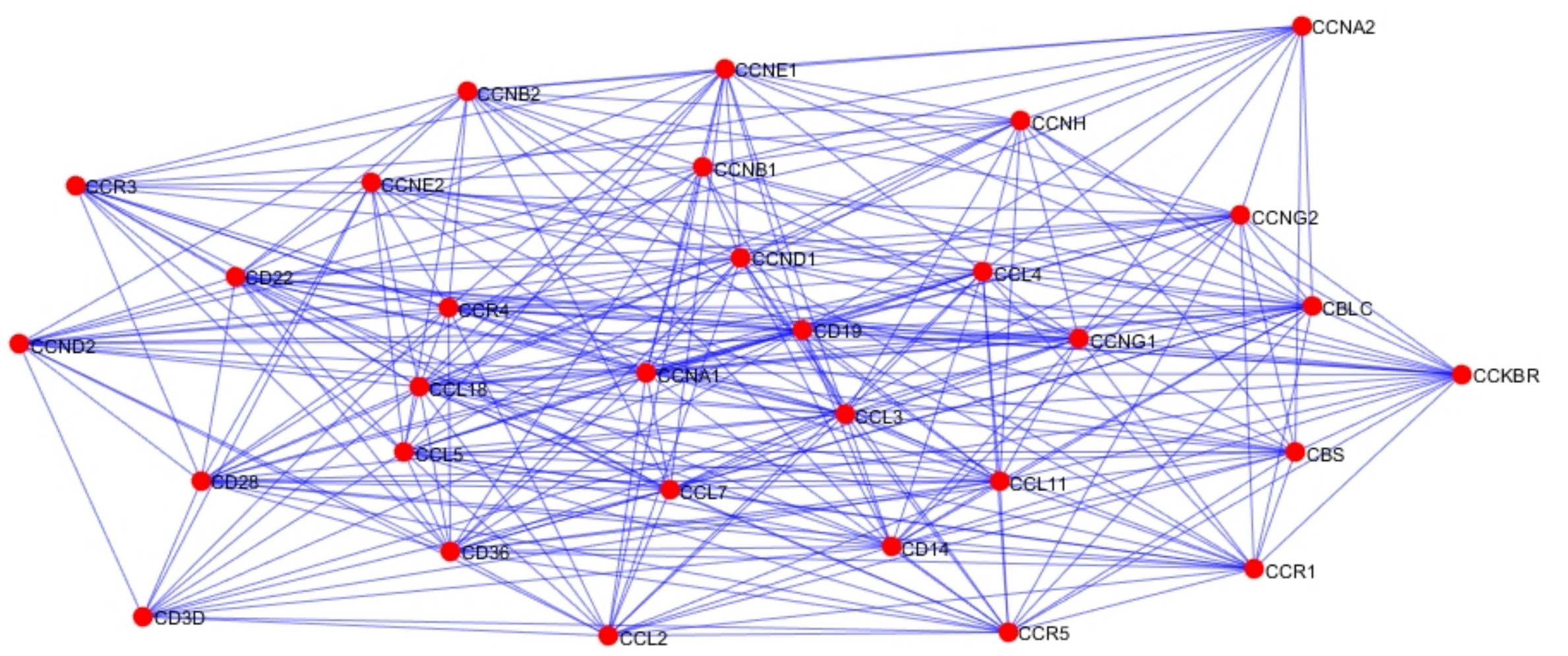}
      \caption{Network layout of grouped genes identified by SSONA for the breast cancer data set. Subfigure corresponds to a densely connected component in Figure~\ref{figbreast} for the structured matrix $Z_4+Z_4^\top$.}
      \label{impgen}
\end{figure}
\begin{figure}[!ht]
  \centering
  \includegraphics[width=.75\textwidth]{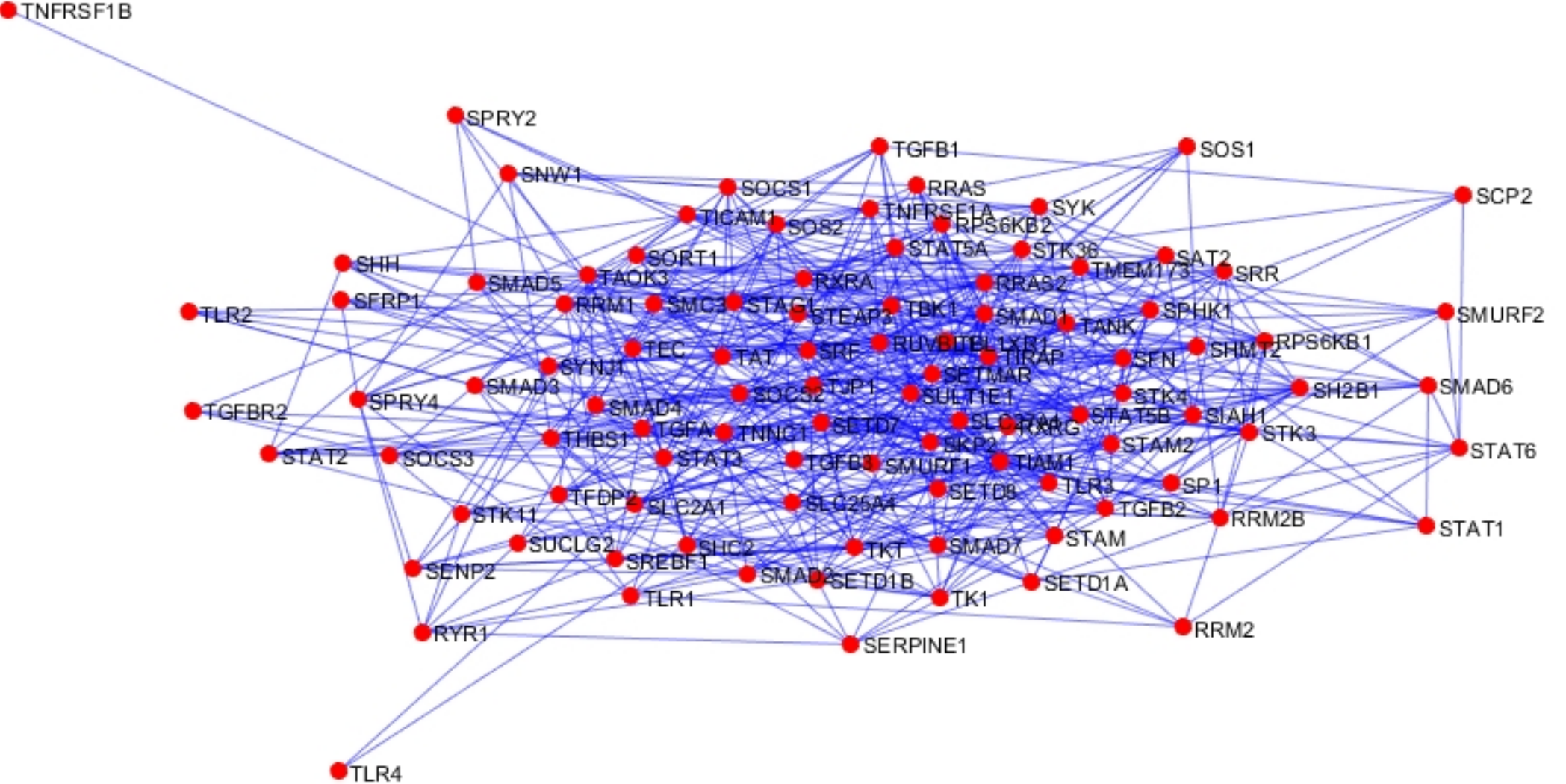}
      \caption{Network layout of grouped genes identified by SSONA for the breast cancer data set. Subfigure corresponds to a densely connected component in Figure~\ref{figbreast} for the structured matrix $Z_6+Z_6^\top$.}
      \label{areas}
\end{figure}

\section{Conclusion}
In this paper, a new structured norm minimization method for solving multi-structure graphical model selection problems is proposed. Using the proposed SSON, we can efficiently and accurately recover the underlying network structure. Our method utilizes a class of sparse structured norms in order to achieve higher order accuracy in approximating the decomposition of the parameter matrix in Markov Random Field and Gaussian Covariance Graph models.
We also provide a brief discussion of its application to regression and classification problems. Further, we introduce a linearized
multi-block ADMM algorithm to solve the resulting optimization problem. The global convergence of the algorithm is established without any upper bound on the penalty parameter. We applied the proposed methodology to a number of real and synthetic data sets that establish its overall usefulness and superior performance to competing methods in the literature.

\acks{The authors would like to thank the Editor and three anonymous referees for many constructive comments and suggestions that
improved significantly the structure and readability of the paper.
This work was supported in part by NSF grants DMS-1545277, DMS-1632730, NIH grant
1R01-GM1140201A1 and by the UF Informatics Institute.
}
\vskip 0.2in
\appendix

\section{Update for $\Theta$}
In each iteration of Algorithm~\ref{alg:1} the update for $\Theta$  depends on the form of the loss function $g(\Theta)$. We consider the following cases to update $\Theta$:
\begin{enumerate}
  \item
The update for $\Theta_1$ in Algorithm~1 (step~2(a)) can be obtained by minimizing
$$
\trace (\hat{\Sigma}\Theta_1)-\log \det \Theta_1 + \dfrac{\gamma}{2} \| \Theta_1 - (\sum_{i=1}^{n}Z_i^{k}+{Z_i^{k}}^\top+ E^k+ \frac{1}{\gamma}\Lambda^k )\|^2_{F},
$$
with respect to $\Theta_1$ (note that the constraint $\Theta_1 \in \Sc$ in  \eqref{eq:5} is treated as an implicit constraint, due to the domain of definition of the $\log \det$ function).  This can be shown to have the solution
$$ \Theta_1 = \frac{1}{2}U\Big(D+\sqrt{D^2+\frac{4}{\gamma}I}\Big) U^T,$$
where $UDU^T$ stands for the eigen-decomposition of $\sum_{i=1}^{n}Z_i^{k}+{Z_i^{k}}^\top + E^k+\frac{1}{\gamma}\Lambda^k-\frac{1}{\gamma}\hat{\Sigma}$.

\item Update for $\Theta_2$ in Step~2(a) of Algorithm~\ref{alg:1} leads to the following optimization problem
\begin{eqnarray} \label{eq:13-1}
\nonumber
\minimize_{\Theta_3 \in \Sc}~\Phi(\Theta_2)&=& \sum_{j=1}^{p}\sum_{j'=1}^{p} \theta_{jj'} (X^T X)_{jj'}-\sum_{i=1}^{m}\sum_{j=1}^{p} \log \Big(1 + \exp [ \theta_{jj}+ \sum_{j'\neq j} \theta_{jj'} x_{ij'}]\Big)\\
 &+& \dfrac{\gamma}{2} \| \Theta_2-(\sum_{i=1}^{n}Z_i^{k}+{Z_i^{k}}^\top + E^k+\frac{1}{\gamma}\Lambda^k)\|^2_{F}.
\end{eqnarray}
We use a novel non-monotone version of the Barzilai-Borwein method \citep{Barz88,Raydan97,Flet05,Ataee14} to solve \eqref{eq:13-1}. The
details are given in Algorithm \ref{alg:2}.
\begin{algorithm}[H]
\caption{Non-monotone Barzilai Borwein Method for solving \eqref{eq:13-1}}\label{alg:2}
\begin{algorithmic}
    \Initialize{The parameters:}
    \begin{enumerate}[(a)]
    \item
        $\Theta^0=I $, $\Theta^{1}= 2 \Theta^0$,  $\alpha^1=1$ and $t^0 = 10$.
    \item
         A positive sequence $\{\eta^t\}$ satisfying $\sum_{k=1}^{\infty} \eta^t = \eta <\infty$.
    \item
        Constants $\sigma > 0$, $\epsilon >0$, and $\nu \in (0,1)$.
    \end{enumerate}
    \Iterate{Until the stopping criterion $\dfrac{\|\Theta^t - \Theta^{t- 1}\|^2_{F}}{\| \Theta^{t-1} \|^2_{F}} \leq \epsilon$ is met:}
    \begin{enumerate}[]
     \item[~1.] $\mathbb{G}^{t} = -\alpha^{t}\nabla \Phi(\Theta^t)$.
     \item[~2.] Set $\rho = 1$.
     \item[~3.] \textbf{If}$~~t > t^0$, \textbf{then}
     \item \quad \textbf{While}~~$\|\Phi(\Theta^t+ \rho^t \mathbb{G}^{t})\|_{F} \leq \Phi(\Theta^t) + \eta^{t}- \sigma \rho^2 {\alpha^t}^2 \|\mathbb{G}^t\|^2_{F}$,~\textbf{do}\\
         \qquad Set $\rho= \nu \rho;$    \\
         \qquad \textbf{EndWhile}\\
     \textbf{EndIf}
     \item [~4.] Define $\rho^{t}=\rho$ and $\Theta^{t+1}=\Theta^t+\rho^t \mathbb{G}^{t}$.
     \item [~5.] Define $\alpha^{t+1}=\dfrac{\trace \Big({(\Theta^t-\Theta^{t+1})}^T {(\Theta^t-\Theta^{t+1})}\Big)}{\trace \Big({(\nabla \Phi(\Theta^t)-\nabla \Phi(\Theta^{t+1}))}^T (\Theta^t-\Theta^{t+1})\Big)}$
     \end{enumerate}
  \end{algorithmic}
\end{algorithm}
\item To update $\Theta_3$ in step~2(a), using \eqref{eq:7}, we have that
\begin{eqnarray*}
  \minimize_{\Theta_3}~\dfrac{1}{2}\|\Theta_3 - \hat{\Sigma}\|^2_{F} &+& \dfrac{\gamma}{2} \| \Theta_3-\Big(\sum_{i=1}^{n}Z_i^{k}+{Z_i^{k}}^\top+ E^k+ \frac{1}{\gamma}\Lambda^k \Big)\|^2_{F}\\&=&
\Big( \frac{1}{1+\gamma}(\hat{\Sigma}+ \gamma (\sum_{i=1}^{n}Z_i^{k}+{Z_i^{k}}^\top+ E^k) +\Lambda^k )\Big)_+
\end{eqnarray*}
 where $V_+ = U_\dag D_+ U_\dag$ such that
 $$
 UDU =\begin{pmatrix}
U_\dag & U_\ddag
 \end{pmatrix}\begin{pmatrix}
  D_+ & 0\\
 0  & D_-
 \end{pmatrix}\begin{pmatrix}
U_\dag \\
 U_\ddag
\end{pmatrix}, $$
is the eigen-decomposition of the matrix $V$ , and $D_+$ and $D_-$ are the nonnegative
and negative eigenvalues of $V$.
\end{enumerate}

\section{Convergence Analysis}\label{apend2}

Before establishing the main result on global convergence of the proposed ADMM algorithm, we provide the necessary definitions used in the proofs (for more details see \citet{Botle14}):

\begin{defn}(Kurdyka- Lojasiewicz property).\\
  The function $f$ is said to have the Kurdyka- Lojasiewicz (K-L) property at point $Z_0$, if there exist $c_1>0$,~$c_2>0$ and $\phi \in \Gamma_{c_2}$ such that for all $$Z \in B(Z_0,c_1) \cap \{Z: f(Z_0) <f(Z)<f(Z_0)+ c_2\},$$
  the following inequality holds
\begin{equation*}
 \phi'\big(f(Z)-f(Z_0)\big) \text{dist} \big(0,\partial f(Z)\big)\geq 1,
\end{equation*}
where $\Gamma_{c_2}$ stands for the class of functions $\phi: [0,c_2]\rightarrow \mathbb{R}^+$ with the properties:
\begin{enumerate}[(i)]
  \item $\phi$ is continuous on $[0,c_2)$;
  \item $\phi$ is smooth concave on $(0,c_2)$;
  \item $\phi(0)=0$, $\nabla \phi(s)> 0, ~\forall~s \in(0,c_2)$.
\end{enumerate}
\end{defn}
\begin{defn}
   (Semi-algebraic sets and functions).
   \begin{enumerate}[(i)]
     \item    A subset $C \in \mathbb{R}^{n \times n}$ is semi-algebraic, if there exists a finite number of real polynomial functions $h_{ij}$, $s_{ij}: \mathbb{R}^{n \times n} \rightarrow \mathbb{R}$ such that $$
     C=\cup_{i=1}^{\bar{p}}\cap_{j=1}^{\bar{q}}\{Z\in \mathbb{R}^{ n \times n} :~g_{ij} (Z)=0 \quad \text{and} \quad s_{ij}(Z)<0\}.
   $$
     \item A function $h: \mathbb{R}^{n \times n}\rightarrow (-\infty, +\infty]$ is called semi-algebraic, if its graph
         $$
         \mathbb{G}(h):= \{(Z,y)\in \mathbb{R}^{n\times n+1}: h(Z)=y\},
      $$
   \end{enumerate}
is a semi-algebraic set in $R^{n \times n+1}$.
\end{defn}

\begin{defn}
   (Sub-analytic sets and functions).
   \begin{enumerate}[(i)]
     \item    A subset $C\in \mathbb{R}^{n \times n}$ is sub-analytic, if there exists a finite number of real analytic functions $h_{ij}$, $s_{ij}: \mathbb{R}^{n \times n} \rightarrow \mathbb{R}$ such that $$
     C=\cup_{i=1}^{\bar{p}}\cap_{j=1}^{\bar{q}}\{ Z\in \mathbb{R}^d : g_{ij}(Z)=0 \quad \text{and} \quad s_{ij}(Z)<0\}.
   $$
     \item A function h: $R^{n \times n}\rightarrow (-\infty, +\infty]$ is called sub-analytic, if its graph
         $$
         \mathbb{G}(h):= \{(Z,y)\in \mathbb{R}^{n \times n+1}: h(Z)=y\}
      $$
   \end{enumerate}
is a sub-analytic set in $R^{n \times n+1}$.
\end{defn}

 It can be easily seen that both real analytic and semi-algebraic functions are sub-analytic. In general, the sum of two sub-analytic functions is not necessarily sub-analytic. However, it is easy to show that for two sub-analytic functions, if at least one function maps bounded sets to bounded sets, then their sum is also sub-analytic \citep{Botle14}.

\begin{rem}\label{rem22}
Each $f_i$ in \eqref{eq:3.1} is a convex semi-algebraic function (see, example~5.3 in \citep{Botle14}), while the loss function $\Gc$ in \eqref{eq:5}, \eqref{eq:13}, \eqref{eq:7}, \eqref{eq:10}, and \eqref{eqaut:7} is sub-analytic (even analytic).
Since each function $f_i$ maps bounded sets to bounded sets, we can conclude that the augmented Lagrangian function
\begin{eqnarray*}
   \Lc_{\gamma} (\Theta, Z_1, \dots, Z_n,E; \Lambda) &=& \Gc(X,\Theta)+ f_1(Z_1)+\dots +f_n(Z_n)+ f_e(E)\\
 &-& \langle \Lambda, \Theta-\sum_{i=1}^{n}Z_i+Z_i^\top -E \rangle\\
 &+& \frac{\gamma}{2} \| \Theta-\sum_{i=1}^{n}Z_i+Z_i^\top -E \|^2_{F},
\end{eqnarray*}
which is the summation of sub-analytic functions is itself sub-analytic. All sub-analytic functions which are continuous over their domain satisfy a K-L inequality, as well as some, but not all, convex functions (see \citealp{Botle14} for details and a counterexample).
Therefore,  the augmented Lagrangian function $\Lc_{\gamma}$ satisfies the K-L property.
\end{rem}

Next, we establish a series of lemmas used in the proof of Theorem \ref{thm1}.
\begin{lem}\label{lem1}
Let $U^k:=(\Theta^k, Z^k_1, \dots, Z^k_n,E^k; \Lambda^k)$ be a sequence generated by Algorithm~\ref{alg:1}, then there exists a positive constant $\vartheta$ such that
\small{
 \begin{eqnarray}\label{h5}
 \nonumber
 \Lc_{\gamma}(U^{k+1}) &\leq& \Lc_{\gamma} (U^k) - \frac{\vartheta}{2}\Big(\|\Theta^{k}-\Theta^{k+1}\|_{F}\\
 &+&\sum_{i=1}^{n}\|Z^k_i-Z^{k+1}_i\|_{F}
+\|E^k-E^{k+1}\|_{F}+\|\Lambda^k-\Lambda^{k+1}\|_{F} \Big).
\end{eqnarray}
}
\end{lem}
\textit{Proof.}
Using the first-order optimality conditions for \eqref{eq:3.3a0} and the convexity of $\Gc(X,\Theta)$, we obtain
\begin{eqnarray}\label{a1}
 \nonumber 
  0 &=& \big \langle \Theta^{k}-\Theta^{k+1},\nabla \Gc(X,\Theta^{k+1})-\Lambda^{k} + \gamma (\Theta^{k+1}-\sum_{i=1}^{n} Z^{k}_i + {Z^{k}_i}^\top-E^k)\big \rangle \\
  \nonumber
   &\leq& \Gc(X, \Theta^{k})-\Gc(X,\Theta^{k+1})-\langle \Theta^k-\Theta^{k+1} , \Lambda^k \rangle \\
   \nonumber
   &+& \gamma \langle \Theta^k-\Theta^{k+1}, \Theta^{k+1}-\sum_{i=1}^{n} Z^{k}_i + {Z^{k}_i}^\top -E^{k}\rangle\\
   \nonumber
   &=&\Gc(X,\Theta^{k})-\langle\Theta^k,\Lambda^k\rangle + \frac{\gamma}{2}\sum_{i=1}^{n}\|\Theta^k -\sum_{i=1}^{n} Z^{k}_i + {Z^{k}_i}^\top-E^{k}\|^2_{F} -\frac{\gamma}{2}\|\Theta^k-\Theta^{k+1}\|^2_{F}\\
   \nonumber
   &-&\Big(\Gc(X,\Theta^{k+1}) - \langle\Theta^{k+1},\Lambda^k\rangle + \frac{\gamma}{2}\|\Theta^{k+1}-\sum_{i=1}^{n} Z^{k}_i + {Z^{k}_i}^\top-E^{k}\|^2_{F} \Big)\\
   &=& \Lc_\gamma(U^k)-\Lc_\gamma(\Theta^{k+1}, Z^k_1,\dots, Z^{k}_n,E^k;\Lambda^k)
   -\frac{\gamma}{2}\|\Theta^k-\Theta^{k+1}\|^2_{F},
  \end{eqnarray}
where the second equality follows from the fact that
  $$(u_1-u_2)^T(u_3-u_1)=\frac{1}{2}\Big( \|u_2-u_3\|^2_{F}-\|u_1-u_2\|^2_{F}-\|u_1-u_3\|^2_{F}\Big).$$

Using \eqref{eq:3.3a1}, \eqref{eq:3.3a2} and  Lemma~\ref{lemdec}, we have that

\begin{eqnarray}\label{a22}
\nonumber
 \Lc_\gamma(\Theta^{k+1}, Z^{k}_{1}, Z^k_2,\dots,E^k;\Lambda^k)&-& \Lc_\gamma(\Theta^{k+1}, Z^{k+1}_{1}, Z^k_{2},\dots, E^k;\Lambda^k)\\
 \nonumber
 &-& \frac{(\gamma \varrho- L_{H_1})}{2}\|Z^k_1-Z^{k+1}_1\|^2_{F}\\
 \nonumber
 &\geq& 0,\\
 \nonumber
 \Lc_\gamma(\Theta^{k+1},\dots, Z^{k+1}_{i-1}, Z^k_i,\dots,E^k;\Lambda^k)&-& \Lc_\gamma(\Theta^{k+1},\dots, Z^{k+1}_{i}, Z^k_{i+1},\dots, E^k;\Lambda^k)\\
 \nonumber
 &-& \frac{(\gamma \varrho- L_{H_i})}{2}\|Z^k_i-Z^{k+1}_i\|^2_{F}\\
 &\geq& 0,  \qquad\qquad i=2, \dots, n,
\end{eqnarray}
where $L_{H_i}$ is a Lipschitz constant of the gradient $\nabla H_i(Z_i)$, and $\varrho\geq \frac{L_{H_i}}{\gamma}, ~(i=1,\dots,n)$ is a proximal parameter.

Following the same steps as \eqref{a1}, we have that

\begin{eqnarray}\label{a2}
\nonumber
 \Lc_\gamma(\Theta^k,Z^{k+1}_1,\dots,Z^{k+1}_n,E^k;\Lambda^k)&-& \Lc_\gamma(\Theta^{k+1}, Z^{k+1}_1,\dots,Z^{k+1}_n,E^{k+1};\Lambda^k)\\
 \nonumber
&-& \frac{\gamma}{2}\|E^k-E^{k+1}\|^2_{F}\\
&\geq&0,
\end{eqnarray}
and
\begin{eqnarray}\label{a222}
\nonumber
\Lc_\gamma(\Theta^{k+1},Z^{k+1}_1,\dots,Z^{k+1}_n,E^{k+1};\Lambda^k)&-& \Lc_\gamma(\Theta^{k+1},Z^{k+1}_1,\dots,Z^{k+1}_n,E^{k+1};\Lambda^{k+1})\\
\nonumber
&-& \frac{\lambda^2_e}{\gamma}\|E^k-E^{k+1}\|^2_{F}\\
&\geq&0.
\end{eqnarray}

Let $$\hat{\gamma} :=  \max(\gamma\varrho- L_{H_1}, \dots,  \gamma\varrho- L_{H_n}), \quad  \bar{\gamma}:=\frac{\gamma^2-2\lambda^2_e}{\gamma(1+\lambda^2_e)}, \quad  \vartheta :=\max(\hat{\gamma},\bar{\gamma},\gamma).$$

Then, using (\ref{a1})-- (\ref{a222}), and $\gamma \geq \sqrt{2} \lambda_e$, we have

\begin{eqnarray*}\label{a3}
 \nonumber 
    &&\Lc_\gamma(U^k)- \Lc_\gamma(U^{k+1}) \geq\frac{\gamma}{2}\|\Theta^k-\Theta^{k+1}\|^2_{F}\\
    \nonumber
    &+& \frac{\hat{\gamma}}{2}\sum_{i=1}^{n}\|Z^k_i-Z^{k+1}_i\|^2_{F}+\frac{\gamma^2-2\lambda^2_e}{2\gamma}\|E^k-E^{k+1}\|^2_{F},\\
    \nonumber
   &=&\frac{\gamma}{2}\|\Theta^k-\Theta^{k+1}\|^2_{F}+ \frac{\hat{\gamma}}{2} \sum_{i=1}^{n}\|Z^k_i-Z^{k+1}_i\|^2_{F}+ \frac{\bar{\gamma}}{2}\|E^k-E^{k+1}\|^2_{F}+ \frac{\lambda^2_e \bar{\gamma}}{2}\|E^k-E^{k+1}\|^2_{F},\\
    \nonumber&=&\frac{\gamma}{2}\|\Theta^k-\Theta^{k+1}\|^2_{F}+ \frac{\hat{\gamma}}{2}\sum_{i=1}^{n}\|Z^k_i-Z^{k+1}_i\|^2_{F}+ \frac{\bar{\gamma}}{2}\Big(\|E^k-E^{k+1}\|^2_{F}+\|\Lambda^k-\Lambda^{k+1}\|^2_{F}\Big),\\ &\geq&\frac{\vartheta}{2}\Big(\|\Theta^k-\Theta^{k+1}\|^2_{F}+\sum_{i=1}^{n}\|Z^k_i-Z^{k+1}_i\|^2_{F}+\|E^k-E^{k+1}\|^2_{F}+\|\Lambda^k-\Lambda^{k+1}\|^2_{F}\Big).
\end{eqnarray*}
\hfill$\Box$

\begin{lem}\label{lem2}
Let $ U^k=(\Theta^k, Z^k_1,\dots X^k_n,E^k,\Lambda^k)$ be a sequence generated by Algorithm~\ref{alg:1}. Then, there exists a subsequence $U^{k_s}$ of $\{U^k\}$, such that

\begin{equation*}
\lim_{s\rightarrow\infty} \Gc(X, \Theta^{k_s})= g(\Theta^*), \quad \lim_{s\rightarrow\infty} f_i(Z^{k_s}_i)= f_i(Z^*_i), \quad \lim_{s\rightarrow\infty} f_e(E^{k_s}_i) = f_e(E^*_i),
\end{equation*}
where
$$
\lim_{s\rightarrow\infty} U^{k_s}=(\Theta^*, Z^*_1,\dots, Z^*_n,E^*,\Lambda^*).
$$
\end{lem}

\textit{Proof}. Let $\Upsilon^{k+1} = \Theta^{k+1}-\sum_{i=1}^{n}Z^{k+1}_i + {Z^{k+1}_i}^\top -E^{k+1}$. Using the quadratic function $f_e(E)=\frac{\lambda_e}{2}\|E\|^2_{F}$, we have that
\begin{eqnarray}\label{a505}
\nonumber
f_e(E^{k+1}-\Upsilon^{k+1})&=&\frac{\lambda_e}{2}\|E^{k+1}-\Upsilon^{k+1}\|^2_{F}\\
&=& \frac{\lambda_e}{2}\|E^{k+1}\|^2- \lambda_e\langle
E^{k+1}, \Upsilon^{k+1}\rangle + \frac{\lambda_e}{2}\|\Upsilon^{k+1}\|^2_{F}.
\end{eqnarray}

Using \eqref{a505} and the fact that each function $f_i$ is lower bounded, there exists $\underline{\Lc}$, such that

\begin{eqnarray} \label{a33}
\nonumber
\Lc_{\gamma} (U^{k+1}) &=& \Gc(X,\Theta^{k+1})+f_1(Z^{k+1}_1)+\dots f_n(Z^{k+1}_n)+ \frac{\lambda_e}{2}\|E^{k+1}-\Upsilon^{k+1}\|^2_{F}\\
&+& \frac{\gamma-\lambda_e}{2} \| \Upsilon^{k+1}\|^2_{F} \geq  \underline{g} + \underline{f_1}+ \dots+ \underline{f_n}\geq\underline{\Lc},
\end{eqnarray}
since $\Gc(X,\Theta^{k+1})$ and $f_i(Z^{k+1}_i)(i=1,\dots, n)$ are all lower bounded.

Now, using Lemma~\ref{lem1}, we have that

\begin{eqnarray}\label{a34}
 \nonumber 
\frac{\vartheta}{2}\sum_{k=0}^{K}\Big(\|\Theta^k-\Theta^{k+1}\|^2_{F}+\sum_{i=1}^{n}\|Z^k_i-Z^{k+1}_i\|^2_{F}&+&\|E^k-E^{k+1}\|^2_{F}+\|\Lambda^k-\Lambda^{k+1}\|^2_{F}\Big) \\ &\leq&\Lc_\gamma(U^0)- \underline{\Lc}.
 \end{eqnarray}

Lemma~\ref{lem1} together with \eqref{a34} shows that $\Lc_{\gamma}(U^k)$ converges to $\Lc_{\gamma}(U^*)$. Note that \eqref{a34} and the coerciveness of $\Gc(X,\Theta)$ and $f_i~(i=1,\dots, n)$ imply that $\{(\Theta^k, Z^k_1,\dots,Z^k_n)\}$ is a bounded sequence.
This together with the updating formula of $\Lambda^{k+1}$ and \eqref{a34} yield the boundedness of $E^{k+1}$.  Moreover, the fact that $\Lambda^k = -\lambda_e E^k$, gives the boundedness of $\Lambda^{k}$, which implies that the
entire sequence $\{U^k\}$ is a bounded one. Therefore, there exists a subsequence $$ U^{k_s} = (\Theta^{k_s}, Z^{k_s}_1, \dots, Z^{k_s}_n,E^{k_s}; \Lambda^{k_s}), \quad \quad s=0,1,\dots$$ such that $U^{k_s} \rightarrow U^*$ as $s \rightarrow \infty$.

Now, using the fact that $\Gc(X,\Theta)$, $f_i(Z_i)~(i=1,\dots, n)$ and $f_e(E)$ are continuous functions, we have that
\begin{equation*}
\lim_{s\rightarrow\infty} \Gc(X,\Theta^{k_s})= g(\Theta^*), \quad \lim_{s\rightarrow\infty} f_i(Z^{k_q}_i)= f_i(Z^*_i), \quad \lim_{s\rightarrow\infty} f_e(E^{k_q}_i) = f_e(E^*_i).
\end{equation*}
\hfill$\Box$

\begin{lem}\label{lem4}
Algorithm~\ref{alg:1} either stops at a stationary point of the problem \eqref{eq:3.1} or generates an infinite sequence $\{U^k\}$, so that any limit point of $\{U^k\}$ is a critical point of  $\Lc_{\gamma}(U^k)$ \eqref{eq:3.1}.
\end{lem}

\textit{Proof}. From the definition of the augmented Lagrangian function in \eqref{eq:3.1}, we have that
\begin{eqnarray} \label{a55}
\nonumber
&&\nabla \Gc(X,\Theta^{k+1})-\Lambda^{k+1} +\gamma\Upsilon^{k+1}  = \nabla_{\Theta} \Lc_\gamma(U^{k+1}), \\
\nonumber
&&\partial f_i(Z^{k+1}_i) - \Lambda^{k+1}-{\Lambda^{k+1}}^\top - \gamma(\Upsilon^{k+1}+ {\Upsilon^{k+1}}^\top) \in \partial_{Z_i} \Lc_\gamma(U^{k+1}),\quad i=1,\dots, n,\\
\nonumber
&& \lambda_e E^{k+1}+ \Lambda^{k+1} - \gamma \Upsilon^{k+1} = \nabla_{E} \Lc_\gamma(U^{k+1}),\\
&& \gamma \Upsilon^{k+1}=- \nabla_{\Lambda} \Lc_\gamma(U^{k+1}),
\end{eqnarray}
where $\Upsilon^{k+1} = \Theta^{k+1}-\sum_{i=1}^{n}Z^{k+1}_i + {Z^{k+1}_i}^\top -E^{k+1}$.

Moreover, the updating formula of $\Lambda^{k+1}$, \eqref{eq:3.2} and \eqref{eq:3.3lin} yields that

\begin{eqnarray} \label{a56}
\nonumber
\nabla \Gc(X,\Theta^{k+1})- \Lambda^{k+1} &=& \gamma \Big(\Theta^{k+1} - \Theta^{k} \\
\nonumber
&+& \sum_{i=1}^{n} Z^{k}_i- Z^{k+1}_i + (Z^{k}_i-Z^{k+1}_i)^\top + E^{k}-E^{k+1}\Big)\\
\nonumber
\partial f_1(Z^{k+1}_1)- \Lambda^{k+1}- {\Lambda^{k+1}}^\top &=& \gamma\varrho(Z^{k}_1- Z^{k+1}_1)
+ \gamma \Big(\Theta^{k+1} - \Theta^{k}\\
&+& (\Theta^{k+1} - \Theta^{k})^\top + \sum_{i=1}^{n} Z^{k}_i-Z^{k+1}_i \\
\nonumber
& +& (Z^{k}_i-Z^{k+1}_i)^\top
+ E^{k}-E^{k+1}+ (E^{k}-E^{k+1})^\top\Big)\\
\nonumber
\partial f_i(Z^{k+1}_i)- \Lambda^{k+1}- {\Lambda^{k+1}}^\top &=& \gamma\varrho(Z^{k}_i- Z^{k+1}_i)\\
\nonumber
&+& \gamma\Big( \Theta^{k+1} - \Theta^{k}  + (\Theta^{k+1} - \Theta^{k})^\top\\
\nonumber
&+& \sum_{j=i}^{n} Z^{k}_i-Z^{k+1}_i + (Z^{k}_i-Z^{k+1}_i)^\top\\
\nonumber
&+& E^{k}-E^{k+1}+ (E^{k}-E^{k+1})^\top\Big) \qquad i=2,\dots, n,\\
\lambda_e E^{k+1}&+&\Lambda^{k+1} = 0.
\end{eqnarray}

Combining \eqref{a55}, \eqref{a56}, and the updating formula of $\Lambda^{k+1}$, we have that

\begin{equation}\label{a59}
(\hbar^{k+1}_\Theta, \hbar^{k+1}_1, \dots, \hbar^{k+1}_n, \hbar^{k+1}_E, \hbar^{k+1}_\Lambda)\in \partial \Lc_{\gamma}(U^{k+1}),
\end{equation}
where
\begin{eqnarray} \label{a57}
\nonumber
\hbar^{k+1}_{\Theta}&:=&   \Lambda^{k}- \Lambda^{k+1} + \gamma \Big(\Theta^{k+1} - \Theta^{k} + \sum_{i=1}^{n} Z^{k}_i- Z^{k+1}_i + (Z^{k}_i-Z^{k+1}_i)^\top + E^{k}-E^{k+1} \Big)\\
\nonumber
\hbar^{k+1}_{Z_1} &:=&  \Lambda^{k}- \Lambda^{k+1}  +  (\Lambda^{k}- {\Lambda^{k+1}})^\top + \gamma\varrho(Z^{k}_1- Z^{k+1}_1) \\
\nonumber
&+& \gamma \Big(\Theta^{k+1} - \Theta^{k} + (\Theta^{k+1} - \Theta^{k})^\top + \sum_{i=1}^{n} Z^{k}_i-Z^{k+1}_i + (Z^{k}_i-Z^{k+1}_i)^\top \\
\nonumber
&+& E^{k}- E^{k+1}+ (E^{k}-E^{k+1})^\top  \Big) \\
\nonumber
\hbar^{k+1}_{Z_i}&:=&  \Lambda^{k}- \Lambda^{k+1}  +  (\Lambda^{k}- {\Lambda^{k+1}})^\top+ \gamma\varrho(Z^{k}_i- Z^{k+1}_i)\\
\nonumber
&+& \gamma\Big( \Theta^{k+1} - \Theta^{k}  + (\Theta^{k+1} - \Theta^{k})^\top+ \sum_{j=i}^{n} Z^{k}_i-Z^{k+1}_i + (Z^{k}_i-Z^{k+1}_i)^\top\\
\nonumber
&+& E^{k}-E^{k+1}+ (E^{k}-E^{k+1})^\top\Big), \qquad i=2,\dots, n,\\
\nonumber
\hbar^{k+1}_E &:=& \Lambda^k-\Lambda^{k+1},\\
\hbar^{k+1}_\Lambda &:=& \frac{1}{\gamma}(\Lambda^{k+1}-\Lambda^k),
\end{eqnarray}

Now, using \eqref{a34}, we obtain that
\begin{eqnarray} \label{a66}
\lim_{k\rightarrow\infty}(\|\hbar^{k+1}_\Theta\|_{F},\|\hbar^{k+1}_{Z_1}\|_{F},\dots, \|\hbar^{k+1}_{Z_n}\|_{F}, \|\hbar^{k+1}_E\|_{F}; \|R^{k+1}_{\Lambda}\|_{F})=(0,\dots,0).
\end{eqnarray}

Suppose that Algorithm~\ref{alg:1} does not stop at a stationary point. Using Lemma~\ref{lem2}, there exists a subsequence ${U^{k_s}}$, such that $U^{k_s}\rightarrow U^*$ as $s\rightarrow \infty.$  Using \eqref{a59} and \eqref{a66}, we conclude that $(0,\dots,0)\in \partial \Lc_{\gamma}(U^*)$.
\hfill$\Box$

\textbf{Proof of Theorem~\ref{thm1}.}
Lemmas~\ref{lem2} and ~\ref{lem4} imply that $\{U^k\}$ is a bounded sequence and the set of limit points of $\{U^k\}$ starting from $U^0$ is non-empty, respectively.
Moreover, Lemma~5 and Remark~5 of \citep{Botle14} imply that the set of limit points of $\{U^k\}$ starting from $U^0$ is compact. The remainder of the proof of this Theorem
follows along similar lines to the proof of Theorem~1 in \citep{Botle14}, by utilizing the K-L property of the problem \eqref{eq:3.1} (see, Remark~\ref{rem22}).
\hfill$
\Box$

\vskip 0.2in
\bibliography{16-486-ref}

\end{document}